\documentclass[11pt,reqno]{amsart}
\pdfoutput=1
\numberwithin{equation}{section} \oddsidemargin=-.0cm
\usepackage{newtxtext}
\usepackage[scaled=0.90]{helvet} 

\usepackage{bm}
\usepackage{ stmaryrd }

\usepackage{dsfont}
\usepackage{subcaption}

\usepackage{amssymb}
\usepackage{amsmath}
\usepackage{color}
\usepackage{graphicx, graphics}
\usepackage{amsfonts, amssymb, amsthm}
\usepackage{mathrsfs}

\usepackage{sidecap}
\usepackage{needspace} 
\usepackage{cases} 
\sidecaptionvpos{figure}{c}

\usepackage[colorlinks=true, pdfstartview=FitV, linkcolor=blue,
            citecolor=blue, urlcolor=blue]{hyperref}  

\usepackage[sort]{cite}
  \usepackage{paralist}
  \usepackage{graphics} 
  \usepackage{epsfig} 
\usepackage{graphicx}  
\usepackage{epstopdf}

\usepackage{epstopdf}
\usepackage{bbm} 
\usepackage{amsfonts, amssymb, amsthm}
\usepackage{color}
\usepackage{needspace} 

\usepackage{booktabs} 

 \usepackage[colorlinks=true]{hyperref}
\hypersetup{urlcolor=blue, citecolor=blue}

\numberwithin{equation}{section}

\newcommand{\ep}{\varepsilon}

\newcommand{\betrag}[1]{\left\lvert #1 \right\rvert}	        

\newtheorem{thm}{Theorem}[section]
\newtheorem{lem}{Lemma}[section]
\newtheorem{defi}{Definition}[section]
\newtheorem{ex}{Example}[section]
\newtheorem{prop}[thm]{Proposition}

\newtheorem{rem}{Remark}[section]
\newtheorem{cor}{Corollary}[section]
\newenvironment{Proof}[1][Proof]{\begin{trivlist} \item[\hskip \labelsep {\it \bf  #1.}]}{\end{trivlist}}
\def\bpp{\begin{Proof}}
\def\epp{\end{Proof}}

\def\bt{\begin{thm}}
\def\et{\end{thm}}

\def\bex{\begin{ex}}
\def\eex{\end{ex}}

\def\bl{\begin{lem}}
\def\el{\end{lem}}
\def\bd{\begin{defi}}
\def\ed{\end{defi}}
\def\bc{\begin{cor}}
\def\ec{\end{cor}}
\def\bp{\begin{proof}}
\def\ep{\end{proof}}
\def\br{\begin{rem}}
\def\er{\end{rem}}

\def\bprop{\begin{prop}}
\def\eprop{\end{prop}}

\def\d{\mathrm{d}}
\def\Vert{\: \vert \:}
\def\be{\begin{equation}}
\def\ee{\end{equation}}
\def\bes{\begin{equation*}}
\def\ees{\end{equation*}}
\def\bea{\begin{equation} \begin{aligned}}
\def\eea{\end{aligned} \end{equation}}
\def\beas{\begin{equation*} \begin{aligned}}
\def\eeas{\end{aligned} \end{equation*}}

\def\bi{\begin{itemize}}
\def\ei{\end{itemize}}
\def\ben{\begin{enumerate}}
\def\een{\end{enumerate}}

\def\s{\mathfrak{s}}
\def\c{\mathfrak{c}}

\def\span{\mathrm{span}}

\definecolor{rred}{rgb}{0.7,0,0.1}
\definecolor{ccyan}{rgb}{0,.5,1}
\definecolor{greenrb}{rgb}{0.2,0.6,0.2}

\newcommand{\mkr}{\color{black}}

\usepackage[normalem]{ulem}

\newcommand{\m}{\ensuremath{\mathfrak{m}}}
\newcommand{\n}{\ensuremath{\mathfrak{n}}}

\newcommand{\R}{\ensuremath{\mathbb{R}}}

\newcommand{\N}{\ensuremath{\mathbb{N}}}

\newcommand{\W}{\ensuremath{\mathbb{P}}}
\newcommand{\E}{\ensuremath{\mathbb{E}}}

\title[A Girsanov approach to slow parameterizing manifolds in the presence of noise]{A Girsanov approach to slow parameterizing manifolds in the presence of noise}

\author[M. D. Chekroun]{Micka\"el D. Chekroun}
\address[MDC]{Department of Atmospheric \& Oceanic Sciences, University of California, Los Angeles, USA, and Department of Earth and Planetary Sciences, Weizmann Institute, Rehovot, Israel} 
\email{mchekroun@atmos.ucla.edu}

\author[J. S. W. Lamb]{Jeroen S.W. Lamb}
\address[JSWL]{Department of Mathematics, Imperial College, London, UK} 

\author[C. J. Pangerl]{Christian J. Pangerl}
\address[CJP]{Department of Mathematics, Imperial College, London, UK} 

\author[M. Rasmussen]{Martin Rasmussen}
\address[MR]{Department of Mathematics, Imperial College, London, UK} 

\thanks{Corresponding author: Micka\"el D. Chekroun}


\begin{document}

\begin{abstract}

This work investigates a three-dimensional slow-fast stochastic system with quadratic nonlinearity and additive noise, inspired by fluid dynamics. 
The deterministic counterpart exhibits a periodic orbit and a slow manifold. 
We demonstrate that, under specific parameter regimes, this deterministic slow manifold can serve as an approximate parameterization of the fast variable by the slow variables within the stochastic system.

Building upon this parameterization, we derive a two-dimensional reduced model, a stochastic Hopf normal form, that captures the essential dynamics of the original system. Both the original and the reduced systems possess ergodic invariant measures, characterizing their long-term behavior. 

We quantify the discrepancy between the original system and its slow approximation by deriving error estimates involving the Wasserstein distance between the marginals of these invariant measures along the radial component. These error bounds are shown to be controlled by a parameterization defect, which measures the quality of the fast-slow variable parameterization.

A key technical innovation lies in the application of Girsanov's theorem to obtain these error estimates in the presence of oscillatory instabilities.  Furthermore, we extend our analysis to regimes exhibiting an "inverted" timescale separation, where the variable to be parameterized evolves on a slower timescale than the resolved variables. To address these more challenging scenarios, we introduce path-dependent coefficients in the parameterizing manifold, enabling the derivation of robust error bounds for the corresponding reduced model. Numerical simulations complement our theoretical findings, providing insights into the model's behavior and exploring parameter regimes beyond the scope of our analytical results.

\end{abstract}

\maketitle 

\tableofcontents
\section{Introduction and Motivations}
We consider the following system
\begin{subnumcases}{\label{Eq_3DHopf}}
\d x =(\lambda x- f y -\gamma xz) \d t + \sigma \d W^1_t\\
\d y =(f x+\lambda y -\gamma yz) \d t+ \sigma \d W^2_t\\
\d z =-\frac{1}{\epsilon}(z-x^2-y^2)\d t+ \frac{\sigma}{\sqrt{\epsilon}} \d W^3_t \label{Eq_3DHopf_z}.
\end{subnumcases}

The stochastic processes $W_t^1, W_t^2, W_t^3$ are independent Brownian motions on the probability space $(\Omega, \mathcal{B}, \mathbb{P})$, where $\Omega=C_0([0,\infty),\mathbb{R})$ denotes the space of continuous functions $\omega$ defined on the non-negative real line,  $\mathcal{B}$ denotes the Borel $\sigma$-algebra on $\Omega$,  and $\mathbb{P}$ denotes the Wiener measure on $\Omega$ \cite{arnold2013random}. The parameters $\lambda, f, \gamma$ and $\epsilon$ are assumed to be positive throughout this article.

This model is known for $\sigma=0$ to provide a low-dimensional reduced system  of a nonlinear fluid flow past a circular cylinder which is characterized by vortex shedding for moderate Reynolds number \cite{noack_al_2003,lusch2018deep}. For such regimes, system \eqref{Eq_3DHopf} is a mean-field model that exhibits a stable limit cycle corresponding to von Karman vortex shedding, and an unstable equilibrium corresponding to a low-drag condition.  The noise terms  can be interpreted for e.g.~accounting for turbulent perturbations in a Lagrangian setting \cite{rodean1996stochastic}.

 Mathematically, system \eqref{Eq_3DHopf} is a slow-fast system driven by additive noise.
In the deterministic limit ($\sigma = 0$), the theory of slow-fast systems is well-established, particularly when a strong timescale separation exists between the "slow" (resolved) and "fast" (unresolved) variables, as controlled by the parameter $\epsilon$; see \cite{Jones1995,nipp2013invariant, wiggins2013normally,Kuehn:2015tol}. This separation is evident as $\epsilon$ diminishes, leading to rapid fluctuations in the $z$-variable while the $x$- and $y$-variables evolve more slowly.

A cornerstone of deterministic slow-fast dynamics is the existence of a slow manifold that often describes an approximate "slaving" relationship between the slow and fast variables. This manifold, typically represented as the graph of a function $h$ of the slow variables, attracts the system's trajectories within an "$\epsilon$-neighborhood" for sufficiently small $\epsilon$. For system \eqref{Eq_3DHopf} with $\sigma=0$, the slow manifold is explicitly given by:
\bea \label{Eq_slow_mnf_int}
h: \, & \mathbb{R}\times \mathbb{R} \longrightarrow \mathbb{R},\\
&(x,y) \mapsto x^2 + y^2.
\eea
More precisely, any solution $(x(t),y(t),z(t))$ to Eq.~\eqref{Eq_3DHopf} is attracted exponentially fast to the manifold given by $\mathcal{M}_\epsilon=\mbox{graph} (h) + \mathcal{O}(\epsilon)$, and in particular  the near to slaving relationship $z(t)=h(x^2(t)+y^2(t)) + \mathcal{O}(\epsilon)$ holds for $t$ sufficiently large. 
 This dynamical behavior is grounded in the foundational work of Fenichel \cite{FENICHEL197953,Fenichel_Moser71} and Tikhonov \cite{tikhonov1952systems}, which established these principles for a broader class of systems beyond Eq.~\eqref{Eq_3DHopf}. For a comprehensive introduction to invariant manifolds (IMs) in singularly perturbed ordinary differential equations (ODEs), refer to \cite[Chap.~10]{nipp2013invariant}and \cite{Kuehn:2015tol}.

From a physical perspective, a crucial objective is to identify the "reduced" variables that effectively capture the long-term behavior of the full system. In the presence of pronounced timescale separation, the slow manifold provides a compelling solution to this problem. For deterministic systems with sufficiently small $\epsilon$, the accuracy of reductions based on the slow manifold is well-justified, as the manifold closely approximates the underlying IM. Various methods have been developed for computing such manifolds in deterministic slow-fast systems, including those by \cite{Guckenheimer_Kuehn09}  and \cite{England_Krauskopf_Osinga07}.

 However, in the absence of time-scale separation when e.g.~the slow and fast motions co-exist for each variable, the concept of slow manifolds has to be revisited. The determination of reduced equations is no longer tied to invariant or nearly invariant manifolds such as slow manifolds but to other manifolds that still capture a slow motion albeit only in an averaged sense \cite{chekroun2017emergence,CLM19_closure}. Slow manifolds have e.g.~a long history in atmospheric dynamics \cite{baer1977complete,leith1980nonlinear,lorenz1986existence}, and the study of their existence and role in the corresponding physical models has been a topic of active research and debate \cite{lorenz1992slow}. Breakdown of slaving as manifested by burst of energetic fast oscillations is indeed observed in such models \cite{vautard1986invariant,chekroun2017emergence} for certain parameter regimes, and recently the notion of  Optimal Parameterizing Manifolds (OPMs) has been introduced to cope with such a mixed dynamical behavior; see \cite{chekroun2017emergence,CLM19_closure,chekroun2021stochastic}.  
The OPM provides the manifold on which lies the average motion of the fast variables to parameterize as a function of the resolved variables, it averages out the fast fluctuations in an optimal  least-squares sense; see \cite[Sec.~3.2]{CLM19_closure} and \cite{chekroun2023optimal}.

In the stochastic context, the concept of the slow manifold has been generalized in various ways.
 One approach consists of using the deterministic slow manifold for reducing the corresponding stochastic system and to analyze its a posteriori  performance.
Such an approach is discussed in \cite{BERGLUND20031}, where finite-time estimates relative to the associated deterministic system are  derived.  A review of  reduced order modeling techniques for deterministic as well as stochastic slow-fast systems is provided in \cite{0951-7715-17-6-R01}. 

An alternative approach is grounded in the theory of random dynamical systems (RDSs)  \cite{arnold2013random}. Within the RDS framework, the fast variables, now stochastic, are approximated by a manifold that is still a function of the slow variables but now depends on the specific noise realization, thus rendering the manifold stochastic. This approach leads, for stochastic slow-fast systems, to a notion of stochastic slow manifold that is tied to the theory of random IMs \cite{Schmalfuss2008,WANG2013822}.  As such, these  stochastic slow manifolds are subject to spectral gap conditions encountered in the theory of random IMs \cite{CLW15_vol1} and that express, for slow-fast systems, a notion of time-scale separation; see \cite{debussche1991inertial} for the deterministic case. We refer also to  \cite{MR2808288,MR2501313} for computational and theoretical aspects in the related context of nonautonomous systems and to \cite{MR2229984} (resp.~\cite{CLW15_vol1}) for explicit Taylor approximations of non-autonomous (resp.~random) IMs.

However, situations where the slaving relationship between slow and fast variables breaks down have been less studied in stochastic systems. To address such issues, the concept of Stochastic Parameterizing Manifolds (SPMs) was introduced in \cite{CLW15_vol2}. An SPM aims to provide a random manifold that, while not necessarily invariant, approximates the variables to be parameterized (unresolved variables, 
$u_\s$) in terms of the resolved variables ($u_\c$) in a least-squares sense for almost every noise realization. This concept aligns with the notion of "fuzzy slow manifolds" \cite{warn1986nonlinear,warn1997nonlinear}, where only approximate slaving relationships are expected to hold, a concept shown to be more suitable for systems like planetary flows \cite{vanneste2004exponentially,temam2011slow}.

The authors in \cite{CLW15_vol2} demonstrated through numerical examples that Taylor formulas, traditionally employed for the rigorous approximation of random Invariant Manifolds (IMs) near criticality (such as random center manifolds \cite[Chaps.~6 and 7]{CLW15_vol1}), can still be effectively used to obtain approximate parameterizations, 
$h_{\textrm{pm}}$, of the unresolved variables, $u_\s$, in terms of the resolved variables, $u_\c$, even away from the instability onset; see also \cite{chekroun2023optimal}. Notably, \cite{CLW15_vol2} showed that the quality of such parameterizations is governed by a parameterization defect of the form $\overline{|u_\s(t)-h_{\textrm{pm}}(u_\c(t),\omega(t))|^2}$, where $\overline{(\cdot)}$ denotes time-averaging over some characteristic time of the dynamics, and $\omega(t)$ signifies the dependence on the noise realization.

Crucially, compared to previous works on random IMs, the SPM approach introduced in \cite{CLW15_vol2}, building upon the approximation formulas of \cite{CLW15_vol1}, allows for the explicit characterization of the dependence on $\omega(t)$. This dependence reflects the influence of the "past" of the driving noise, rendering the corresponding reduced systems path-dependent and non-Markovian; see \cite[Chap.~5]{CLW15_vol2}.

This article aims to rigorously analyze the SPM approach in the context of the stochastic system \eqref{Eq_3DHopf}, with a particular focus on deriving error estimates for the approximation of long-term statistics.
 In a first step, the reduction of Eq.~\eqref{Eq_3DHopf} using  the deterministic slow manifold \eqref{Eq_slow_mnf_int} is analyzed near criticality, for small noise and  $\epsilon \ll 1$. The case $\epsilon >1$ corresponding to an absence of  time-scale separation, is then analyzed in a second step. In this latter case, the deterministic slow manifold is no-longer a valid parameterization. Rather, we show that path-dependent coefficients depending on the past of the noise need to be included in the substitute parameterizations to derive efficient reduced systems; see Section \ref{Sec_stoch_pm} below.

To compare the long-term statistics of the original and reduced systems, we focus on their invariant measures. For analytical tractability, the analysis is conducted in polar coordinates.
We establish the existence of ergodic invariant measures, denoted by $\mu$ and $\nu$, for the original and reduced systems, respectively (see Appendices \ref{App_ex_uniq_org} and \ref{Appendix_E}). To quantify the discrepancy between these probability measures, we employ a metric $D_{\mathcal{F}}$ that is equivalent to the $L^1$-Wasserstein metric on an appropriate space of probability measures on the real line (see \eqref{Eq_IPM}). Specifically, denoting the marginals of $\mu$ and $\nu$ along the radial direction as $\mu_r$ and $\nu_r$, respectively, Theorems \ref{Thm_hopf_1} and \ref{Thm_hopf_2} provide conditions under which the following error estimate holds:
\begin{equation} \label{Eq_est_int}
D_{\mathcal{F}}(\mu_r, \nu_r) \leq     C(\lambda, \gamma, \sigma) +c \left(\int_{\mathbb{R}+\times \mathbb{R}} |z - r^2|^{4}  , \mu{r,z}(\d r,\d z)\right)^{\frac{1}{4}}.
\end{equation}
Here, $r^2=x^2+y^2$ represents the analytical expression of the slow manifold \eqref{Eq_slow_mnf_int} in polar coordinates, and $\mu_{r,z}$ denotes the marginal of $\mu$ onto the $(r,z)$-plane. The additive constant $C(\lambda, \gamma, \sigma)$ in Theorem \ref{Thm_hopf_1} depends explicitly on the invariant measures $\mu$ and $\nu$, whereas in Theorem \ref{Thm_hopf_2} it depends solely on $\nu$. Remark \ref{Rmk_epsilon_tails}-(iii) below elucidates that for parameter regimes near criticality, with sufficiently small noise and $\epsilon$, and a well-balanced nonlinear dissipation, $C(\lambda, \gamma, \sigma)$ remains small.

The integral term in the right-hand side (RHS) of \eqref{Eq_est_int} is the key controlling factor. This term quantifies the parameterization defect---the error incurred by approximating the variable $z$ using the slow manifold.
Therefore, inequality \eqref{Eq_est_int} establishes that, for specific parameter regimes, the parameterization defect directly controls the accuracy of the reduced system in approximating the statistics of the original system.  
Previous studies, including numerical and analytical results in \cite[Corollary 4.1 and Proposition 5.1]{CLW15_vol2}, have suggested that the magnitude of such parameterization defects plays a critical role in the efficacy of reduced-order models for stochastic partial differential equations (SPDEs). Furthermore, \cite[Theorem 1]{Chekroun2015}  provides error estimates involving the parameterization defect for the optimal control of PDEs using reduced systems. Our analysis provides a rigorous foundation for these observations.

A crucial step in deriving the error estimate \eqref{Eq_est_int} involves obtaining an auxiliary finite-time error estimate, as detailed in \eqref{Eq_aux_expl} below. The primary challenge in proving \eqref{Eq_aux_expl} stems from our working assumption: the existence of a limit cycle in the deterministic system ($\sigma = 0$), resulting from an oscillatory instability with a positive growth rate ($\lambda>0$) that is subsequently saturated by nonlinear effects.\footnote{It is noteworthy that when $\sigma>0$, the slow manifold reduction of system \eqref{Eq_3DHopf} formally corresponds to the Hopf normal form subject to additive noise. This fundamental system has been extensively studied in the literature. For instance, \cite{doan2018hopf} investigates it from a random dynamical systems perspective, while  \cite{tantet2020ruelle}  employs methods from the theory of Markov semigroups. Furthermore, \cite{engel2019bifurcation} provides a bifurcation analysis of systems with stochastically perturbed limit cycles.} For scenarios involving a single eigenvalue losing stability (pitchfork bifurcation), pathwise error estimates can be derived with high probability for a wide range of stochastic differential equations (SDEs), including SPDEs relevant to fluid dynamics, by leveraging the stochastic invariant manifold reduction formulas derived in \cite{CLW15_vol1} (see \cite[Theorem 4.1]{Chekroun_al2023}).

To address the presence of oscillatory instabilities, potentially away from the instability onset, we employ a novel approach based on probability law transformations. First, we modify the drift term of the original system \eqref{Eq_3DHopf} in polar coordinates by introducing a coupling term that depends on the reduced system (see Remark \ref{Connection_Hairer} below). This coupling term is designed to counteract the linear unstable effects. While the transition probabilities of this modified system generally differ from those of the original system, we demonstrate that a Girsanov transformation exists (see Lemma \ref{Lemma_girsanov}). This allows us to identify a new probability measure $\tilde{\mathbb{P}}$ on the underlying probability space $(\Omega, \mathcal{B}, \mathbb{P})$ such that the transition probabilities of the modified system, when viewed under $\tilde{\mathbb{P}}$, coincide with those of the original system under the original probability measure $\mathbb{P}$. This approach enables us to prove the finite-time error estimate \eqref{Eq_aux_expl} for the radial component of the transformed system with respect to the new probability measure $\tilde{\mathbb{P}}$, while preserving the statistics of the original system.

This technique draws inspiration from \cite{Hairer2001ExponentialMP}, where Girsanov's theorem was employed to derive exponential mixing bounds for a class of SPDEs through a coupling argument.
Our work is also related to the findings of \cite{B2004} and \cite{B2005}, which investigate a broad class of SPDEs with cubic nonlinearity. These studies establish upper bounds on the Wasserstein and total variation distances between the invariant measures of the original SPDE and its finite-dimensional amplitude equation under the assumption of timescale separation.

Another originality of this article is the derivation of error estimates for parameter regimes in which such assumptions about time-scale separation do not necessarily hold. In that respect, the deterministic slow manifold is replaced by a family of SPMs. The construction of these stochastic manifolds is operated as follows.  For a given $\tau$ in $(0,\infty)$ we define the mapping 
\bea \label{Eq_stoch_mnf_int}
h_{\tau}: \, & \mathbb{R} \times \mathbb{R}_+ \longrightarrow \mathbb{R},\\
&(m, r) \mapsto m +  c_{\tau} r^2,\,\,\,\; c_{\tau} = \left(1 - e^{-\frac{\tau}{\epsilon}}\right).
\eea 
The  family of SPMs, $(h_{\tau}(M_t, \cdot))_{t\geq0}$, is then obtained by replacing $m$ in \eqref{Eq_stoch_mnf_int} by an auxiliary stationary process $(M_t)_{t\geq 0}$ depending on the past of the noise such as defined in  Section \ref{Sec_intro_stoch_pm} below.

The parameter $\tau$ plays the role of a tuning parameter to adjust the shape of the deterministic part of the manifold $h_{\tau}$. This way, SPMs offer an extension to the slow manifold, $h(r)=r^2$, as the latter is recovered when $m=0$ and $\tau \rightarrow \infty$, but allow for more flexibility in the design of the relevant parameterizations, in particular by finding an optimal $\tau$ that minimizes the parameterization defect; see \cite[Sec.~4]{CLM19_closure}.

In Section \ref{Sec_stoch_pm}, we demonstrate both analytically and numerically that Stochastic Parameterizing Manifolds (SPMs) often serve as effective reduced-order models for system \eqref{Eq_3DHopf}, particularly in the challenging regimes of $\epsilon>1$ near criticality ($\lambda$ small). We derive bounds on the discrepancy between the marginals of the invariant measures of the original and reduced systems along the radial component, analogous to \eqref{Eq_est_int}. These bounds are formalized in Theorems \ref{Thm_hopf_stoch} and \ref{Thm_hopf_22}, counterparts to Theorems \ref{Thm_hopf_1} and \ref{Thm_hopf_2}, respectively. Section \ref{Sec_numerics_2} presents numerical results that support these findings.

Importantly, the Girsanov approach presented here offers a unified framework for deriving error bounds of the form \eqref{Eq_est_int} for both deterministic slow manifolds (Theorems \ref{Thm_hopf_1} and \ref{Thm_hopf_2}) and SPMs (Theorems \ref{Thm_hopf_stoch} and \ref{Thm_hopf_22}).

This article is structured as follows. At the beginning of Section \ref{Sec_slow_pm} we give a brief introduction to system \eqref{Eq_3DHopf} and its slow manifold reduction. In Section \ref{Sec_main_idea} we discuss the main ideas for the proof of Theorems \ref{Thm_hopf_1} and \ref{Thm_hopf_2}, the main results of Section \ref{Sec_slow_pm}. A rigorous treatment of Girsanov's theorem applied to our situation is provided in Section \ref{Sec_girs}. The proofs of Theorems \ref{Thm_hopf_1} and \ref{Thm_hopf_2} are then detailed in Sections \ref{Sec_proof_thm_1} and \ref{Sec_proof_thm_3}, respectively. Numerical results in Section \ref{Sec_numerics_1} conclude Section \ref{Sec_slow_pm}. Section \ref{Sec_stoch_pm} first introduces the concept of SPMs, before the main results of this section (Theorems \ref{Thm_hopf_stoch} and \ref{Thm_hopf_22}) are stated and proven in Section  \ref{Sec_main_thm_stoch_pm} and \ref{Sec_main_thm_stoch_pm2}, respectively. 
 Numerical simulations in Section \ref{Sec_numerics_2} complement our theoretical findings, providing valuable insights into the model's behavior, particularly in the context of inverted timescale separation. These simulations explore parameter regimes beyond the scope of our analytical results, offering further evidence of the efficacy of our approach in these challenging scenarios.

\section{Slow parameterizing manifolds: Close to criticality and $\epsilon \ll 1$}
\label{Sec_slow_pm}

\subsection{Preliminaries}\label{Sec_prelim}
As recalled in Introduction, in the absence of noise (i.e. $\sigma=0$) system \eqref{Eq_3DHopf} exhibits a deterministic {\it slow manifold}, that is obtained by setting the drift term in {\mkr Eq.~\eqref{Eq_3DHopf_z} to zero and solving for $z$, yielding to} a function $h$
defined by
\bea \label{Eq_slow_mnf}
h: \, & \mathbb{R}\times \mathbb{R} \longrightarrow \mathbb{R},\\
&(x,y) \mapsto x^2 + y^2.
\eea
For $\sigma=0$ and $\epsilon \ll 1$, the solutions to system \eqref{Eq_3DHopf} are expected to rapidly converge to this slow manifold up to a small error dominated by $\epsilon$; see e.g.~ \cite[Theorem 10.1]{nipp2013invariant}. 
We investigate hereafter, in the case $\sigma > 0$, whether there exist parameter regimes where the deterministic slow manifold remains relevant, providing an approximate slaving relationship between the $z$-variable and the $x$-, $y$-variables.
 If one takes $h$ as a parametrization of the  "fast" variable $z$ in terms of the "slow" variables $x$ and $y$, we arrive formally at the following reduced system
\bea\label{Eq_2DHopf}
\d u&=\big(\lambda u- f  v -\gamma  u ( u^2 +  v^2) \big)\d t + \sigma \d W^1_t\\
\d v &=\big(f  v +\lambda  v -\gamma  v (u^2 +  v^2) \big)\d t+ \sigma \d W^2_t.
\eea
The reduced system \eqref{Eq_2DHopf} is referred to as the slow manifold reduced system.

For $\mathbb{T}^1$ {\mkr denoting the unit circle identified with the interval} $[0, 2\pi)$ (modulo $2\pi$) we introduce the {\mkr state spaces} 
\be \label{Def_M_hat_M}
\hat{\mathcal{M}} = \mathbb{R}_{+} \times \mathbb{T}^1\,\,\,\text{ and }\,\,\, \mathcal{M} = \hat{\mathcal{M}} \times \mathbb{R}.
\ee
When $\sigma = 0$, system \eqref{Eq_2DHopf} {\mkr corresponds to a particular case of the} Hopf normal form \cite[p.~86-90]{Kuznetsov:1998:EAB:289919}. For $\sigma > 0$ system \eqref{Eq_2DHopf} {\mkr corresponds thus to} the Hopf normal form subject to additive noise; we refer to \cite{doan2018hopf,tantet2020ruelle} for  analyses of such a system.

In a first step we write the reduced system \eqref{Eq_2DHopf} in polar coordinates, {\mkr leading to a system on} $\hat{\mathcal{M}}$ of the form
\begin{subnumcases}{\label{Eq_2DHopf_polar_base}}
\d \hat{r} = \left[\lambda \hat{r} -\gamma \hat{r}^3 + \frac{\sigma^2}{2\hat{r}}\right]\d t  + \sigma \d W^{\hat{r}}_t \label{Eq_r_hat}\\
\d \hat{\theta} = f \d t +  \frac{\sigma}{\hat{r}} \d W^{\hat{\theta}}_t \label{Eq_theta_hat} \mod (2\pi).
\end{subnumcases}
Here the stochastic processes $W^{\hat{r}}_t$ and $W^{\hat{\theta}}_t$  satisfy the following system, for all $t\geq0$,
\bea\label{Eq_Wr}
\d W^{\hat{r}}_t &= \cos(\hat{\theta}_t) \d W^1_t + \sin(\hat{\theta}_t) \d W^2_t\\
\d W^{\hat{\theta}}_t &= - \sin(\hat{\theta}_t) \d W^1_t + \cos(\hat{\theta}_t) \d W^2_t.
\eea

To compare the qualitative behavior of the original system \eqref{Eq_3DHopf} {\mkr with that of} the reduced system \eqref{Eq_2DHopf}, we write system \eqref{Eq_3DHopf} in {\mkr cylindrical} coordinates and obtain on $\mathcal{M}$
\begin{subnumcases}{\label{Eq_3DHopf_polar_base}}
\d r = \left[\lambda r -\gamma r z + \frac{\sigma^2}{2 r}\right]\d t + \sigma \d W^r_t \label{Eq_3DHopf_polar_base_r1}\\
\d \theta = f \d t +  \frac{\sigma}{r} \d W^{\theta}_t \mod (2\pi) \\
\d z = -\frac{1}{\epsilon}\left[z -r^2\right]\d t + \frac{\sigma}{\sqrt{\epsilon}} \d W^3_t \label{Eq_3DHopf_polar_base_z1}.
\end{subnumcases}
Here $W^r_t$ and $W^{\theta}_t$ satisfy Eq.~\eqref{Eq_Wr} with $\theta_t$ replacing $\hat{\theta}_t$. 	The Brownian motion $W^3_t$ is as in Eq.~\eqref{Eq_3DHopf_z}.

 {\mkr To the Brownian motions $W^r_t,W^{\theta}_t$, and $W^3_t$ driving Eq.~\eqref{Eq_3DHopf_polar_base}, we associate} the filtration $\{\mathcal{B}_t\}_{t\geq 0}$  {\mkr (i.e.~an increasing family of sub $\sigma$-algebras of (the completion of) $\mathcal{B}$)} defined by
\be \label{Eq_filtration}
\mathcal{B}_t = \sigma\left\{W^r_s, W^{\theta}_s, W^3_s: 0\leq s\leq t \right\}, \,\,\, t \geq 0.
\ee
It corresponds to the  smallest filtration with respect to which $W^r_t,W^{\theta}_t$, and $W^3_t$ are adapted, i.e.~$\mathcal{B}_t$-measurable. The filtration gives a direction to time and allows one to speak of the past etc.~of $t$, thus getting closer to a physical concept of time. Invariant measure used below means that the probability measure is stationary with respect to the corresponding Markov transition semigroup.

Obviously, system \eqref{Eq_3DHopf_polar_base} is a slow-fast system, where the slow variables are now formed by the radial variable, $r$, and the angular variable, $\theta$.  Within this new coordinate system, the slow manifold is given by
	\bea \label{Eq_slow_mnf_polar}
	h_{\mathrm{polar}}: \, & \mathbb{R}_+ \longrightarrow \mathbb{R}_+,\\
	&r \mapsto r^2.
	\eea
Geometrically, as $\lambda >0$,  the origin is unstable and due to the dissipative and rotation effects caused respectively by the nonlinear and the linear terms,  system \eqref{Eq_2DHopf_polar_base} exhibits for $\sigma=0$ a globally stable limit cycle $\mathcal{C}$, whose 
radius is given by
	\be\label{Eq_r_det}
	r_{\det} = \sqrt{\frac{\lambda}{\gamma}}.
	\ee
When $\epsilon\ll 1$ and $\sigma=0$, system \eqref{Eq_3DHopf_polar_base} exhibits also a globally stable limit given as $\Gamma=\mathcal{C}\times \{r_{\det}^2\}$ in $\mathcal{M}$.
Intuitively, when the noise is turned on, at least for small noise intensity, one expects to have the dynamics associated with  Eq.~\eqref{Eq_3DHopf_polar_base} to meander near the limit cycle $\Gamma$. To analyze rigorously this idea, and in particular in order to identify parameter regimes for which the slow manifold provides a valid approximation, we aim at comparing  the invariant measure $\mu$ of the original system \eqref{Eq_3DHopf_polar_base} --- after projection to the $(r,\theta)$-plane --- with  the invariant measure, $\nu$,  of the slow manifold reduced system \eqref{Eq_2DHopf_polar_base}; see Remark \ref{Rmk_problems_to_overcome}-(i) below for the existence of such invariant measures. Since the $\theta$-equation is the same as in Eqns.~\eqref{Eq_3DHopf_polar_base} and \eqref{Eq_2DHopf_polar_base}, it boils down to comparing only the marginals $\mu_r$ and $\nu_r$ along the $r$-direction of the invariant measures $\mu$ and $\nu$.

For this purpose, we consider the following integral probability metric for probability measures on $\mathbb{R}^n$, for $n$ in $\mathbb{N}$.
The set of test functions (against probability measures) that we consider is given by
\be \label{Eq_test_fcts}
\mathcal{F}^{(n)}=\{ \varphi:\mathbb{R}^n \rightarrow \mathbb{R}\, \mbox{ measurable }|\, \|\varphi\|_{\textrm{Lip}} \leq 1, \; \varphi(0)=0\},
\ee
where $\|\cdot\|_{\textrm{Lip}}$ is defined as
\be \label{Eq_lip}
\|\varphi\|_{\textrm{Lip}} = \sup_{x,y \in \mathbb{R}^n, x\neq y} \left\{\frac{|\varphi(x) - \varphi(y)|}{|x - y|} \right\}.
\ee
Hereafter we only consider the spaces $\mathcal{F}^{(n)}$ for $n=1$. For brevity of notation we therefore set 
\be \label{Eq_test_fcts_1}
\mathcal{F} = \mathcal{F}^{(1)}.
\ee
 We set also  
\be\label{Eq_IPM}
D_{\mathcal{F}^{(n)}}(\mathfrak{m},\mathfrak{n})= \underset{\varphi \in \mathcal{F}^{(n)}}\sup\Big|
\int \varphi \d\m -\int \varphi \d \n\Big|,
\ee
for any probability measures $\m$ and $\n$ on $\mathbb{R}^n$.

There is actually an equivalence between the $L^1$-Wasserstein distance $d_{W_1}$ (e.g.~\cite{Villani:2003}) and $D_{\mathcal{F}^{(n)}}$ on the following set of probability measures 
	\be
	{Pr}_0(\mathbb{R}^n)= \left\{\m \in Pr(\mathbb{R}^n): \int_{\mathbb{R}^n}|x|\m(\d x) < \infty,\,\, \m(\{0\}) = 0 \right\}.
	\ee
More precisely we have
\bprop \label{Remark_Wasserstein_metric}
	Let $\m$ and $\n$ be probability measures in $Pr_{0}(\R^n)$, then it holds that 
	\be
	D_{\mathcal{F}^{(n)}}(\m, \n) = d_{W_1}(\m, \n).
	\ee
	\eprop
	\begin{proof}
		Let $\m$ and $\n$ be in $Pr_{0}(\R^n)$.
		Note that by removing the condition $\varphi(0) = 0$ in \eqref{Eq_test_fcts}, $D_{\mathcal{F}^{(n)}}$ coincides with $d_{W_1}(\m, \n)$. We therefore have
		\bes
		D_{\mathcal{F}^{(n)}}(\m, \n) \leq d_{W_1}(\m, \n).
		\ees
		Furthermore, it is easy to see that for all Lipschitz continuous functions $\varphi:\R^n \rightarrow \R$ with $\|\varphi\|_{\textrm{Lip}}\leq 1$ there exists a function $\tilde{\varphi}$ in $\mathcal{F}^{(n)}$ such that 
		\bes
		\varphi = \tilde{\varphi} + \varphi(0).
		\ees
		Hence
		\bes
		\left|\int_{\R^n} \varphi \d \m - \int_{\R^n}\varphi \d \n \right| = \left|\int_{\R^n} \tilde{\varphi} \d \m - \int_{\R^n}\tilde{\varphi} \d \n \right|,
		\ees
		and 
		\be \label{Eq:wasserstein_is_my_metricI}
		\left|\int_{\R^n} \varphi \d \m - \int_{\R^n}\varphi \d \n \right| \leq \sup_{\tilde{\varphi} \in \mathcal{F}^{(n)}} \left|\int_{\R^n} \tilde{\varphi} \d \m - \int_{\R^n}\tilde{\varphi} \d \n \right|.
		\ee
		Inequality \eqref{Eq:wasserstein_is_my_metricI} holds for all Lipschitz continuous $\varphi:\R^n \rightarrow \R$ with $\| \varphi \|_{\textrm{Lip}}\leq 1$. By taking the supremum over all such $\varphi$ on both sides of \eqref{Eq:wasserstein_is_my_metricI} we arrive at
		\bes
		d_{W_1}(\m, \n) \leq D_{\mathcal{F}^{(n)}}(\m, \n),
		\ees
		which completes the proof.	
	\end{proof}

\needspace{1\baselineskip}
\br\label{Rmk_problems_to_overcome}
\hspace*{2em}  \vspace*{0.1ex}
\bi
\item[(i)] Lemmas \ref{Lemma_invm_ex_org} and \ref{Lemma_invm_exuniq_red} recalled in Appendix \ref{App_ex_uniq_org} below, ensure the existence of unique ergodic invariant measures $\mu$ and $\nu$, 
	for the original system \eqref{Eq_3DHopf_polar_base} and for the slow manifold reduced system \eqref{Eq_2DHopf_polar_base}, respectively.
	Working with ${Pr}_0(\mathbb{R}^n)$ instead of $Pr(\mathbb{R}^n)$ does not pose any serious limitation, as the invariant measures $\mu$ and $\nu$ also admit smooth Lebesgue-densities according to Lemmas \ref{Lemma_invm_ex_org} and \ref{Lemma_invm_exuniq_red}. Therefore, $\mu(\{0\}) = 0 = \nu(\{0\})$.
\item[(ii)]There is a variety of metrics of type \eqref{Eq_IPM} that are often tailored to specific tasks, see e.g. \cite{2009arXiv0901.2698S}. 
The Wasserstein metric $D_\mathcal{F}$ turned out to be relevant to derive useful estimates that relate explicitly to the system's parameters.  
\ei
\er

Given a reduced state space $V$ of an ambient Euclidean space $X$, we {\mkr denote by  $\Pi_V $ the canonical projection, namely}
\bea \label{Def_proj}
	\Pi_V :\, & X\rightarrow V\\
	& \bm{x} \mapsto \Pi_V  \bm{x} =\bm{v}.
\eea
Let $e_r, e_\theta$ and $e_z$ denote the vectors of $\mathbb{R}^3$ forming its cylindrical coordinate system.
Associated with the subspaces {\mkr $V_{r,z} = \span\{e_r, e_z\}$ and $V_r = \span\{e_r\}$,  we define the  following 
marginal probability measures
\be \label{Eq_marginals_1}
\mu_{r,z} = \Pi_{V_{r,z}}^* \mu, \,\, \mu_r = \Pi_{V_r}^* \mu, \,\, \nu_r = \Pi_{V_r}^* \nu,
\ee
where  $\mu$ and $\nu$ denote the  ergodic invariant measures of Eq.~\eqref{Eq_3DHopf_polar_base} and  Eq.\eqref{Eq_2DHopf_polar_base}, respectively (see Remark \ref{Rmk_problems_to_overcome}-($i$)),  while 
$\Pi_V$ denotes the projector onto $V=V_{r,z}$ or $V=V_r$, as defined in \eqref{Def_proj}.

{\mkr Theorems \ref{Thm_hopf_1} and \ref{Thm_hopf_2} proved below assert that, under mild conditions, 
the (Wasserstein) distance $D_{\mathcal{F}}(\mu_r, \nu_r)$ between the marginal probabilities along the radial direction, of the original and the reduced dynamics, are 
controlled  by constants $C(\lambda, \gamma,\sigma)>0 $ and $c>0$  
such that} 
\be \label{Eq_est}
D_{\mathcal{F}}(\mu_r, \nu_r) \leq C(\lambda, \gamma, \sigma) + c \left(\int_{\mathbb{R}_+\times \mathbb{R}} |z - r^2|^{4}  \, \mu_{r,z}(\d r,\d z)\right)^{\frac{1}{4}}.
\ee
{\mkr Furthermore, Theorems \ref{Thm_hopf_1} and \ref{Thm_hopf_2}  provide explicit expressions of these constants involving the original system's parameters as well as its invariant measure $\mu$. A discussion of such a result is postponed to Section \ref{Subsec_statement_thm_1} in which the main theorems of this article are presented. In the meantime, we present next the main ideas behind the proof of these theorems.}

\subsection{Main ideas behind the proof of Theorems \ref{Thm_hopf_1} and \ref{Thm_hopf_2}}
\label{Sec_main_idea}
In this subsection we outline the main difficulties, as well as the main ideas and techniques used in the proof of {\mkr the estimate} \eqref{Eq_est}. 
These ideas constitute the core of the proofs of Theorems \ref{Thm_hopf_1} and \ref{Thm_hopf_2}, the main results of this section.  Our analysis takes place for the filtered probability space $(\Omega, \mathcal{B}, \{\mathcal{B}_t\}_{t\geq 0} , \mathbb{P})$ as recalled in Section \ref{Sec_prelim}.

A crucial step in deriving error estimate \eqref{Eq_est} involves establishing the following auxiliary finite-time error estimate:
\be  \label{Eq_aux_expl}
\mathbb{E}_{\mathbb{P}} \left[|\hat{r}_t - r_t|^2\right] \leq e^{-qt}|\hat{r}_0 - r_0|^2 + \frac{\gamma^2}{q} \int_0^t e^{-q(t-s)} \mathbb{E}_{\mathbb{P}}\left[ r^2_s(z_s-r^2_s)^2\right]\d s,\,\,\,\text{for all}\,\,\,t\in [0,T],
\ee
where $q$ is a suitable positive constant defined below.

A direct approach to estimating $\hat{r}_t-r_t$ based on Eq.~\eqref{Eq_3DHopf_polar_base} and Eq.~\eqref{Eq_2DHopf_polar_base} by taking the difference and applying Gronwall's lemma yields an estimate similar to \eqref{Eq_aux_expl} but with $-q=\lambda >0$, where $\lambda$ is the positive growth rate associated with the oscillatory instability. This leads to a significant issue: the RHS of the inequality blows up as the time horizon $T$ approaches infinity.\footnote{As we demonstrate below, deriving \eqref{Eq_est} from \eqref{Eq_aux_expl} requires taking the limit as $T\rightarrow \infty$ via an ergodic argument.

This limitation arises directly from the inherent instability of the origin, as characterized by $\lambda>0$.
To overcome this challenge, we introduce a novel approach. We modify the radial equation of the original system by incorporating a coupling term that explicitly depends on the radial part, $\hat{r}_t$, of the solution to the slow manifold reduced system \eqref{Eq_2DHopf_polar_base}. This modification is carefully designed to eliminate the linear instability when taking the difference between the modified and reduced systems. This leads to a differential inequality of the form $\dot{u}\leq -q u +f(u)$, where $q>0$, which can be effectively analyzed using Gronwall's lemma; see Appendix \ref{Appendix_aux_ineq}.}

We now elucidate the construction of this coupling term. Let $q$ be a positive constant, which is explicitly defined in \eqref{Eq_def_q_girs} below, and let $g$ be the function defined as follows
\bea \label{Eq_g}
g:\,& \mathbb{R}_+\times \mathbb{R}_+ \longrightarrow \mathbb{R},\\
&(r_1,r_2) \mapsto  \frac{q + \lambda}{\sigma} (r_2 - r_1).
\eea

For $\hat{r}_t$ solving Eq.~\eqref{Eq_2DHopf_polar_base}, we consider then the following stochastic system obtained by adding the term  $\sigma g(\tilde{r}_t, \hat{r}_t)$  to the drift part of Eq.~\eqref{Eq_3DHopf_polar_base}, namely 
\begin{subnumcases}{\label{Eq_3DHopf_polar_base_transf}}
\d \tilde{r} = \left[\lambda \tilde{r} -\gamma \tilde{r} \tilde{z} + \frac{\sigma^2}{2 \tilde{r}}\right]\d t + \sigma g(\tilde{r}, \hat{r}) \d t + \sigma \d W^r_t \label{Eq_r_transf}\\
\d \tilde{\theta} = f \d t +  \frac{\sigma}{\tilde{r}} \d W^{\theta}_t \mod (2\pi) \\
\d \tilde{z} = -\frac{1}{\epsilon}\left[\tilde{z} -\tilde{r}^2\right]\d t + \frac{\sigma}{\sqrt{\epsilon}} \d W^3_t.
\end{subnumcases}
The function $g$ in \eqref{Eq_r_transf} introduces a coupling with $\hat{r}_t$ from Eq.~\eqref{Eq_r_hat}, which results in a modified version of $(r_t,\theta_t,z_t)$ as denoted by the process $(\tilde{r} _t,\tilde{\theta}_t, \tilde{z}_t)$ satisfying thus \eqref{Eq_3DHopf_polar_base_transf}.
The particular choice of $g$ in \eqref{Eq_g} then allows us to show that
\be \label{Eq_lemma_grw_5_2}
\frac{\d}{\d t}\mathbb{E}_{\mathbb{P}}\left[|\hat{r}_t - \tilde{r}_t|^2 \right] \leq  -q\mathbb{E}_{\mathbb{P}}\left[|\hat{r}_t - \tilde{r}_t|^2\right]  + \frac{\gamma^2}{q} \mathbb{E}_{\mathbb{P}}\left[ \tilde{r}^2_t\left(\tilde{z}_t-\tilde{r}^2_t\right)^2 \right],\,\,\,\text{for all}\,\,\,t \in [0,T];
\ee 
see Appendix \ref{Appendix_aux_ineq} for more details.

By applying Gronwall's inequality to \eqref{Eq_lemma_grw_5_2} we arrive at the finite-time error estimate \eqref{Eq_aux_expl}  satisfied by $\tilde{r}_t$ and $\tilde{z}_t$.  We have that
\be  \label{Eq_aux_expl2}
\mathbb{E}_{\mathbb{P}} \left[|\hat{r}_t - \tilde{r}_t|^2\right] \leq e^{-qt}|\hat{r}_0 - \tilde{r}_0|^2 + \frac{\gamma^2}{q} \underbrace{\int_0^t e^{-q(t-s)} \mathbb{E}_{\mathbb{P}}\left[ \tilde{r}^2_s(\tilde{z}_s-\tilde{r}^2_s)^2\right]\d s}_{\mathcal{I}(\tilde{r},\tilde{z})},\,\,\,\text{for all}\,\,\,t\in [0,T].
\ee

Note that the integral term $\mathcal{I}(\tilde{r},\tilde{z})$ on the RHS of inequality \eqref{Eq_aux_expl2} exhibits similarities to the parameterization defect in \eqref{Eq_est}. However, several key distinctions exist:
\bi
 \item[(i)] The expectation in \eqref{Eq_aux_expl2} is taken with respect to the probability measure $\mathbb{P}$, whereas the expectation in \eqref{Eq_est} is taken with respect to the invariant measure $\mu_{r,z}$.
 \item[(ii)]  $\mathcal{I}(\tilde{r},\tilde{z})$ involves the variables of the transformed system.
 \item[(iii)] The integrand in \eqref{Eq_aux_expl2} involves a product between the variables $\tilde{r}^2_s$ and $(\tilde{z}_s-\tilde{r}^2_s)^2$, in contrast to the integrand defined in \eqref{Eq_est}, involving  $(z-r^2)^4$.
 \ei
Our primary objective is to  guide the reader on how to effectively replace $\mathcal{I}(\tilde{r},\tilde{z})$ with $\mathcal{I}(r,z)$ in \eqref{Eq_aux_expl2}. With this replacement performed, the issue $(i)$ is then addressed in Step 3 of the proof of Theorem \ref{Thm_hopf_1} in Section \ref{Sec_proof_thm_1} while $(iii)$ is handled through standard Cauchy-Schwarz estimates; see estimates after \eqref{Eq_proof_thm_1_dep_phi} therein. 

It is important to note though that a direct replacement of $\mathcal{I}(\tilde{r},\tilde{z})$ by $\mathcal{I}(r,z)$  in \eqref{Eq_aux_expl2} is not feasible. This is because the transition probabilities of the transformed system \eqref{Eq_3DHopf_polar_base_transf} do not necessarily coincide with those of the original system \eqref{Eq_3DHopf_polar_base} when both are evaluated under the original probability measure $\mathbb{P}$. In other words, the transformation involving the function $g$ in \eqref{Eq_g} does not generally preserve the statistics of the system when viewed under the original probability measure $\mathbb{P}$.

In order to overcome this difficulty, we employ Girsanov's theorem, which is stated in its general form in Theorem \ref{Lemma_girsanov_general} below and is shown to apply to our context in Lemma \ref{Lemma_girsanov}.  
More concretely, let $\big(R_t^{R_0}\big)_{t\geq0}$ denote the unique strong solution (e.g.~\cite[Definition 1.6, p. 163]{watanabe_ikeda}). to the original system \eqref{Eq_3DHopf_polar_base} emanating from an initial datum $R_0$ in $\mathcal{M}$. Analogously, $\big(\tilde{R}_t^{R_0}\big)_{t\geq0}$ denotes the unique strong solution of the transformed  system \eqref{Eq_3DHopf_polar_base_transf} emanating from the same initial datum $R_0$ in $\mathcal{M}$. We refer to Remark \ref{Remark_ex_str_sol} below for a discussion on the existence and uniqueness of strong solutions to systems \eqref{Eq_3DHopf_polar_base}, \eqref{Eq_2DHopf_polar_base} and \eqref{Eq_3DHopf_polar_base_transf}. 
Girsanov's theorem ensures that for every $T>0$ there exists a probability measure $\tilde{\mathbb{P}} = \tilde{\mathbb{P}}_T$ on $(\Omega, \mathcal{B})$, which is equivalent to $\mathbb{P}$ on $\mathcal{B}_T$ such that
\be \label{Eq_trans_prob}
\tilde{\mathbb{P}}\Big(\tilde{R}_t^{R_0} \in \Gamma\Big) = \mathbb{P}\Big(R_t^{R_0} \in \Gamma\Big),\,\,\, \Gamma \in \mathcal{B}( \mathcal{M}), \; R_0 \in \mathcal{M}, \,\,\,t\in [0,T].
\ee

Girsanov's theorem allows us to preserve the transition probabilities, provided the probability measure $\tilde{\mathbb{P}}$ is used to assess the transition probabilities for the transformed system. 
Thus, with this change of probability measure, we arrive at
\be \label{Eq_aux_expl3}
\mathbb{E}_{\tilde{\mathbb{P}}} \left[|\hat{r}_t - \tilde{r}_t|^2\right] \leq e^{-qt}|\hat{r}_0 - \tilde{r}_0|^2 + \frac{\gamma^2}{q} \int_0^t e^{-q(t-s)} \mathbb{E}_{ \mathbb{P}}\left[ r^2_s(z_s-r^2_s)^2\right]\d s,\,\,\,\text{for all}\,\,\,t\in[0,T],
\ee
with $\mathbb{E}_{\tilde{\mathbb{P}}} \left[|\hat{r}_t - \tilde{r}_t|^2\right]=\mathbb{E}_{\mathbb{P}} \left[|\hat{r}_t - r_t|^2\right]$ according to \eqref{Eq_trans_prob}. We refer to Lemma \ref{Lemma_aux_ineq_1} in Appendix \ref{Appendix_aux_ineq} for a detailed proof of the auxiliary inequality \eqref{Eq_aux_expl3}. 

The fact that Girsanov's theorem just holds on finite-time intervals does not pose a serious difficulty for obtaining the main results in Theorems \ref{Thm_hopf_1} and \ref{Thm_hopf_2}. In the final step in the proof of Theorems \ref{Thm_hopf_1} and \ref{Thm_hopf_2} we exploit that the final time $T>0$ can be chosen large enough for our purposes. This, together with the fact that the original and reduced systems are ergodic eventually yields the proof.

\br\label{Connection_Hairer}
The introduction of the transformed  system \eqref{Eq_3DHopf_polar_base_transf} is inspired by \cite{Hairer2001ExponentialMP}, where a probabilistic coupling approach was used to obtain exponential mixing bounds. Therein the author considers two copies of an original system, where the equation of one copy\footnote{In  \cite{Hairer2001ExponentialMP} a copy consists of the same original system with different initial data.} is explicitly coupled to the second copy. By virtue of Girsanov's theorem this construction gives rise to a suitable joint distribution between the laws of both systems, which is referred to as probabilistic coupling.
An important difference and novelty of our approach consists of coupling the original system to the reduced system \eqref{Eq_2DHopf_polar_base} {\mkr (through the $g$-term)}, leading to the transformed system \eqref{Eq_3DHopf_polar_base_transf}. Using this coupling idea outlined above we are then able to derive estimates of the form \eqref{Eq_est} that relate the discrepancy between the long-term statistics of the original and reduced system to the parameterization defect.
\er

\subsection{Preserving transition probabilities via the Girsanov transformation} 
\label{Sec_girs}
In this subsection we provide a detailed treatment of the ideas relating to the Girsanov approach outlined in Section \ref{Sec_main_idea}.
The following theorem is a restatement of Girsanov's theorem in \cite[Theorem 10.14, p. 290]{daprato_zabczyk_2014} for finite dimensions.
\bt \label{Lemma_girsanov_general}
	Let $d$ be a positive integer and let $W_t$ be a $\mathbb{R}^d$-valued standard Brownian motion on a filtered probability space $(\Omega, \mathcal{B}, \{\mathcal{B}_t\}_{t\geq0}, \mathbb{P})$. 
Assume that $U_t$ is a $\mathbb{R}^d$-valued and predictable in the sense of \cite[Definition 5.2, p. 21]{watanabe_ikeda}. stochastic process such that 
	the  random variable $\mathfrak{D}_T$ ($T>0$),
	\bes
	\mathfrak{D}_T = \exp\left(\int_0^T U_s\d W_s - \frac{1}{2} \int_0^T \left|U_s\right|^2\d s \right),
	\ees
	is integrable with respect to $\mathbb{P}$, and 
	\be \label{Eq_girs_identity_1}
	\mathbb{E}_{\mathbb{P}}\left[\mathfrak{D}_T\right] =1.
	\ee

Then $\tilde{\mathbb{P}}_T(\d \omega) = \mathfrak{D}_T(\omega) \mathbb{P}(\d \omega)$ defines a probability measure on $(\Omega, \mathcal{B})$, which is equivalent\footnote{Meaning that, for all $A$ in $\mathcal{B}$ it holds that $\mathbb{P}(A) = 0$, if an only if $\tilde{\mathbb{P}}(A) = 0$.} to $\mathbb{P}$. Furthermore, the stochastic process $(\tilde{W}_t)_{t\in [0,T]}$ given by
	\bes
	\tilde{W}_t = W_t - \int_0^tU_s\d s, \,\,\, t \in [0,T],
	\ees
	is a Brownian motion with respect to $\tilde{\mathbb{P}}_T$, i.e.~it is a Brownian motion on the filtered probability space $(\Omega, \mathcal{B}, \{\mathcal{B}_t\}_{t\in [0,T]}, \tilde{\mathbb{P}}_T)$.
\et
Lemma \ref{Lemma_girsanov} below expresses that Theorem \ref{Lemma_girsanov_general} applies to the transformed system \eqref{Eq_3DHopf_polar_base_transf}. For that purpose, identity \eqref{Eq_girs_identity_1} is shown hereafter with $U_t=\betrag{g(\tilde{r}_s, \hat{r}_s)}^2$ ith $g$ defined in \eqref{Eq_g}. It follows that transition probabilities of the original system \eqref{Eq_3DHopf_polar_base} are preserved for the transformed  system \eqref{Eq_3DHopf_polar_base_transf}, {\mkr when the transition probabilities of the latter system are assessed with respect to the transformed probability measure $\tilde{\mathbb{P}}$.}The following Remark \ref{Remark_ex_str_sol} ensures the existence of global unique strong solutions to all systems under consideration in this section.
	\br \label{Remark_ex_str_sol}
	Note that the drift and diffusion coefficients of the original systems, both in Cartesian and polar coordinates in \eqref{Eq_3DHopf} and \eqref{Eq_3DHopf_polar_base}, respectively, are locally Lipschitz continuous. The same {\mkr observation} holds for the reduced systems in \eqref{Eq_2DHopf} and \eqref{Eq_2DHopf_polar_base}, respectively. This yields to the existence and uniqueness of (local) strong solutions to {\mkr each of} these systems on the filtered probability space $(\Omega, \mathcal{B}, \{\mathcal{B}_t\}_{t\geq 0} , \mathbb{P})$ by means of standard results of the theory of SDEs; see e.g.~\cite[Theorem 2.3, p.~173]{watanabe_ikeda} and \cite[Theorem 3.1, p.~178-179]{watanabe_ikeda}. Furthermore, regarding the reduced systems, we note that \cite[Theorem 3.6, p.~58-59]{MaoXuerong2008Sdea} ensures global existence and uniqueness of strong solutions. 

However, the existence of global solutions for the original systems is not immediate. This is in particular due to the fact that that standard arguments tied to energy-preservation are not applicable here, as the quadratic nonlinearity {\mkr of e.g.~the original system  \eqref{Eq_3DHopf}} is not energy preserving {\mkr when $\gamma \neq 1/\epsilon$}, i.e.
	\be
	\left\langle\begin{pmatrix} 
		-\gamma xz \\
		-\gamma yz \\
		\frac{1}{\epsilon}\left(x^2 + y^2\right) 
	\end{pmatrix},\begin{pmatrix} 
		x \\
		y \\
		z 
	\end{pmatrix}\right\rangle
	\neq 0.
	\ee
It thus precludes the direct application of \cite[Theorem 3.6, p.~58-59]{MaoXuerong2008Sdea} to ensure global existence of solutions.  To bypass this difficulty, the Lyapunov function
approach (see e.g.~\cite[Theorem 3.5, p.~75]{khasminskii2011stochastic}) offers a well-established alternative. Essentially, given an SDE with generator $L$, this approach asserts that global solutions exist if a Lyapunov function $V$ can be found that satisfies the following inequality:	
\be \label{Eq_lyapunov_ex_glob}
	LV \leq cV, \;\; c>0.
	\ee
	
For the original system \eqref{Eq_3DHopf}, a Lyapunov function is given by
	\bea \label{Eq_V_ex1}
	V:\,& \R\times \R\times \R \longrightarrow \mathbb{R},\\
	&(x, y, z) \mapsto  x^2 + y^2 + \sqrt{\epsilon \gamma}z^2 + 1.
	\eea
This result naturally extends to the original system \eqref{Eq_3DHopf_polar_base} expressed in polar coordinates (see Appendix \ref{App_ex_uniq_org} for details). This is further corroborated by inequality \eqref{Eq_KV_leq_aV} derived in Step 4 of the proof of Lemma \ref{Lemma_invm_ex_org}. The existence of unique global strong solutions for the transformed system \eqref{Eq_3DHopf_polar_base_transf} can be established analogously.
 \er 

\bl \label{Lemma_girsanov}
Let $T>0$. Recall the filtration $\{\mathcal{B}_t\}_{t\geq0}$ introduced in \eqref{Eq_filtration} and the corresponding filtered probability space $(\Omega, \mathcal{B}, \{\mathcal{B}_t\}_{t\geq0},\mathbb{P})$. The {\mkr state spaces}, $\mathcal{M}$ and $\hat{\mathcal{M}}$, are as in \eqref{Def_M_hat_M}.   For any initial datum $R_0=(r_0, \theta_0, z_0)$ in $\mathcal{M}$,
we denote by $R_t^{R_0}=(r_t, \theta_t, z_t)$, $t \in [0,T]$, the corresponding unique strong solution to the original system \eqref{Eq_3DHopf_polar_base}.
The unique strong solution $(\hat{r}_t, \hat{\theta}_t)$ to the reduced system \eqref{Eq_2DHopf_polar_base} emanating from the initial datum $\hat{R}_0 = (r_0,\theta_0)$ in $\hat{\mathcal{M}}$ is denoted by $\hat{R}_t^{R_0}.$	

Let $q>0$ be given by
\be \label{Eq_def_q_girs}
q = \frac{1}{\epsilon \gamma}\left(1 + \frac{1}{\sigma^2}\right),
\ee
and let $g$ be the function defined in \eqref{Eq_g}. Consider the transformed system \eqref{Eq_3DHopf_polar_base_transf}. We denote by $\tilde{R}_t^{R_0}$ its unique strong solution emanating from $R_0$.

Then, the stochastic exponential $\mathfrak{D}_T$ at time $T$ on $(\Omega, \mathcal{B}, \mathbb{P})$ given by
	\be \label{Eq_D_girs}
	\mathfrak{D}_T = \exp\left(-\int_0^T g(\tilde{r}_s, \hat{r}_s)\d W^r_s - \frac{1}{2}\int_0^T \betrag{g(\tilde{r}_s, \hat{r}_s)}^2\d s \right),
	\ee
	is integrable with respect to $\mathbb{P}$ and satisfies 
	\be \label{Eq_D_1}
	\mathbb{E}_{\mathbb{P}} [\mathfrak{D}_T] = 1.
	\ee
	Hence, $\mathfrak{D}_T$ induces a probability measure $\tilde{\mathbb{P}} = \tilde{\mathbb{P}}_T$ on $(\Omega, \mathcal{B})$ given by
	\be \label{Eq_P_tilde_girs}
	\tilde{\mathbb{P}}(\d \omega) = \mathfrak{D}_T(\omega)\mathbb{P}(\d \omega).
\ee
Furthermore, $\mathbb{P}$ and $\tilde{\mathbb{P}}$  are equivalent and
	the stochastic process $(\tilde{W}_t)_{t \in [0,T]} = (\tilde{W}^r_t, \tilde{W}^{\theta}_t, \tilde{W}^3_t)_{t\in [0,T]}$ defined by
	 \bea \label{Eq_girsanov_tilde_W}
	 \begin{pmatrix} 
	\tilde{W}^r_t \\
	\tilde{W}^{\theta}_t \\
	\tilde{W}^3_t
	\end{pmatrix} = 
	\int_0^t G(\tilde{r}_s, \hat{r}_s)\d s +
	\begin{pmatrix} 
	W^r_t\\
	W^{\theta}_t \\
	W^3_t
	\end{pmatrix},
	\eea
	where 
	\bes
	G(\tilde{r}_t, \hat{r}_t)  = \begin{pmatrix} 
		g(\tilde{r}_t, \hat{r}_t) \\
		0 \\
		0
	\end{pmatrix}, \,\,\, t \in [0,T], 
	\ees
	is an $\mathbb{R}^3$-valued Brownian motion.  
	
Finally it holds that
\be \label{Eq_trans_prob_girs}
\tilde{\mathbb{P}}\Big(\tilde{R}_t^{R_0} \in \Gamma\Big) = \mathbb{P}\Big(R_t^{R_0} \in \Gamma\Big),\,\,\, \Gamma \in \mathcal{B}( \mathcal{M}), \; R_0 \in \mathcal{M}, \,\,\,t\in [0,T].
\ee
\el

\begin{proof} 
	{\bf{Step 1:} Proof of \eqref{Eq_D_1}.} 
	First, note that the stochastic processes $(\tilde{R}_t^{R_0})_{t\in [0,T]}$ and $(\hat{R}_t^{\hat{R}_0})_{t\in [0,T]}$ are continuous and adapted to the filtration $\{\mathcal{B}_t\}_{t\in [0,T]}$ in \eqref{Eq_filtration}, as they are unique strong solutions to systems \eqref{Eq_3DHopf_polar_base_transf} and \eqref{Eq_2DHopf_polar_base}, respectively. Therefore, $(\tilde{R}_t^{R_0})_{t\in [0,T]}$ and $(\hat{R}_t^{\hat{R}_0})_{t\in [0,T]}$ are predictable.
	According to \cite[Proposition 10.17, p. 294]{daprato_zabczyk_2014},
for \eqref{Eq_D_1} {\mkr to hold, it is sufficient to prove that there exists $\delta>0$ such that} 
	\be \label{Eq_novikov_rep_2}
	\sup_{t\in[0,T]} \mathbb{E}_{\mathbb{P}}\left[\exp\left(\delta |g(\tilde{r}_t, \hat{r}_t)|^2\right)\right] < \infty.
	\ee
	To this end, Lemma \ref{Lemma_moment_bounds} in Appendix \ref{Appendix_lemmas} provides us with the following exponential moment bound. For the parameters
	\bes
	\alpha = \frac{1}{\sigma^2 \epsilon \gamma},\,\,\,\beta = \frac{1}{\sigma^2}\,\,\,
	\text{and}\,\,\,\eta = \alpha,
	\ees
	we introduce the function $V$ given by
	\bea \label{Eq_V_girs_appl}
	V:\,& \mathbb{R}_+\times \mathbb{R} \times \mathbb{R}_+ \longrightarrow \mathbb{R},\\
	&(\tilde{r}, \tilde{z}, \hat{r}) \mapsto  \alpha \tilde{r}^2 + \beta \tilde{z}^2 + \eta \hat{r}^2.
	\eea
	Furthermore, let
	\bes
	Q(\tilde{r}, \tilde{z}, \hat{r}) = \exp\left(V(\tilde{r}, \tilde{z}, \hat{r})\right).
	\ees
	Lemma \ref{Lemma_moment_bounds} then asserts that for all $T>0$ there exists a constant $C_T>0$ such that 
	\be \label{Eq_to_show_B1_1}
	\mathbb{E}_{\mathbb{P}}\left[Q(\tilde{r}_t, \tilde{z}_t, \hat{r}_t) \right] \leq C_T,\,\,\,\text{for all}\,\,\,t \in [0,T].
	\ee
	By choosing 
	\bes
	\delta = \left(\frac{\sigma}{q+\lambda}\right)^2 \frac{\alpha}{2},
	\ees
	where $q$ is as in \eqref{Eq_def_q_girs},
	we obtain for all $t$ in $[0,T]$ that
	\bea \label{Eq_exp_mom_b_1}
	\mathbb{E}_{\mathbb{P}}\left[\exp\left(\delta |g(\tilde{r}_t,\hat{r}_t)|^2\right)\right] &= \mathbb{E}_{\mathbb{P}}\left[\exp\left(\frac{\alpha}{2}\hat{r}_t^2 -\alpha \tilde{r}_t\hat{r}_t + \frac{\alpha}{2} \tilde{r}_t^2\right)\right] \\
	&\leq \mathbb{E}_{\mathbb{P}}\left[\exp\left(\frac{\alpha}{2}\hat{r}_t^2 +\alpha |\tilde{r}_t||\hat{r}_t| + \frac{\alpha}{2} \tilde{r}_t^2\right)\right] \\
	&\leq \mathbb{E}_{\mathbb{P}}\left[\exp\left(\alpha\hat{r}_t^2 + \alpha \tilde{r}_t^2\right)\right] &(\text{due to Young's inequality})  \\
	&\leq \mathbb{E}_{\mathbb{P}}\left[\exp\left(\alpha\hat{r}_t^2 +\beta \tilde{z}_t^2 + \alpha \tilde{r}_t^2\right)\right] \\
	&=\mathbb{E}_{\mathbb{P}}\left[Q(\tilde{r}_t, \tilde{z}_t, \hat{r}_t)\right] \\
	&\leq C_T.
	\eea
	From \eqref{Eq_exp_mom_b_1} it follows that
	\bes
	\sup_{t \in [0,T]} \mathbb{E}_{\mathbb{P}}\left[\exp\left(\delta |g(\tilde{r}_t,\hat{r}_t)|^2\right)\right] \leq C_T,
	\ees
i.e.~\eqref{Eq_novikov_rep_2}.\\
	
		{\bf{Step 2:} Proof of \eqref{Eq_trans_prob_girs}.}
	Recall Remark \ref{Remark_ex_str_sol}, which ensures the existence and uniqueness of strong solutions to the original system \eqref{Eq_3DHopf_polar_base} on the probability space $(\Omega, \mathcal{B}, \mathbb{P})$. This property, together with \cite[Theorem 1.1, p. 163]{watanabe_ikeda} yields that for every initial value $R_0$ in $\mathcal{M}$ there exists a solution in the sense of \cite[Definition 1.2, p. 160-161]{watanabe_ikeda}. $(R_t^{R_0})_{t\geq0}$ to system \eqref{Eq_3DHopf_polar_base}, which is pathwise unique in the sense of \cite[Definition 1.5, p. 162]{watanabe_ikeda}.

	According to \cite[Corollary, p.~166]{watanabe_ikeda}, the solution $(R_t^{R_0})_{t\geq0}$ is also \textit{unique} in the following sense. If $(\bar{R}_t^{R_0})_{t\geq 0}$ is another solution for system \eqref{Eq_3DHopf_polar_base} to the same initial value $R_0$ on a probability space $(\Omega, \mathcal{B}, \bar{\mathbb{P}})$, then the law of the processes $(R_t^{R_0})_{t\geq0}$ and $(\bar{R}_t^{R_0})_{t\geq 0}$ coincide. This implies in particular that
	\bes \label{Eq_trans_prob_girs_hat}
	\bar{\mathbb{P}}\Big(\bar{R}_t^{R_0} \in \Gamma\Big) = \mathbb{P}\Big(R_t^{R_0} \in \Gamma\Big),\,\,\, \Gamma \in \mathcal{B}( \mathcal{M}), \; R_0 \in \mathcal{M}, \,\,\,t\geq 0.
	\ees
	We choose a $T>0$ and consider solutions only on the interval $[0,T]$. In the first step of this proof we have shown that the process \eqref{Eq_girsanov_tilde_W} is a Brownian motion on $(\Omega, \mathcal{B}, \tilde{\mathbb{P}})$, where $\tilde{\mathbb{P}}$ is defined as in \eqref{Eq_P_tilde_girs}. 
	From this it follows that the unique strong solution $(\tilde{R}_t^{R_0})_{t\in [0,T]}$ to the transformed  system \eqref{Eq_3DHopf_polar_base_transf} is also a solution in the sense of \cite[Definition 1.2, p. 160-161]{watanabe_ikeda} for the original system \eqref{Eq_3DHopf_polar_base} when viewed on the probability space $(\Omega, \mathcal{B}, \tilde{\mathbb{P}})$. Note that this implies that the law of $(R_t^{R_0})_{t\in[0,T]}$ with respect to $\mathbb{P}$ and the law of $(\tilde{R}_t^{R_0})_{t\in[0,T]}$ with respect to $\tilde{\mathbb{P}}$ coincide. In other words it holds that
	\bes \label{Eq_trans_prob_girs_conclusion}
	\tilde{\mathbb{P}}\Big(\tilde{R}_t^{R_0} \in \Gamma\Big) = \mathbb{P}\Big(R_t^{R_0} \in \Gamma\Big),\,\,\, \Gamma \in \mathcal{B}( \mathcal{M}), \; R_0 \in \mathcal{M}, \,\,\,t\in [0,T],
	\ees
	which is \eqref{Eq_trans_prob_girs}.
\end{proof}

\br \label{Remark_ex_inv_m_1}
For the original system \eqref{Eq_3DHopf_polar_base} there exists a unique ergodic invariant measure $\mu$ in $Pr(\mathcal{M})$ with smooth density with respect to the Lebesgue measure on $\mathcal{M}$; see Lemma \ref{Lemma_invm_ex_org} in Appendix \ref{App_ex_uniq_org}. Similarly, due to Lemma \ref{Lemma_invm_exuniq_red}, there is also a unique ergodic invariant measure $\nu$ in $Pr(\mathcal{\hat{M}})$ for the reduced system \eqref{Eq_2DHopf_polar_base}. Moreover, $\nu$ has a smooth Lebesgue-density on $\mathcal{\hat{M}}$.
\er

\subsection{Main results: Statements and interpretations}
\label{Subsec_statement_thm_1}

\bt \label{Thm_hopf_1}
Let $\mu(\d r,\d \theta, \d z)$ be the ergodic invariant measure for the original system \eqref{Eq_3DHopf_polar_base} and denote by $\mu_{r,z}$ and $\mu_r$ its marginals in the $(r,z)$-plane and along {\mkr the radial direction}, respectively. Let $\nu(\d r,\d \theta)$ be the ergodic invariant measure for the reduced system \eqref{Eq_2DHopf_polar_base} and $\nu_r$ its marginal along the radial direction.
Let
\be \label{Eq_r_*}
r_* = \left(\frac{\lambda}{2\gamma} + \frac{1}{2}\sqrt{\left(\frac{\lambda}{\gamma}\right)^2 + \frac{2\sigma^2}{\gamma}}\right)^\frac{1}{2},
\ee
{\mkr denote} the unique positive root of 
\bes
\lambda r - \gamma r^3 + \frac{\sigma^2}{2r} =0.
\ees
Recall the definitions of $\mathcal{F}$ in \eqref{Eq_test_fcts_1} and of the {\mkr Wasserstein} metric $\mathcal{D}_{\mathcal{F}}$ {\mkr given by} \eqref{Eq_IPM}.
Then there exist constants $c>0$ and $C(r_*, \mu_r,\nu_r)>0$ such that 
\be \label{Eq_main_result_3}
\boxed{
D_{\mathcal{F}}(\mu_r, \nu_r) \leq C(r_*, \mu_r, \nu_r) + c_\epsilon \left(\int_ {\mathbb{R}_+ \times \mathbb{R}}|z - r^2|^4\,\mu_{r, z}(\d r, \d z)\right)^{\frac{1}{4}},
}
\ee
where
\be \label{Eq_closeness_to_crit_3}
C(r_*, \mu_r, \nu_r) = 2r_* + \int_{r_*}^{\infty} r \nu_r(\d r) + \int_{r_*}^{\infty} r \mu_r(\d r),
\ee
and
\be\label{Eq_c_epsilon}
{\mkr c_\epsilon = \epsilon \left(\frac{\gamma^2 \sigma^2}{1+\sigma^2}\right)\left(\int_0^{\infty}r^4\mu_r(\d r)\right)^{\frac{1}{4}}}.
\ee

\et
\needspace{1\baselineskip}
\br\label{Rmk_epsilon_tails}
\hspace*{2em}  \vspace*{-0.4em}
\bi
\item[(i)] A key quantity of Theorem \ref{Thm_hopf_1} is $r_*$ in \eqref{Eq_r_*}. It coincides, for $\sigma = 0$, with the radius of the limit cycle {\mkr $r_{\textrm{det}}$ given by \eqref{Eq_r_det} of the reduced system \eqref{Eq_2DHopf_polar_base}}. Moreover, for any chosen $\gamma$, we observe that the quantity $r_*$ increases, if either $\lambda$ controlling the repulsiveness of the origin, or the noise {\mkr intensity $\sigma$ increases}.  In all the cases, we have $r_\ast >r_{\textrm{det}}$.

\item[(ii)] We call the quantity 
\be \label{Eq_def_defect}
Q(h)=\left(\int_ {\mathbb{R}_+ \times \mathbb{R}}|z - r^2|^4\,\mu_{r, z}(\d r, \d z)\right)^{\frac{1}{4}},
\ee
which appears on the RHS of inequality \eqref{Eq_main_result_3} the \textit{parameterization defect} or {\mkr simply} \textit{defect}. It measures the error incurred by parameterizing the fast variable $z$ in \eqref{Eq_3DHopf_polar_base} by the slow manifold \eqref{Eq_slow_mnf_polar}.

\item[(iii)]  If $\sigma$ is small enough, one expects the dynamics of the original system to meander around the limit cycle  $\Gamma=\mathcal{C}\times \{r_{\textrm{det}}^2\}$ (see discussion right after \eqref{Eq_r_det}), so as $r_\ast >r_{\textrm{det}}$, it is expected that $\int_{r_\ast} r \mu_r(\d r)$ is small, since most of the "mass" is expected to be located near this limit cycle. A similar observation holds regarding $\int_{r_*}^{\infty} r \nu_r(\d r)$, when considering the reduced system (replacing $\Gamma$ by $\mathcal{C}$) and therefore for $\sigma$ sufficiently small. One is thus left with $D_{\mathcal{F}}$ controlled essentially by let us say $3 r_\ast + c_\epsilon Q(h)$. Thus $D_{\mathcal{F}}$ is small if $r_\ast$ is small and $c_\epsilon$ is sufficiently small. 
Given the expression of $r_\ast$ in \eqref{Eq_r_*}, $r_\ast$ is ensured to be small if
\bea\label{Cond_Dsmall}
&\frac{\lambda}{\gamma} <\eta \; \; \mbox{ and }\\
&\frac{2 \sigma^2}{\gamma}<\eta^3,
\eea
with $\eta$ sufficiently small. 

As a consequence, $D_{\mathcal{F}}$ is guaranteed to be small for situations {\bf close to criticality} ($\lambda$ small) for which the nonlinear dissipative effects are sufficiently important ($\gamma$ large), the noise
intensity $\sigma$ is sufficiently small, while the {\bf time-scale separation is sufficiently strong} ($\epsilon$ small).  
The numerical results presented in Section \ref{Sec_numerics_1} demonstrate that the slow manifold reduced system \eqref{Eq_2DHopf} can still provide efficient reductions even beyond these restrictive conditions.
\ei  
\er

The following theorem provides an alternative bound of the form \eqref{Eq_main_result_3}. The main difference is that the additive constant in the RHS of inequality \eqref{Eq_main_result_2} below is no longer dependent on $\mu_r$.
\bt \label{Thm_hopf_2}
	Let $\mu(\d r,\d \theta, \d z)$ be the ergodic invariant measure for the original system \eqref{Eq_3DHopf_polar_base} and denote by $\mu_{r,z}$ and $\mu_r$ its marginals in the $(r,z)$-plane and along $r$, respectively. Let $\nu(\d r,\d \theta)$ be the ergodic invariant measure for the reduced system \eqref{Eq_2DHopf_polar_base} and $\nu_r$ its marginal along $r$.
	Furthermore, let
	\be \label{Eq_r_det_thm}
	r_{\mathrm{det}} = \sqrt{\frac{\lambda}{\gamma}},
	\ee
	which is the unique positive root of 
	\bes
	\lambda r - \gamma r^3 =0.
	\ees
	Recall the definition of $\mathcal{F}$ in \eqref{Eq_test_fcts_1}. 
	Then there exist constants $c>0$ and $C(\lambda, \gamma, \sigma, \nu_r)>0$ such that
	\be \label{Eq_main_result_2}
	D_{\mathcal{F}}(\mu_r, \nu_r) \leq C(\lambda, \gamma, \sigma, \nu_r) + c \left(\int_ {\mathbb{R}_+ \times \mathbb{R}}|z - r^2|^4\,\mu_{r, z}(\d r, \d z)\right)^{\frac{1}{4}},
	\ee
	where 
	\be \label{Eq_closeness_to_crit_2}
	C(\lambda, \gamma, \sigma, \nu_r) = \sqrt{2\left(r_{\mathrm{det}}^2 + \frac{\sigma^2}{\lambda}\right)}+\int_0^{\infty} r \nu_r(\d r),
	\ee
	and
	\bes
	c =\frac{q+2\lambda}{qr_{\mathrm{det}}^2}\left(\int_0^{\infty}r^4\mu_r(\d r)\right)^{\frac{1}{4}},
	\ees
	with $q$ as in \eqref{Eq_def_q_girs}.
\et

\subsection{Proof of Theorem \ref{Thm_hopf_1}}
\label{Sec_proof_thm_1}
\begin{proof} 
{\bf{Step 1:}} In this preparatory first step we introduce the main objects needed to carry out the proof.
For this purpose, we recall the original and reduced system \eqref{Eq_3DHopf_polar_base} and \eqref{Eq_2DHopf_polar_base}, respectively, and introduce their Markov semigroups. 

First of all, the original system \eqref{Eq_3DHopf_polar_base} is given---using cylindrical coordinates---by
\begin{subnumcases}{\label{Eq_3DHopf_polar_base_proof}}
\d r = \left[\lambda r -\gamma r z + \frac{\sigma^2}{2 r}\right]\d t + \sigma \d W^r_t \label{Eq_r_proof}\\
\d \theta = f \d t +  \frac{\sigma}{r} \d W^{\theta}_t  \mod(2\pi)\label{Eq_theta_proof} \\
\d z = -\frac{1}{\epsilon}\left[z -r^2\right]\d t + \frac{\sigma}{\sqrt{\epsilon}} \d W^3_t.
\end{subnumcases}

The unique strong solution to system \eqref{Eq_3DHopf_polar_base_proof} for an initial value $R_0$ in  $\mathcal{M}$ is denoted by $\big(R_t^{R_0}\big)_{t\geq 0}$, where
\bes
R_t^{R_0} = (r_t, \theta_t, z_t), \text{ for all } t \geq 0. 
\ees
The Markov semigroup
$(P_t)_{t \geq 0}$ associated with system \eqref{Eq_3DHopf_polar_base_proof} is defined as
\be \label{Eq_P_proof}
[P_t\psi](R_0) = \mathbb{E}_{\mathbb{P}}\big[\psi\big(R_t^{R_0}\big)\big], \,\,\, \psi \in C_b(\mathcal{M}), \,\,\, R_0 \in \mathcal{M},\,\,\, t \geq 0,
\ee
where $C_b(\mathcal{M})$ denotes the space of real-valued, bounded and continuous functions on $\mathcal{M}$.
The existence of a unique ergodic invariant measure $\mu$ for $P_t$ with smooth Lebesgue-density on $\mathcal{M}$ (see Remark \ref{Remark_ex_inv_m_1}), together with \cite[Theorem 5.8, p. 73]{daprato_inf_dim} implies that $P_t$ can be extended to a strongly continuous semigroup on $ L^2_{\mu}(\mathcal{M})$. We denote by $L^2_{\mu}(\mathcal{M})$ the space of all real-valued measurable functions on $\mathcal{M}$, which are square integrable {\mkr with respect to} $\mu$.
As discussed in Section \ref{Sec_main_idea}, the slow manifold \eqref{Eq_slow_mnf_polar} leads to the following slow manifold reduced system:
\begin{subnumcases}{\label{Eq_2DHopf_polar_base_proof}}
\d \hat{r} = \left[\lambda \hat{r} -\gamma \hat{r}^3 + \frac{\sigma^2}{2\hat{r}}\right]\d t  + \sigma \d W^{r}_t \label{Eq_r_hat_proof}\\
\d \hat{\theta} = f \d t +  \frac{\sigma}{\hat{r}} \d W^{\theta}_t \mod(2\pi) \label{Eq_theta_hat_proof}. 
\end{subnumcases}
Crucially, we require that  Eqns.~\eqref{Eq_r_hat_proof}-\eqref{Eq_theta_hat_proof} be driven by the same realizations of the Brownian motions $W^{r}_t$ and $ W^{\theta}_t$ as those driving Eq.~\eqref{Eq_3DHopf_polar_base_proof}. 


Furthermore, given an initial value $\hat{R}_0$ in $\hat{\mathcal{M}}$, the unique strong solution to system \eqref{Eq_2DHopf_polar_base_proof} is $\big(\hat{R}_t^{\hat{R}_0}\big)_{t\geq 0}$. 
Let $(Q_t)_{t\geq 0}$ be the Markov semigroup of system \eqref{Eq_2DHopf_polar_base_proof}, where 
\be \label{Eq_Q_proof}
[Q_t\psi](\hat{R}_0) = \mathbb{E}_{\mathbb{P}}\big[\psi\big(\hat{R}_t^{\hat{R}_0}\big)\big], \,\,\, \psi \in C_b(\hat{\mathcal{M}}), \,\,\, \hat{R}_0 \in \hat{\mathcal{M}},\,\,\, t \geq 0.
\ee

As pointed out in Remark \ref{Remark_ex_inv_m_1}, the existence of a unique ergodic invariant measure $\nu$ for $Q_t$ with smooth density with respect to the Lebesgue measure on $\hat{\mathcal{M}}$, is here ensured. Moreover, thanks to \cite[Theorem 5.8, p. 73]{daprato_inf_dim}, $Q_t$ is defined for all $\psi$ in $L^2_{\nu}(\hat{\mathcal{M}})$ and all $t\geq0$.

Let the parameter $q$ be as defined in \eqref{Eq_def_q_girs} and the function $g$ be given by \eqref{Eq_g}.
Given the radial component $(\hat{r}_t)_{t\geq0}$ of the unique strong solution $(\hat{R}_t^{\hat{R}_0})_{t \geq 0}$ to system \eqref{Eq_2DHopf_polar_base_proof} {\mkr (emanating from $\hat{R}_0 = (r_0,\theta_0)$ in $\hat{\mathcal{M}}$)},
we consider the following transformation of the original system \eqref{Eq_3DHopf_polar_base_proof}:
\begin{subnumcases}{\label{Eq_3DHopf_polar_base_transf_proof}}
\d \tilde{r} = \left[\lambda \tilde{r} -\gamma \tilde{r} \tilde{z} + \frac{\sigma^2}{2 \tilde{r}}\right]\d t + \sigma g(\tilde{r}, \hat{r})\d t + \sigma \d W^r_t \label{Eq_r_transf_proof}\\
\d \tilde{\theta} = f \d t +  \frac{\sigma}{\tilde{r}} \d W^{\theta}_t \mod(2\pi) \\
\d \tilde{z} = -\frac{1}{\epsilon}\left[\tilde{z} -\tilde{r}^2\right]\d t + \frac{\sigma}{\sqrt{\epsilon}} \d W^3_t.
\end{subnumcases}
{\mkr The function $g$ acts thus as a coupling term between Eq.~\eqref{Eq_3DHopf_polar_base_transf_proof} and the reduced system \eqref{Eq_2DHopf_polar_base_proof}; the solution of the latter acting as a forcing on the former (through $g$) but not reciprocally.} The unique strong solution to this system {\mkr originating from} $R_0$ in $\mathcal{M}$ is denoted by $(\tilde{R}_t^{R_0})_{t\geq 0}$. 

Lemma \ref{Lemma_girsanov} ensures that for every $T_*>0$ there exists a probability measure $
\tilde{\mathbb{P}} = \tilde{\mathbb{P}}_{T_*}$ on $(\Omega, \mathcal{B})$, which is equivalent to $\mathbb{P}$ on the $\sigma$-algebra $\mathcal{B}_{T_*}$ such that 
\be \label{Eq_trans_prob_proof}
\tilde{\mathbb{P}}\Big(\tilde{R}_t^{R_0} \in \Gamma\Big) = \mathbb{P}\Big(R_t^{R_0} \in \Gamma\Big),\,\,\, \Gamma \in \mathcal{B}( \mathcal{M}), \; R_0 \in \mathcal{M}, \,\,\,t\in [0,T_*].
\ee
{\mkr Identity \eqref{Eq_trans_prob_proof} expresses the idea that} transition probabilities of the original system \eqref{Eq_3DHopf_polar_base_proof} are preserved in its transformed counterpart \eqref{Eq_3DHopf_polar_base_transf_proof} up to a final time $T_*$, {\mkr  as long as the transition probabilities are assessed under the new probability measure $\tilde{\mathbb{P}}$ for \eqref{Eq_3DHopf_polar_base_transf_proof}.}

We set
\be \label{Eq_P_T_*}
[\tilde{P}^{T_*}_t\psi](R_0) = \mathbb{E}_{\tilde{\mathbb{P}}}\big[\psi\big(\tilde{R}_t^{R_0}\big)\big], \,\,\, \psi \in C_b(\mathcal{M}), \,\,\, R_0 \in \mathcal{M},\,\,\, t \in [0,T_*].
\ee
Because of identity \eqref{Eq_trans_prob_proof} and the fact that $P_t$ is defined for all $\psi$ in $L^2_{\mu}(\mathcal{M})$, we can set
\be \label{Eq_proof_1}
 [\tilde{P}^{T_*}_t\psi](R_0) = [P_t\psi](R_0),\,\,\, \psi \in L^2_{\mu}(\mathcal{M}), \,\,\, R_0 \in \mathcal{M},\,\,\, t \in [0,T_*]. 
\ee
{\bf{Step 2:}} In a second step we introduce the marginal measures $\mu_r, \nu_r$ and $\mu_{r,z}$ appearing in e.g.~\eqref{Eq_main_result_1}. Moreover, averages associated with the Markov semigroups $P_t$ and $Q_t$ are defined.

Recall the projections $\Pi_V$ as defined in \eqref{Def_proj}, for $V=V_r$ and $V=V_{r,z}$.
The marginal measures $\mu_{r}$, $\nu_r$ and $\mu_{r,z}$ associated with the subspaces $V_r$ and $V_{r,z}$ are then given by 
\be \label{Eq_marginal_proof}
\mu_r = \Pi_{V_r}^* \mu, \,\, \nu_r = \Pi_{V_r}^* \nu \,\,\,\text{and}\,\,\, \mu_{r,z} = \Pi_{V_{r,z}}^* \mu.
\ee
For the Markov semigroup $P_t$, we introduce the average operator $M_T$, $T\geq0$, {\mkr given by}
\be \label{Eq_erg_1}
M_T \psi = \frac{1}{T} \int_0^T P_t\psi \d t, \,\,\,\text{for all }\psi \text{ in } L^2_{\mu}(\mathcal{M}) \text{ and } T\geq0.
\ee
Analogously, for the Markov semigroup $Q_t$ defined in \eqref{Eq_Q_proof}, {\mkr we introduce}  
\be \label{Eq_erg_2}
N_T \psi = \frac{1}{T} \int_0^T Q_t\psi \d t, \,\,\, \text{for all }\psi \text{ in } L^2_{\nu}(\hat{\mathcal{M}}) \text{ and } T\geq0.
\ee
{\mkr Finally}, we define the averages $\tilde{M}_T$ {\mkr of $\tilde{P}^{T_*}_t$ in \eqref{Eq_proof_1}, for $T\leq T_*$,  as}
\be \label{Eq_erg_3}
\tilde{M}_T \psi = \frac{1}{T}\int_0^T \tilde{P}_t^{T_*}\psi\d t, \,\,\, \text{for all }\psi \text{ in } L^2_{\mu}(\mathcal{M}) \text{ and } T\leq T_*.
\ee

Our goal is to compare the marginals $\mu_r$ and $\nu_r$ of the invariant measures $\mu$ and $\nu$, respectively. For this purpose, let $\varphi$ be a real-valued measurable function defined on $\mathbb{R}_+$ and $\psi_{\varphi}$ an observable given by
\bea \label{Eq_psi}
\psi_{\varphi}: \, &\mathcal{M} \longrightarrow \mathbb{R},\\
&(r, \theta, z) \mapsto  (\varphi \circ \Pi_{V_r})(r, \theta, z).
\eea
Note that 
\be \label{Eq_proof_erg_2}
\int_{\mathcal{M}} \psi_{\varphi}(r, \theta, z)\,\mu(\d r, \d \theta, \d z) = \int_{\mathbb{R}_+} \varphi(r) \, \mu_r(\d r).
\ee
Furthermore, for $\varphi$ in $\mathcal{F}$ it holds that 
\bea \label{Eq_proof_erg_1_1}
\int_{\mathcal{M}} |\psi_{\varphi}(r, \theta, z)|^2\,\mu(\d r, \d \theta, \d z) &= \int_{\mathbb{R}_+} |\varphi(r)|^2 \, \mu_r(\d r) \\
&\leq \int_{\mathbb{R}_+} r^2 \, \mu_r(\d r) < \infty, &\text{(due to Lemma \ref{Lemma_invm_ex_org} in Appendix \ref{App_ex_uniq_org})}.
\eea
Hereafter, we assume $\varphi$ in $\mathcal{F}$, for which \eqref{Eq_proof_erg_1_1} shows that $\psi_{\varphi}$ is an element of $L^2_{\mu}(\mathcal{M})$.
Moreover, we have for $\mu$-almost all $R_0=(r_0,\theta_0,z_0)$ in $\mathcal{M}$ and all $t \geq 0$ that
\be \label{Eq_proof_erg_1}
 [P_t \psi_{\varphi}](R_0) = \mathbb{E}_{\mathbb{P}}[\psi_{\varphi}(R_t^{R_0})] = \mathbb{E}_{\mathbb{P}}[\varphi(r_t)].
\ee
The ergodicity of $\mu$ and identity \eqref{Eq_proof_erg_2} imply that 
\bea \label{Eq_proof_erg_3}
\lim_{T \rightarrow \infty}M_{T}\psi_{\varphi} = \int_{\mathcal{M}} \psi_{\varphi}(r, \theta, z)\,\mu(\d r, \d \theta, \d z) = \int_{\mathbb{R}_+} \varphi(r) \, \mu_r(\d r)\,\,\,\text{in}\,\,\,L^2_{\mu}(\mathcal{M}).
\eea
From \eqref{Eq_proof_erg_3} it follows that there exists a sequence 
\be \label{Eq_sub_Tn}
(T_n)_{n\in\mathbb{N}}\,\,\,\text{with}\,\,\, T_n\rightarrow \infty\,\,\,\text{as}\,\,\,n\rightarrow \infty,
\ee
and a measurable set $E_{\varphi}$ in $\mathcal{B}(\mathcal{M})$ with $\mu(E_{\varphi}) = 1$ such that for all $R_0$ in $E_{\varphi}$ it holds that
\be \label{Eq_girs_M_T_fin}
\lim_{n \rightarrow \infty}[M_{T_n}\psi_{\varphi}](R_0) = \lim_{n\rightarrow \infty}\frac{1}{T_n} \int_0^{T_n}[P_t \psi_{\varphi}](R_0)\d t =  \lim_{n \rightarrow \infty} \frac{1}{T_n} \int_0^{T_n} \mathbb{E}_{\mathbb{P}}[\varphi(r_t)]\d t = \int_{\mathbb{R}_+} \varphi(r) \, \mu_r(\d r),
\ee
where the second to last equality is due to \eqref{Eq_proof_erg_1}. 
Furthermore, let the function $f$ be defined by
\bea \label{Eq_f_proof}
f: \, & \mathbb{R}_+\times \mathbb{R} \longrightarrow \mathbb{R},\\
&(r,z) \mapsto  \left(\frac{\gamma}{q}\right)^2r^2 (z - r^2)^2.
\eea
Additionally, for $r_*$ given in \eqref{Eq_r_*} we introduce
\bea \label{Eq_f_aux}
\bar{\varphi}: \, & \mathbb{R}_+ \longrightarrow \mathbb{R},\\
&r \mapsto  r\mathds{1}_{\{r>r_*\}}.
\eea

{\mkr It is not difficult to observe that} $f$ is an element of $L^2_{\mu}(\mathcal{M})$; {\mkr see e.g.~Remark \ref{Remark_f_mom_bound}}.  Thus, there exists a subsequence $(T_{n_k})_{k\in\mathbb{N}}$ of $(T_n)_{n\in \mathbb{N}}$ and a measurable set $E_f$ in $\mathcal{B}(\mathcal{M})$ with $\mu(E_f) = 1$ such that for all $R_0$ in $E_f$ it holds that 
\be \label{Eq_girs_M_T_fin_f}
\lim_{k \rightarrow \infty}[M_{T_{n_k}}f](R_0) = \lim_{k\rightarrow \infty}\frac{1}{T_{n_k}} \int_0^{T_{n_k}}[P_t f](R_0)\d t  = \int_{\mathbb{R}_+\times \mathbb{R}} f(r,z) \, \mu_{r,z}(\d r, \d z).
\ee
Similarly, there is a subsequence of $(T_{n_k})_{k\in \mathbb{N}}$ also denoted by 
\be \label{Eq_sub_Tnk}
(T_{n_k})_{k\in\mathbb{N}}\,\,\,\text{with}\,\,\, T_{n_k}\rightarrow \infty\,\,\,\text{as}\,\,\,k\rightarrow \infty,
\ee
and a measurable set $E_{\bar{\varphi}}$ in $\mathcal{B}(\mathcal{M})$ with $\mu(E_{\bar{\varphi}})=1$ such that for all $R_0$ in $E_{\bar{\varphi}}$ we have that
\be \label{Eq_girs_M_T_fin_f_2}
\lim_{k \rightarrow \infty}[M_{T_{n_k}}\psi_{\bar{\varphi}}](R_0) = \lim_{k\rightarrow \infty}\frac{1}{T_{n_k}} \int_0^{T_{n_k}}[P_t \psi_{\bar{\varphi}}](R_0)\d t  = \int_{r_*}^{\infty} r \,\mu_{r}(\d r).
\ee
An analogous reasoning holds for the Markov semigroup $Q_t$ in \eqref{Eq_Q_proof}. 
More specifically, we define an observable $\hat{\psi}_{\varphi}$ for a real-valued and measurable function $\varphi$ on $\mathbb{R}_+$ by
\bea \label{Eq_psi_hat}
\hat{\psi}_{\varphi}: \, &\hat{\mathcal{M}} \longrightarrow \mathbb{R},\\
&(r, \theta) \mapsto  (\varphi \circ \Pi_{V_r})(r, \theta).
\eea
Observe that
\be \label{Eq_proof_erg_2_2}
\int_{\hat{\mathcal{M}}} \hat{\psi}_{\varphi}(r, \theta)\,\nu(\d r, \d \theta) = \int_{\mathbb{R}_+} \varphi(r) \, \nu_r(\d r),
\ee
and for $\varphi$ in $\mathcal{F}$
\bea \label{Eq_proof_erg_2_3}
\int_{\hat{\mathcal{M}}} |\hat{\psi}_{\varphi}(r, \theta)|^2\,\nu(\d r, \d \theta) &= \int_{\mathbb{R}_+} |\varphi(r)|^2 \, \nu_r(\d r) \\
&\leq \int_{\mathbb{R}_+} r^2 \, \nu_r(\d r) < \infty, & \text{(due to Lemma \ref{Lemma_invm_exuniq_red} in Appendix \ref{App_ex_uniq_org})}.
\eea
Let $\varphi$ in $\mathcal{F}$. From \eqref{Eq_proof_erg_2_3} we have that $\hat{\psi}_{\varphi}$ is an element of $L^2_{\nu}(\hat{\mathcal{M}})$.
For an initial value $\hat{R}_0 = (\hat{r}_0, \hat{\theta}_0)$ in $\hat{\mathcal{M}}$, 
it holds for all $t \geq 0$ that
\be \label{Eq_proof_erg_5}
[Q_t \hat{\psi}_{\varphi}](\hat{R}_0) = \mathbb{E}_{\mathbb{P}}[\hat{\psi}_{\varphi}(\hat{R}_t^{\hat{R}_0})] = \mathbb{E}_{\mathbb{P}}[\varphi(\hat{r}_t)].
\ee
The ergodicity of $\nu$ and identity \eqref{Eq_proof_erg_2_2} imply, analogously to \eqref{Eq_proof_erg_3}, that
\be \label{Eq_proof_erg_6}
\lim_{T\rightarrow \infty} N_T\hat{\psi}_{\varphi} =\int_{\hat{\mathcal{M}}} \hat{\psi}_{\varphi}(r, \theta)\,\nu(\d r, \d \theta) = \int_{\mathbb{R}_+} \varphi(r) \, \nu_r(\d r)\,\,\,\text{in}\,\,\,L^2_{\nu}(\hat{\mathcal{M}}),
\ee
and for the sequence $(T_{n_k})_{k\in \mathbb{N}}$ in \eqref{Eq_sub_Tnk}
\be \label{Eq_girs_M_T_fin_2}
\lim_{k\rightarrow \infty} N_{T_{n_k}}\hat{\psi}_{\varphi} =\int_{\hat{\mathcal{M}}} \hat{\psi}_{\varphi}(r, \theta)\,\nu(\d r, \d \theta) = \int_{\mathbb{R}_+} \varphi(r) \, \nu_r(\d r)\,\,\,\text{in}\,\,\,L^2_{\nu}(\hat{\mathcal{M}}).
\ee
Hence, there exists a subsequence of $(T_{n_k})_{k\in \mathbb{N}}$ also denoted by $(T_{n_k})_{k\in \mathbb{N}}$ and a measurable set $\hat{E}_{\varphi}$ in the Borel $\sigma$-algebra  $\mathcal{B}(\hat{\mathcal{M}})$ with $\nu(\hat{E}_{\varphi}) = 1$ such that for all $\hat{R}_0$ in $\hat{E}_{\varphi}$ it holds that
\be \label{Eq_girs_N_T_fin}
\lim_{k\rightarrow \infty} [N_{T_{n_k}}\hat{\psi}_{\varphi}](\hat{R}_0) = \lim_{k\rightarrow\infty}\frac{1}{T_{n_k}}\int_0^{T_{n_k}} [Q_t\hat{\psi}_{\varphi}](\hat{R}_0)\d t =
\lim_{k\rightarrow \infty} \frac{1}{T_{n_k}} \int_0^{T_{n_k}} \mathbb{E}_{\mathbb{P}}[\varphi(\hat{r}_t)] \d t =  \int_{\mathbb{R}_+} \varphi(r) \, \nu_r(\d r),
\ee
where the second to last equality is due to \eqref{Eq_proof_erg_5}. Recall the definition of $\bar{\varphi}$ in \eqref{Eq_f_aux}.
Analogously to \eqref{Eq_girs_M_T_fin_f_2} there exists a further subsequence of $(T_{n_k})_{k\in \mathbb{N}}$ also denoted by
\be \label{Eq_final_subseq}
 (T_{n_k})_{k\in \mathbb{N}},\,\,\,\text{with}\,\,\,T_{n_k} \rightarrow \infty\,\,\,\text{as}\,\,\,k\rightarrow \infty,
 \ee
 and a Borel set $\hat{E}_{\bar{\varphi}}$ in $\mathcal{B}(\hat{\mathcal{M}})$ with $\nu(\hat{E}_{\bar{\varphi}}) = 1$ such that for all $\hat{R}_0$ in $\hat{E}_{\bar{\varphi}}$ it holds that
\be \label{Eq_girs_M_T_fin_f_3}
\lim_{k \rightarrow \infty}[N_{T_{n_k}}\hat{\psi}_{\bar{\varphi}}](\hat{R}_0) = \lim_{k\rightarrow \infty}\frac{1}{T_{n_k}} \int_0^{T_{n_k}}[Q_t \hat{\psi}_{\bar{\varphi}}](\hat{R}_0)\d t  = \int_{r_*}^{\infty} r \,\nu_{r}(\d r).
\ee
\\
{\bf{Step 3:}} 
Our goal in this final step is to show that there exists a constant $C>0$, which is independent of $\varphi$ in $\mathcal{F}$, such that
\be \label{Eq_proof_thm_1_step_3_1}
\left|\int_{\R_+} \varphi(r) \, \mu_r(\d r) - \int_{\R_+} \varphi(r) \, \nu_r(\d r)\right| \leq C.
\ee
Inequality \eqref{Eq_main_result_3} then follows by taking the supremum over all $\varphi$ in $\mathcal{F}$ in \eqref{Eq_proof_thm_1_step_3_1} {\mkr as soon as} $C$  {\mkr is identified} with the RHS of inequality \eqref{Eq_main_result_3}. 

We start by choosing an arbitrary $\delta >0 $ and $\varphi$ in $\mathcal{F}$. 
Recall the definitions of $\psi_{\varphi}$ and $\hat{\psi}_{\varphi}$ in \eqref{Eq_psi} and \eqref{Eq_psi_hat}, respectively.

Let $R_0$ and $\hat{R}_0$ be elements of $E_{\varphi}\cap E_{\bar{\varphi}}\cap E_f$ and $\hat{E}_{\varphi} \cap \hat{E}_{\bar{\varphi}}$, respectively. Identities \eqref{Eq_girs_M_T_fin}, \eqref{Eq_girs_M_T_fin_f} and \eqref{Eq_girs_M_T_fin_f_2} imply that for the subsequence $(T_{n_k})_{k\in\mathbb{N}}$ in \eqref{Eq_final_subseq} there exist indices ${n_k}_1 = {n_k}_1(\delta,R_0, \varphi)$ in $\mathbb{N}$, ${n_k}_2 = {n_k}_2(\delta,R_0,\bar{\varphi})$ in $\mathbb{N}$ and ${n_k}_3 = {n_k}_3(\delta, R_0, f)$ in $\mathbb{N}$ such that 
\be  \label{Eq_limits_1}
\left|[M_{T_{n_k}} \psi_{\varphi}](R_0) - \int_{\mathbb{R}_+} \varphi(r) \, \mu_r(\d r)\right| \leq \delta, \,\,\,\text{for all}\,\,\,{n_k}\geq {n_k}_1,
\ee
\bes
\left|[M_{T_{n_k}} \psi_{\bar{\varphi}}](R_0) - \int_{r_*}^{\infty} r \, \mu_r(\d r)\right| \leq \delta, \,\,\,\text{for all}\,\,\,n_k\geq {n_k}_2,
\ees
and
\be \label{Eq_f_delta}
\left|[M_{T_{n_k}} f](R_0) - \int_{\mathbb{R}_+ \times \mathbb{R}}f(r,z) \, \mu_{r,z}(\d r, \d z)\right| \leq \delta, \,\,\,\text{for all}\,\,\,n_k\geq {n_k}_3.
\ee
Furthermore, due to equations \eqref{Eq_girs_N_T_fin} and \eqref{Eq_girs_M_T_fin_f_3}, there exist indices ${n_k}_4 = {n_k}_4(\delta,\hat{R}_0, \varphi)$ in $\mathbb{N}$ and ${n_k}_5 = {n_k}_5(\delta, \hat{R}_0, \bar{\varphi})$ in $\mathbb{N}$ for the subsequence $(T_{n_k})_{k\in \mathbb{N}}$ in \eqref{Eq_final_subseq} such that
\be \label{Eq_limits_5}
\left|[N_{T_{n_k}} \hat{\psi}_{\varphi}](\hat{R}_0) - \int_{\mathbb{R}_+} \varphi(r) \, \nu_r(\d r)\right| \leq \delta, \,\,\,\text{for all}\,\,\,n_k\geq {n_k}_4,
\ee
and
\be \label{Eq_limits_6}
\left|[N_{T_{n_k}} \hat{\psi}_{\bar{\varphi}}](\hat{R}_0) - \int_{r_*}^{\infty} r \, \nu_r(\d r)\right| \leq \delta, \,\,\,\text{for all}\,\,\,n_k\geq {n_k}_5.
\ee
We set 
\be \label{Eq_def_T_star}
n_* = \max\{{n_k}_1, {n_k}_2, {n_k}_3, {n_k}_4, {n_k}_5\}\,\,\,\text{and}\,\,\,T_* = T_{n_*}.
\ee
From \eqref{Eq_limits_1} and \eqref{Eq_limits_5} it follows that 
\be \label{Eq_proof_thm_1_3}
\left|\int_{\R_+} \varphi(r) \, \mu_r(\d r) - \int_{\R_+} \varphi(r) \, \nu_r(\d r)\right| \leq 2 \delta + \left|[M_{T_*}\psi_{\varphi}](R_0) - [N_{T_*}\hat{\psi}_{\varphi}](\hat{R}_0)\right|.
\ee
We need to provide an estimate for the second term in the RHS of \eqref{Eq_proof_thm_1_3}. For simplicity we assume that $r_0 = \hat{r}_0$.
An application of the triangular inequality gives
\be \label{Eq_proof_thm_1_1}
 \left|[M_{T_*}\psi_{\varphi}](R_0) - [N_{T_*}\hat{\psi}_{\varphi}](\hat{R}_0)\right| \leq
 \frac{1}{T_*} \int_0^{T_*} \left| \mathbb{E}_{\mathbb{P}}[\varphi(\hat{r}_t)] - \mathbb{E}_{\mathbb{P}} [\varphi(r_t)]\right|\d t.
\ee 
We now apply Girsanov's theorem shown to hold in Lemma \ref{Lemma_girsanov} for $T=T_*$, where $T_*$ is defined in \eqref{Eq_def_T_star}. More specifically, there exists a probability measure $\tilde{\mathbb{P}} = \tilde{\mathbb{P}}_{T_*}$ on $(\Omega, \mathcal{B})$, which is equivalent to $\mathbb{P}$ on $\mathcal{B}_{T_*}$ such that 

Applying Girsanov's theorem, as guaranteed by Lemma \ref{Lemma_girsanov} for $T=T_*$ (where $T_*$ is defined in \eqref{Eq_def_T_star}), we obtain the existence of a probability measure $\tilde{\mathbb{P}} = \tilde{\mathbb{P}}_{T_*}$ on $(\Omega, \mathcal{B})$, equivalent to $\mathbb{P}$ on $\mathcal{B}_{T_*}$, such that
\be \label{Eq_equality_of_E}
\mathbb{E}_{\mathbb{P}}[\varphi(r_t)]  = \mathbb{E}_{\tilde{\mathbb{P}}}[\varphi(\tilde{r}_t)],\,\,\,\text{for all}\,\,\,t\in [0,T_*].
\ee
This equality follows directly from \eqref{Eq_trans_prob_girs}.

Using identity \eqref{Eq_equality_of_E} in the RHS of inequality \eqref{Eq_proof_thm_1_1} yields 
\be \label{Eq_main_proof_girs_1}
\frac{1}{T_*} \int_0^{T_*} \left| \mathbb{E}_{\mathbb{P}}[\varphi(\hat{r}_t)] - \mathbb{E}_{\mathbb{P}} [\varphi(r_t)]\right|\d t = \frac{1}{T_*} \int_0^{T_*} \left| \mathbb{E}_{\mathbb{P}}[\varphi(\hat{r}_t)] - \mathbb{E}_{\tilde{\mathbb{P}}} [\varphi(\tilde{r}_t)]\right|\d t.
\ee
By applying Lemma \ref{Lemma_aux_ineq_2} to the RHS of \eqref{Eq_main_proof_girs_1} and leveraging the properties of $\varphi$ as defined in \eqref{Eq_test_fcts_1} and \eqref{Eq_test_fcts}, we obtain
\bea \frac{1}{T_*} \int_0^{T_*} \left| \mathbb{E}_{\mathbb{P}}[\varphi(\hat{r}_t)] - \mathbb{E}_{\tilde{\mathbb{P}}} [\varphi(\tilde{r}_t)]\right|\d t
 &\leq 2r_* + \frac{1}{T_*} \int_0^{T_*} \mathbb{E}_{\tilde{\mathbb{P}}}[|\hat{r}_t-\tilde{r}_t|]\d t \label{Eq_proof_thm_1_aux_1} \\
 & +\frac{1}{T_*} \int_0^{T_*}\mathbb{E}_{\mathbb{P}}[\hat{r}_t \mathds{1}_{\{\hat{r}_t > r_*\}}]\d t \\
 & +\frac{1}{T_*} \int_0^{T_*} \mathbb{E}_{\mathbb{P}}[r_t\mathds{1}_{\{r_t > r_*\}}]\d t.
\eea 
For the second term in the RHS of \eqref{Eq_proof_thm_1_aux_1} we have
\bea \label{Eq_proof_thm_1_dep_phi}
 \frac{1}{T_*} \int_0^{T_*} & \mathbb{E}_{\tilde{\mathbb{P}}}[|\hat{r}_t-\tilde{r}_t|]\d t 
 \leq \left(\frac{1}{T_*} \int_0^{T_*} \mathbb{E}_{\tilde{\mathbb{P}}}[|\hat{r}_t-\tilde{r}_t|^2]\d t\right)^{\frac{1}{2}} &&\text{(due to the Cauchy-Schwarz inequality)} \\
 & \leq \frac{\gamma}{\sqrt{q}} \left(\frac{1}{T_*} \int_0^{T_*} \left(\int_0^te^{-q(t-s)} \mathbb{E}_{\tilde{\mathbb{P}}}\left[ \tilde{r}^2_s(\tilde{z}_s-\tilde{r}^2_s)^2\right]\d s\right) \d t \right)^{\frac{1}{2}} &&\text{(due to Lemma \ref{Lemma_aux_ineq_1} in Appendix \ref{Appendix_aux_ineq})} \\
 & \leq \frac{\gamma}{q} \left(\frac{1}{T_*} \int_0^{T_*} \mathbb{E}_{\tilde{\mathbb{P}}}\left[ \tilde{r}^2_t(\tilde{z}_t-\tilde{r}^2_t)^2\right]\d t \right)^{\frac{1}{2}} &&\text{(due to Fubini's theorem)}\\
  &=  \left(\frac{1}{T_*} \int_0^{T_*} \mathbb{E}_{\mathbb{P}}\left[\left(\frac{\gamma}{q}\right)^2 r^2_t(z_t-r^2_t)^2\right]\d t \right)^{\frac{1}{2}}, &&\text{(due to \eqref{Eq_trans_prob_girs} in Lemma \ref{Lemma_girsanov})}
\eea
and therefore
\be \label{Eq_proof_thm_1_4}
 \frac{1}{T_*} \int_0^{T_*} \mathbb{E}_{\tilde{\mathbb{P}}}[|\hat{r}_t-\tilde{r}_t|]\d t \leq \left(\frac{1}{T_*}\int_0^{T_*} [P_tf](R_0)\d t \right)^{\frac{1}{2}},
\ee
with $f$ given by \eqref{Eq_f_proof}.
 Note that the RHS of \eqref{Eq_proof_thm_1_4} depends implicitly on $\varphi$ due to the definition of $T_*$ in \eqref{Eq_def_T_star}. Next we show that
 \bes
 \frac{1}{T_*}\int_0^{T_*} [P_tf](R_0)\d t \leq \delta + \int_{\mathbb{R}_+\times \mathbb{R}} f(r,z)\mu_{r,z}(\d r, \d z).
 \ees
From \eqref{Eq_f_delta} and the definition of $T_*$ in \eqref{Eq_def_T_star} it follows that
 \bea \label{Eq_proof_delta_f_conclusion}
 \frac{1}{T_*} \int_0^{T_*} [P_tf](R_0)\d t &= [M_{T_*}f](R_0) \\
  &\leq \left|[M_{T_*}f](R_0) - \int_{\mathbb{R}_+ \times \mathbb{R}}f(r,z) \, \mu_{r,z}(\d r, \d z)\right| + \int_{\mathbb{R}_+ \times \mathbb{R}}f(r,z) \, \mu_{r,z}(\d r, \d z) \\
 &\leq \delta +\int_{\mathbb{R}_+ \times \mathbb{R}}f(r,z) \, \mu_{r,z}(\d r, \d z).
 \eea
 Recall once more the definition of $f$ in \eqref{Eq_f_proof}. Thanks to \eqref{Eq_proof_delta_f_conclusion} we have in \eqref{Eq_proof_thm_1_4}
 \bea \label{Eq_proof_thm_1_2}
 \frac{1}{T_*} \int_0^{T_*} \mathbb{E}_{\tilde{\mathbb{P}}}[|\hat{r}_t-\tilde{r}_t|]\d t 
 & \leq\left(\delta + \int_{\mathbb{R}_+ \times \mathbb{R}} 
 \left(\frac{\gamma}{q}\right)^2 r^2(z-r^2)^2 \mu_{r,z}(\d r, \d z) \right)^{\frac{1}{2}}\\
 &\leq\left(\delta + \left(\frac{\gamma}{q}\right)^2 \left(\int_0^{\infty} r^4\mu_r(\d r)\right)^{\frac{1}{2}}\left(\int_{\mathbb{R}_+ \times \mathbb{R}} 
  \left|z-r^2\right|^4 \mu_{r,z}(\d r, \d z)\right)^{\frac{1}{2}} \right)^{\frac{1}{2}}.
\eea
An analogous reasoning for the third term in \eqref{Eq_proof_thm_1_aux_1} yields, due to \eqref{Eq_girs_M_T_fin_f_3},
\be \label{Eq_tail_r_hat}
\frac{1}{T_*} \int_0^{T_*}\mathbb{E}_{\mathbb{P}}[\hat{r}_t \mathds{1}_{\{\hat{r}_t > r_*\}}]\d t 
\leq \delta +  \int_{r_*}^{\infty} r \nu_r(\d r).
\ee
Similarly, because of \eqref{Eq_girs_M_T_fin_f_2}, we obtain for the fourth term in \eqref{Eq_proof_thm_1_aux_1} that
\be \label{Eq_tail_r}
\frac{1}{T_*} \int_0^{T_*} \mathbb{E}_{\mathbb{P}}[r_t\mathds{1}_{\{r_t > r_*\}}]\d t
\leq \delta + \int_{r_*}^{\infty} r \mu_r(\d r).
\ee
We insert \eqref{Eq_proof_thm_1_2}, \eqref{Eq_tail_r_hat} and \eqref{Eq_tail_r} into \eqref{Eq_proof_thm_1_aux_1}, which yields in \eqref{Eq_proof_thm_1_3} that
\bea \label{Eq_proof_final_1}
\left|\int_{\R_+} \varphi(r) \, \mu_r(\d r) - \int_{\R_+} \varphi(r) \, \nu_r(\d r)\right| 
&\leq 4\delta +2r_* +\int_{r_*}^{\infty} r \nu_r(\d r) +\int_{r_*}^{\infty} r \mu_r(\d r)\\
&+\left(\delta + \left(\frac{\gamma}{q}\right)^2 \left(\int_0^{\infty} r^4\mu_r(\d r)\right)^{\frac{1}{2}}\left(\int_{\mathbb{R}_+ \times \mathbb{R}} 
\left|z-r^2\right|^4 \mu_{r,z}(\d r, \d z)\right)^{\frac{1}{2}} \right)^{\frac{1}{2}}.
\eea
Inequality \eqref{Eq_proof_final_1} holds for all $\delta>0$ and $\varphi$ in $\mathcal{F}$. Furthermore, its RHS is independent of $\varphi$. Hence, we can take the supremum over all $\varphi$ in $\mathcal{F}$ and let $\delta$ tend to zero, which results in 
\bea \label{Eq_thm_1_ineq_1}
D_{\mathcal{F}}(\mu_r, \nu_r) &\leq 2r_* + \int_{r_*}^{\infty} r \nu_r(\d r) + \int_{r_*}^{\infty} r \mu_r(\d r)\\
&+ \left(\frac{\gamma}{q}\right)\left(\int_0^{\infty}r^4\mu_r(\d r)\right)^{\frac{1}{4}} \left(\int_{\mathbb{R}_+\times\mathbb{R}} |z-r^2|^4 \mu_{r,z}(\d r, \d z) \right)^{\frac{1}{4}}.
\eea
By setting
\bes
C(r_*, \mu_r, \nu_r) = 2r_* + \int_{r_*}^{\infty} r \nu_r(\d r) +\int_{r_*}^{\infty} r \mu_r(\d r),
\ees
and
\bes
c = \left(\frac{\gamma}{q}\right)\left(\int_0^{\infty}r^4\mu_r(\d r)\right)^{\frac{1}{4}},
\ees
where $q$ is defined in \eqref{Eq_def_q_girs},
we finally obtain
\bes
D_{\mathcal{F}}(\mu_r, \nu_r) \leq C(r_*, \mu_r, \nu_r) + c \left(\int_ {\mathbb{R}_+ \times \mathbb{R}}|z - r^2|^4\,\mu_{r, z}(\d r, \d z)\right)^{\frac{1}{4}}.
\ees
\end{proof}

\subsection{Proof of Theorem \ref{Thm_hopf_2}}
\label{Sec_proof_thm_3}
\begin{proof}
We use the same notations as introduced in Steps 1 and 2 of the proof of Theorem \ref{Thm_hopf_1}. In addition, the reasoning of Step 3 until identity \eqref{Eq_main_proof_girs_1} applies analogously.
	By using Lemma \ref{lemma_aux_ineq_new_lemma} instead of Lemma \ref{Lemma_aux_ineq_2} in \eqref{Eq_proof_thm_1_aux_1}, we obtain
	\bea \label{Eq_proof_thm_1_aux_2} \frac{1}{T_*} \int_0^{T_*} \left| \mathbb{E}_{\mathbb{P}}[\varphi(\hat{r}_t)] - \mathbb{E}_{\tilde{\mathbb{P}}} [\varphi(\tilde{r}_t)]\right|\d t
	&\leq \frac{1}{T_*}\int_0^{T_*} e^{-\frac{\lambda}{2}t}\hat{r}_0\d t  
	+ \frac{1}{T_*} \int_0^{T_*} \mathbb{E}_{\tilde{\mathbb{P}}}[|\hat{r}_t-\tilde{r}_t|]\d t \\ &+\frac{q+\lambda}{\sqrt{\lambda}} \frac{1}{T_*}\int_0^{T_*} \left(\int_0^t e^{-\lambda(t-s)} \mathbb{E}_{\tilde{\mathbb{P}}}[|\hat{r}_s - \tilde{r}_s|^2]\d s\right)^{\frac{1}{2}}\d t  \\
	&+\sqrt{2\left(r_{\mathrm{det}}^2 + \frac{\sigma^2}{\lambda}\right)} +\frac{1}{T_*} \int_0^{T_*}\mathbb{E}_{\mathbb{P}}[\hat{r}_t]\d t. \\
	\eea 
	The second term in the RHS of \eqref{Eq_proof_thm_1_aux_2} is treated as in \eqref{Eq_proof_thm_1_dep_phi}. In addition, we use the Cauchy-Schwarz inequality and Lemma \ref{Lemma_aux_ineq_1}, where we assume that $r_0 = \hat{r}_0$ for simplicity, to obtain for the third term in the RHS of \eqref{Eq_proof_thm_1_aux_2}
	\beas
	\frac{q+\lambda}{\sqrt{\lambda}} \frac{1}{T_*}\int_0^{T_*} \left(\int_0^t e^{-\lambda(t-s)} \mathbb{E}_{\tilde{\mathbb{P}}}[|\hat{r}_s - \tilde{r}_s|^2]\d s\right)^{\frac{1}{2}}\d t &\leq \frac{q+\lambda}{q}\frac{\gamma}{\lambda}\left(\frac{1}{T_*}\int_0^{T_*} \mathbb{E}_{\tilde{\mathbb{P}}}[\tilde{r}^2_t(\tilde{z}_t - \tilde{r}_t^2)^2] \d t\right)^{\frac{1}{2}}\\
	&=\frac{q+\lambda}{q}\frac{\gamma}{\lambda}\left(\frac{1}{T_*}\int_0^{T_*} \mathbb{E}_{\mathbb{P}}[r^2_t(z_t - r_t^2)^2] \d t\right)^{\frac{1}{2}},
	\eeas
	where the last line is due to the fact that transition probabilities are preserved by virtue of Girsanov's theorem (see identity \eqref{Eq_trans_prob_girs} in Lemma \ref{Lemma_girsanov}).
	The last term in the RHS of \eqref{Eq_proof_thm_1_aux_2} is treated as in \eqref{Eq_tail_r_hat}.
	Following the same reasoning as in lines \eqref{Eq_proof_thm_1_4} -\eqref{Eq_thm_1_ineq_1}, leads to 
	\bes
	D_{\mathcal{F}}(\mu_r, \nu_r) \leq C(\lambda, \gamma, \sigma, \nu_r)+  c\left(\int_{\mathbb{R}_+\times\mathbb{R}} |z-r^2|^4 \mu_{r,z}(\d r, \d z) \right)^{\frac{1}{4}},
	\ees
	where
	\bes
	C(\lambda, \gamma, \sigma, \nu_r) = \sqrt{2\left(r_{\mathrm{det}}^2 + \frac{\sigma^2}{\lambda}\right)} +\int_{0}^{\infty} r \nu_r(\d r),
	\ees
	and 
	\bes
	c = \frac{q+2\lambda}{q r_{\mathrm{det}}^2}\left(\int_0^{\infty}r^4\mu_r(\d r)\right)^{\frac{1}{4}},
	\ees
	which completes the proof.
\end{proof}

\subsection{Numerical results}
\label{Sec_numerics_1}
Theorems \ref{Thm_hopf_1} and \ref{Thm_hopf_2} provide useful insights regarding 
the usage of the slow manifold as a parameterizing manifold \cite{CLM19_closure}. The purpose of this section is to illustrate and discuss this point via numerical experiments. 
Let us set the parameters in Table \ref{Table_CaseI} and differentiate between Cases I and II.

\begin{table}[h] 
	\caption{Parameter regimes: Case I and Case II}
	\label{Table_CaseI}
	\centering
	\begin{tabular}{llllll}
		\toprule\noalign{\smallskip}
		& $\lambda$ & $f$ &  $\gamma$  & $\epsilon$ & $\sigma$\\ 
		\noalign{\smallskip}\hline\noalign{\smallskip}
		Case I  & $10^{-3}$  & $10^2$ & $5.6 \times 10^{-2}$ & $10^{-2}$ & $0.55$ \\ 
		Case II & $10^{-3}$ &  $10$ & $1$ & $10^{-2}$ & $0.2$\\  
		\noalign{\smallskip} \bottomrule 
	\end{tabular}
\end{table}

\begin{figure}[htbp]
	\begin{subfigure}{0.45\textwidth}
		\includegraphics[width=1\textwidth,height=0.7\textwidth]{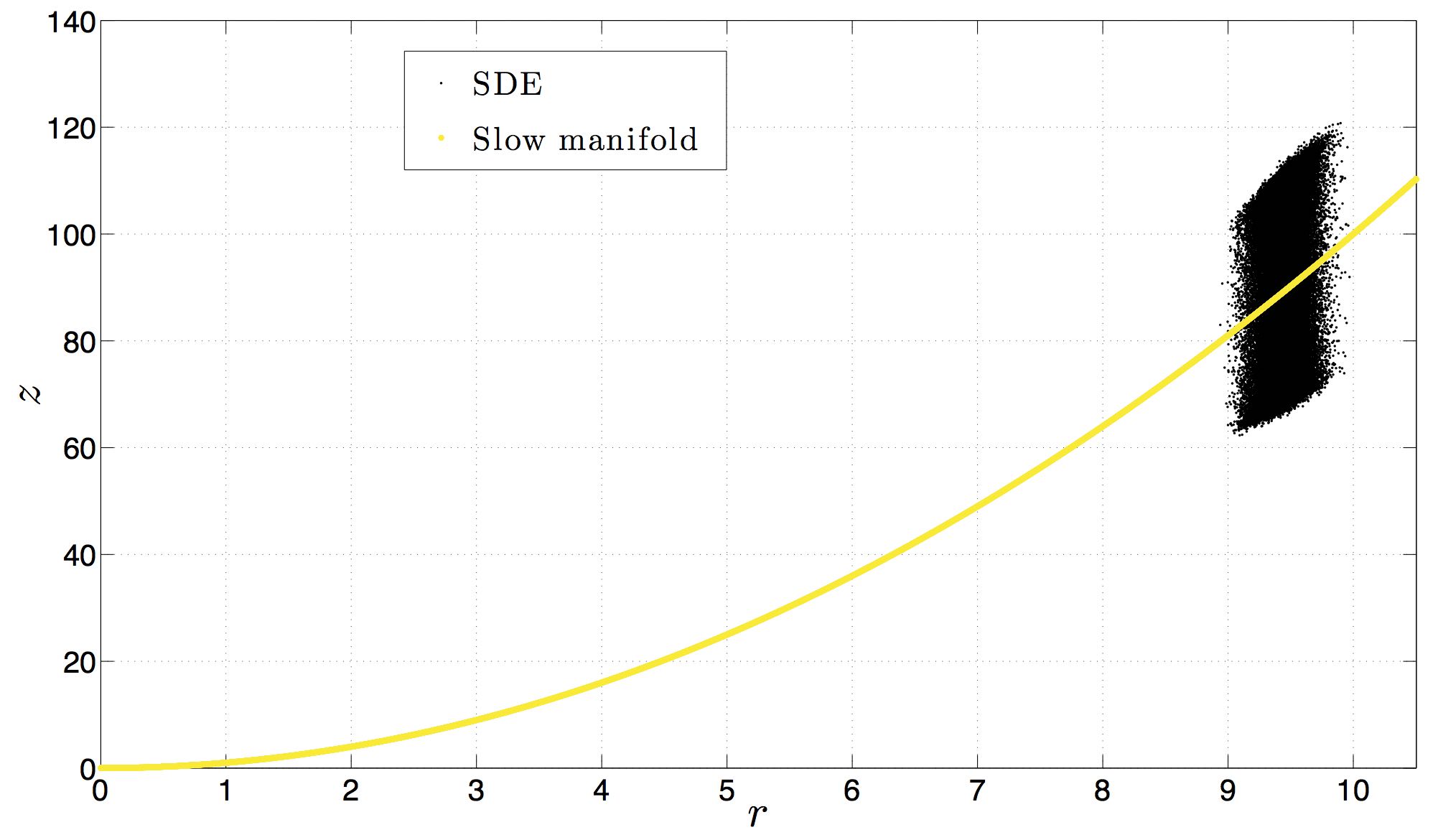}
		\caption{Case I}
	\end{subfigure}
	\begin{subfigure}{0.45\textwidth}
		\		\includegraphics[width=1\textwidth,height=0.7\textwidth]{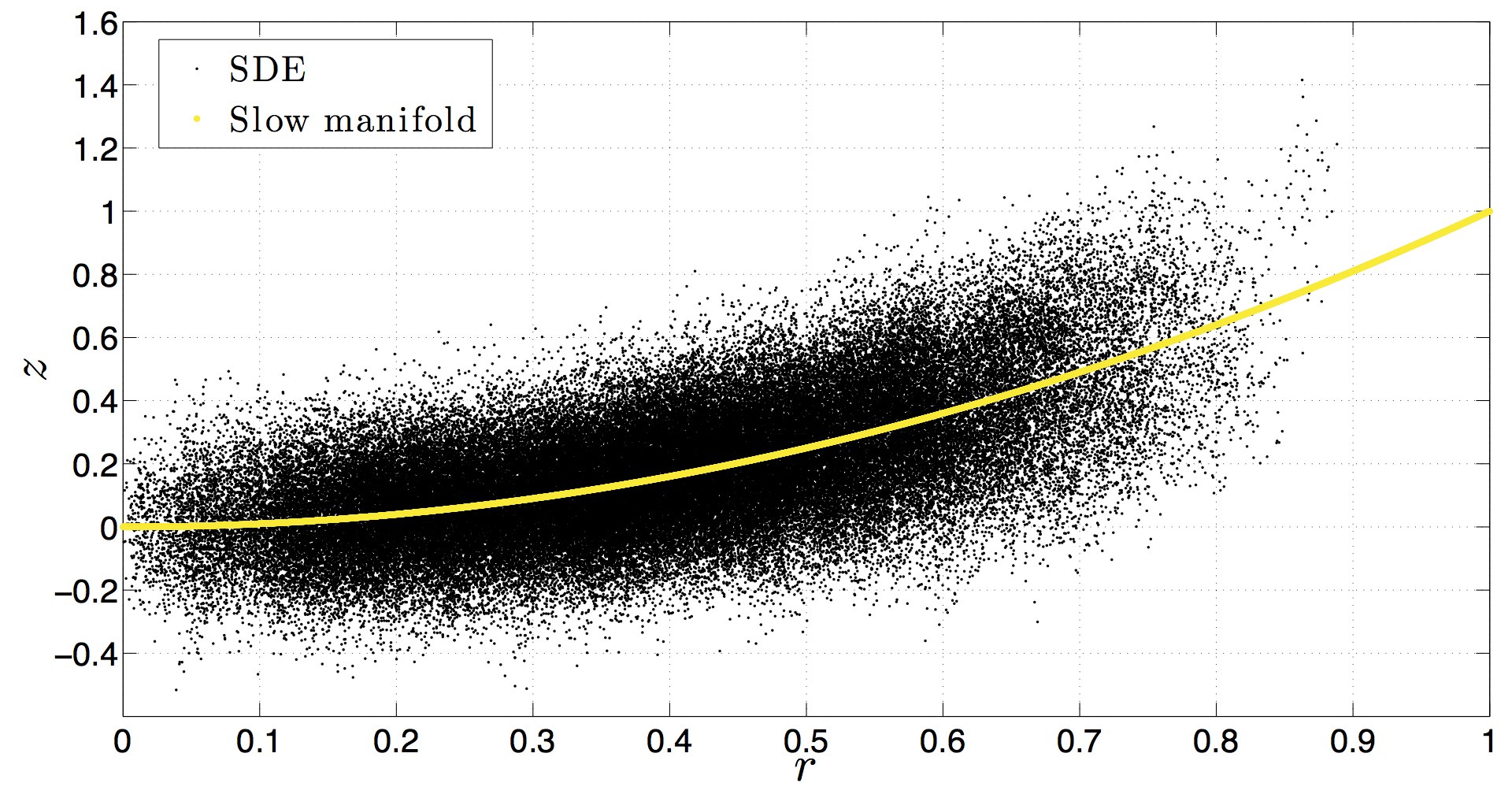}
		\caption{Case II}
	\end{subfigure}
	\caption{\footnotesize{\bf Scatter plots}. In Case I, the normalized parameterization defect for the slow manifold is given by $Q=1.789 \times 10^{-1}$, while in Case II, $Q=5.834 \times 10^{-1}$.} 
	\label{Visualization_Q}
\end{figure}
Figure \ref{Visualization_Q} shows that solutions to the original system \eqref{Eq_3DHopf_polar_base} evolve, after transient dynamics, mostly confined within an annulus  for the parameter regime corresponding to Case I. 
In Case II, the dynamics is less confined and meanders around the origin. In both cases   In both cases the deterministic slow manifold $h$ in \eqref{Eq_slow_mnf_polar} provides the average behavior of the fast variable $z$. 

We then compute the following normalized parameterization defect
\be \label{Eq_param_def_cont}
Q = \frac{\int_{\R_+ \times \R} \betrag{z-r^2}^4\mu_{r,z}(\d r, \d z)}{\int_{\R} \betrag{z}^4\mu_{z}(\d z)},
\ee
where $\mu_z$ denotes the marginal of $\mu$ in the direction of the fast variable $z$. The quantity $Q$ in \eqref{Eq_param_def_cont} can be seen as the asymptotic counterpart to the parameterization defect introduced in \cite{CLW15_vol2}. By numerically approximating $Q$ in \eqref{Eq_param_def_cont} we have in Case I, $Q=1.789 \times 10^{-1}$, while in Case II, $Q=5.834 \times 10^{-1}$. Small values of the parameterization defect that seem auspicious for  a good approximation of the long-term statistics by the reduced system Eq.~\eqref{Eq_2DHopf_polar_base} according to our error estimates.

A more careful inspection of the parameter regimes indicates though that Theorem \ref{Thm_hopf_1} predicts indeed a good approximation in Case I but is unable to conclude for Case II. 

Nevertheless, numerical results show that in both cases the radial probability density function (PDF) of the original system Eq.~\eqref{Eq_3DHopf_polar_base} is well approximated by that of the reduced system  Eq.~\eqref{Eq_2DHopf_polar_base}; see Figures \ref{PDF1} and \ref{PDF2}.  
The numerical results for Case II combined with the smallness of the parameterization defect suggest thus that an inequality analogous to \eqref{Eq_main_result_1} is valid beyond what Theorem \ref{Thm_hopf_1} predicts.

In that respect, we conjecture that holds the following error estimate
\be \label{Eq_main_result_1}
	D_{\mathcal{F}}(\mu_r, \nu_r) \leq C \left(\int_ {\mathbb{R}_+ \times \mathbb{R}}|z - r^2|^4\,\mu_{r, z}(\d r, \d z)\right)^{\frac{1}{4}},
	\ee
for some constant $C>0$.

\begin{figure}[htbp]
	\centering
	\includegraphics[width=0.9\textwidth,height=.45\textwidth]{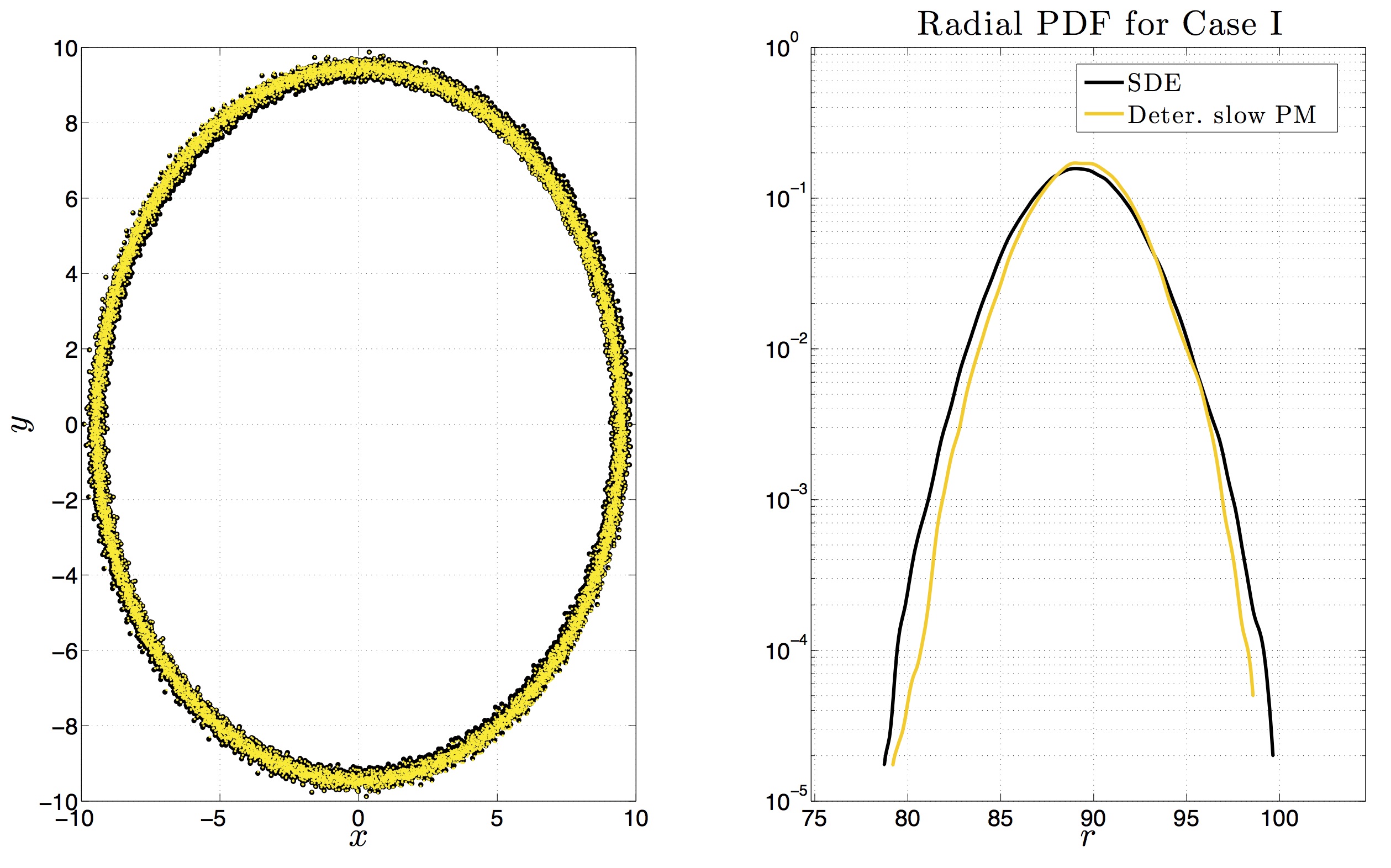}
	\caption{\footnotesize{\bf Modeling skills using the slow manifold $h(x,y)=x^2+y^2$.} Here for Case I, see Table \ref{Table_CaseI}.}
	\label{PDF1}
\end{figure}

\begin{figure}[htbp]
	\centering
	\includegraphics[width=1\textwidth,height=.45\textwidth]{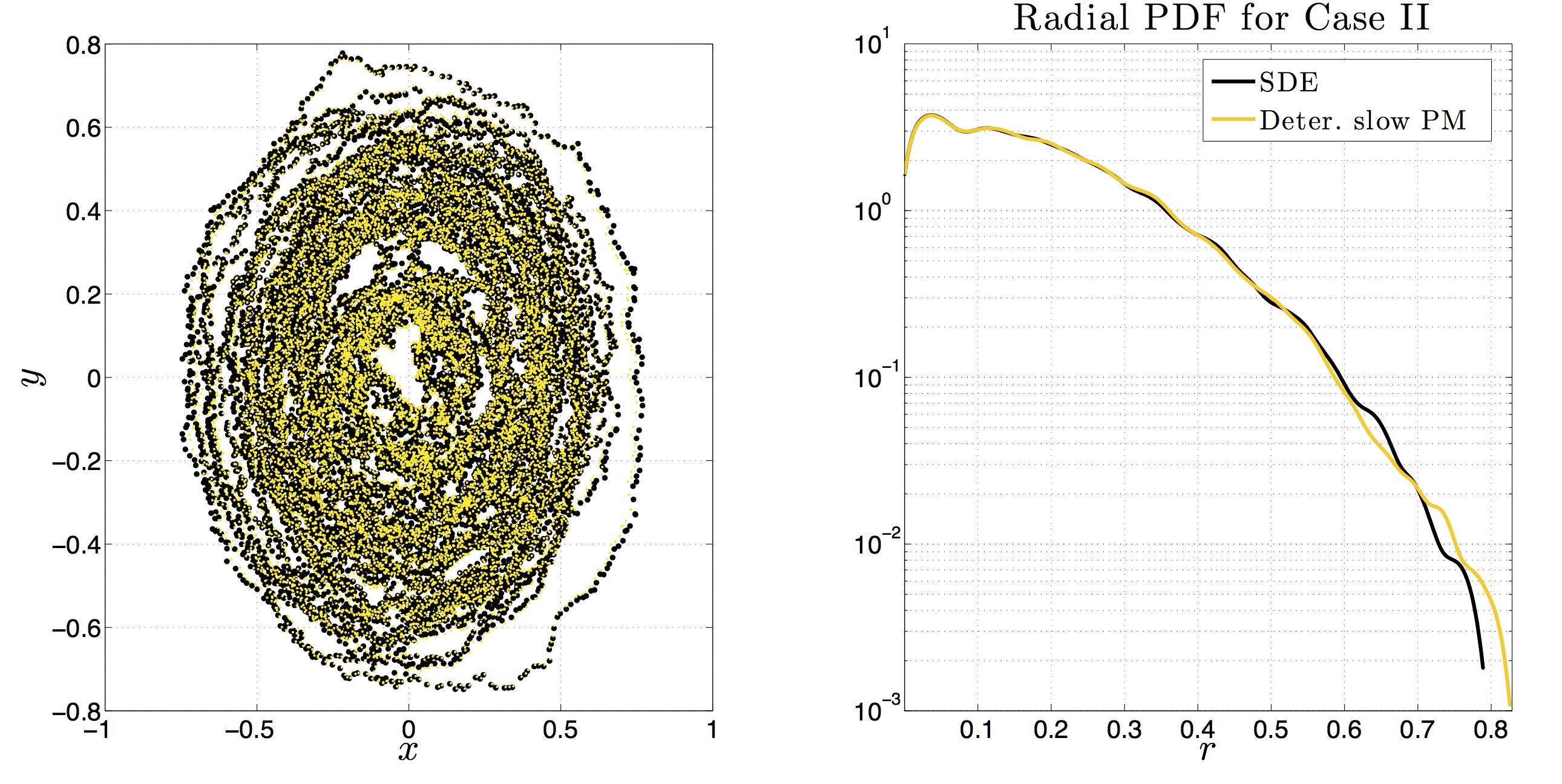}
	\caption{\footnotesize{\bf Modeling skills using the slow manifold $h(x,y)=x^2+y^2$.} Here for Case II, see Table \ref{Table_CaseI}.}
	\label{PDF2}
\end{figure}

\section{Stochastic parameterizing manifolds: Close to criticality and $\epsilon >1$}\label{Sec_stoch_pm}
Before introducing the concept of stochastic parameterizing manifolds, we set some notations. First of all, let $(W^r_t, W^{\theta}_t, W^3_t, W^4_t)_{t \geq 0}$ be a standard $\mathbb{R}^4$-valued Brownian motion. The real-valued Brownian motion $W^4_t$ is independent of $W^r_t$, $W^{\theta}_t$ and $W^3_t$, which are as in Section \ref{Sec_main_idea}.
Furthermore, the filtration on $(\Omega, \mathcal{B}, \mathbb{P})$ defined in \eqref{Eq_filtration} is replaced by the filtration $\{\mathcal{B}_t\}_{t\geq0}$ given by
\be \label{Eq_filtration_2}
\mathcal{B}_t = \sigma\left(W^r_s, W^{\theta}_s, W^3_s, W^4_s: 0\leq s\leq t \right), \,\,\, t \geq 0,
\ee
and we consider the filtered probability space $(\Omega, \mathcal{B}, \{\mathcal{B}_t\}_{t\geq 0} , \mathbb{P})$.

\subsection{"Slow" stochastic parameterizing manifolds}\label{Sec_intro_stoch_pm}

In this section, we introduce generalizations of the deterministic slow manifold \eqref{Eq_slow_mnf_polar} associated with the original system \eqref{Eq_3DHopf_polar_base}. Specifically, we consider the following family of parameterizations:
\bea \label{Eq_stoch_mnf}
h_{\tau}: \, & \mathbb{R} \times \mathbb{R}_+ \longrightarrow \mathbb{R},\\
&(m, r) \mapsto m +  c_{\tau} r^2,\,\,\, c_{\tau} = \left(1 - e^{-\frac{\tau}{\epsilon}}\right),
\eea 
where $\tau$ is a free parameter in $(0, \infty)$.

The primary objective of $h_{\tau}$ is to parameterize  the variable $z$ of the original system \eqref{Eq_3DHopf_polar_base} as a function of the radial part $r$. Since system \eqref{Eq_3DHopf_polar_base} is inherently stochastic, the parameter $m$ in  \eqref{Eq_stoch_mnf} is chosen to be a stochastic process designed to capture the stochastic fluctuations of $z(t)$ away from the deterministic slow manifold. 
This way, when $m$ in \eqref{Eq_stoch_mnf} is chosen to be stochastic, one introduces a family of stochastic parameterizations indexed by 
$\tau$ in $(0,\infty)$.  The parameter $\tau$ effectively controls the curvature of the deterministic component of the manifold $h_{\tau}$, providing flexibility in adapting to different dynamical regimes.

In what follows, we take this stochastic dependence to be driven by the (stationary) stochastic process $(M_t)_{t\geq 0}$ solving the following stochastic differential equation:
\begin{equation} \label{Eq_def_M}
\d M_t = -\frac{1}{\epsilon} M_t^3\d t + \frac{\sigma}{\sqrt{\epsilon}}\d W^4_t, \;\; t \geq0.
\end{equation}
The resulting family of mappings $(h_{\tau}(M_t, \cdot))_{t\geq0}$, forms what we call a Stochastic Parameterization Manifold (SPM) with parameter $\tau$.

Before we state the main result of this section, we formulate a few useful systems. 
Recall the definitions of the state spaces, $\mathcal{M}$ and $\hat{\mathcal{M}}$, defined in \eqref{Def_M_hat_M}.
As in Section \ref{Sec_slow_pm}, we consider the original system on $\mathcal{M}$:
\begin{subnumcases}{\label{Eq_3DHopf_polar_base_stoch}}
\d r = \left[\lambda r -\gamma r z + \frac{\sigma^2}{2 r}\right]\d t + \sigma \d W^r_t \label{Eq_r_stoch}\\
\d \theta = f \d t +  \frac{\sigma}{r} \d W^{\theta}_t \mod (2\pi) \label{Eq_theta_stoch} \\
\d z = -\frac{1}{\epsilon}\left[z -r^2\right]\d t + \frac{\sigma}{\sqrt{\epsilon}} \d W^3_t. \label{Eq_z_stoch} 
\end{subnumcases}
On the extended phase space $\mathcal{M}\times \mathbb{R}$, we introduce the augmented original system
\begin{subnumcases}{\label{Eq_3DHopf_polar_base_aug}}
\d r = \left[\lambda r -\gamma r z + \frac{\sigma^2}{2 r}\right]\d t + \sigma \d W^r_t \label{Eq_r_aug}\\
\d \theta = f \d t +  \frac{\sigma}{r} \d W^{\theta}_t \mod (2\pi) \label{Eq_theta_aug} \\
\d z = -\frac{1}{\epsilon}\left[z -r^2\right]\d t + \frac{\sigma}{\sqrt{\epsilon}} \d W^3_t \label{Eq_z_aug} \\
\d M = -\frac{1}{\epsilon} M^3 \d t + \frac{\sigma}{\sqrt{\epsilon}}\d W^4_t.\label{Eq_M_aug}
\end{subnumcases}
Building upon the parameterization $h_{\tau}$, we consider the following reduced system
\begin{subnumcases}{\label{Eq_2DHopf_polar_base_aug}}
\d \hat{r} = \left[\lambda \hat{r} -\gamma c_{\tau} \hat{r}^3 + \frac{\sigma^2}{2\hat{r}} -\gamma \hat{r} M_t \right]\d t  + \sigma \d W^r_t \label{Eq_r_hat_aug}\\
\d \hat{\theta} = f \d t +  \frac{\sigma}{\hat{r}} \d W^{\theta}_t \mod(2\pi) \label{Eq_theta_hat_aug}\\
\d M = -\frac{1}{\epsilon} M^3 \d t + \frac{\sigma}{\sqrt{\epsilon}}\d W^4_t \label{Eq_M_red_aug},
\end{subnumcases}
in which the term $z$ in Eq.~\eqref{Eq_r_aug} has been replaced by $ M_t + c_\tau \hat{r}^2$. 

For $\hat{R}_0 = (\hat{r}_0, \hat{\theta}_0, M_0)$ in $\hat{\mathcal{M}}\times \mathbb{R}$ the unique strong solution to system \eqref{Eq_2DHopf_polar_base_aug} shall be denoted by $(\hat{R}_t^{\hat{R}_0})_{t\geq 0}$ (see Remark \ref{Remark_ex_str_sol_2} below), with
\be \label{Eq_R_hat_aug}
\hat{R}_t^{\hat{R}_0} = (\hat{r}_t, \hat{\theta}_t, M_t), \,\,\, t \geq 0.
\ee
Analogously to system \eqref{Eq_3DHopf_polar_base_transf} in Section \ref{Sec_slow_pm} we seek an appropriate transformation of the original system \eqref{Eq_3DHopf_polar_base_stoch}. To this end, the quantity $q$ defined in \eqref{Eq_def_q_girs} is replaced by the following constant, which is also denoted by $q$
\be \label{Eq_def_q_2}
q = \frac{1}{\epsilon \gamma}\left(1+\frac{5}{2\sigma^2}\right).
\ee
 We introduce the function
\bea \label{Eq_g_stoch}
g: \, & \mathbb{R} \times \mathbb{R}_+ \times \mathbb{R}_+ \longrightarrow \mathbb{R},\\
&(m, r_1, r_2) \mapsto \frac{q+\lambda }{\sigma} (r_2 - r_1) - \frac{\gamma}{\sigma}m r_2.
\eea
Given the radial component $\hat{r}$ and the process $M_t$ in Eqns.~\eqref{Eq_r_hat_aug} and \eqref{Eq_M_red_aug}, we define the transformed system to the original system \eqref{Eq_3DHopf_polar_base_stoch} as
\begin{subnumcases}{\label{Eq_3DHopf_polar_base_transf_aug}}
\d \tilde{r} = \left[\lambda \tilde{r} -\gamma \tilde{r}\tilde{ z} + \frac{\sigma^2}{2 \tilde{r}}\right]\d t + \sigma g(M,\tilde{r}, \hat{r}) \d t + \sigma \d W^r_t \label{Eq_r_transf_aug}\\
\d \tilde{\theta} = f \d t +  \frac{\sigma}{\tilde{r}} \d W^{\theta}_t \mod(2\pi) \label{Eq_theta_transf_aug} \\
\d \tilde{z} = -\frac{1}{\epsilon}\left[\tilde{z} -\tilde{r}^2\right]\d t + \frac{\sigma}{\sqrt{\epsilon}} \d W^3_t.
\end{subnumcases}
For an initial value $R_0$ in $\mathcal{M}$ let $(\tilde{R}_t^{R_0})_{t\geq0}$ be the unique strong solution to system \eqref{Eq_3DHopf_polar_base_transf_aug} (see Remark \ref{Remark_ex_str_sol_2} below), with
\be \label{Eq_R_tilde_aug}
\tilde{R}_t^{R_0} = (\tilde{r}_t, \tilde{\theta}_t, \tilde{z}_t), \,\,\, t \geq 0.
\ee
\br \label{Remark_ex_str_sol_2}
Existence and uniqueness of global solutions to the systems in this section are due to an analogous reasoning to that of Remark \ref{Remark_ex_str_sol}. In particular \cite[Theorem 2.3, p.~173]{watanabe_ikeda} and \cite[Theorem 3.1, p.~178-179]{watanabe_ikeda} ensure the existence of unique local strong solutions to the augmented original, reduced and transformed  system in \eqref{Eq_3DHopf_polar_base_aug}, \eqref{Eq_2DHopf_polar_base_aug} and \eqref{Eq_3DHopf_polar_base_transf_aug}, respectively.
A Lyapunov function approach then implies global existence and uniqueness of strong solutions; see also definition \eqref{Eq_V_lyap_stoch} in Lemma \ref{Lemma_invm_ex_org_aug} below.
\er
The following lemma is the counterpart of Lemma \ref{Lemma_girsanov} in Section \ref{Sec_slow_pm} for stochastic parameterizing manifolds discussed here.
\bl \label{Lemma_girsanov_stoch_mnf}
	Let $T>0$. For $R_0$ in $\mathcal{M}$, let $(R_t^{R_0})_{t\in[0,T]}$ be the unique strong solution to the original system \eqref{Eq_3DHopf_polar_base_stoch}.
	Furthermore, the augmented reduced system on $\hat{\mathcal{M}}\times \mathbb{R}$ in \eqref{Eq_2DHopf_polar_base_aug} has the unique strong solution $(\hat{R}_t^{\hat{R}_0})_{t\in [0,T]}$ for an initial value $\hat{R}_0$ in $\hat{\mathcal{M}} \times \mathbb{R}$.
	In addition, consider $(\tilde{R}_t^{R_0})_{t\in [0,T]}$, the unique strong solution of the transformed  system \eqref{Eq_3DHopf_polar_base_transf_aug}.
	Recall the definition of $g$ in \eqref{Eq_g_stoch}.

Then, the random variable $\mathfrak{D}_T$ on $(\Omega, \mathcal{B}, \mathbb{P})$ of the form
	\be \label{Eq_D_girs_2}
	\mathfrak{D}_T = \exp\left(-\int_0^T g(M_s, \tilde{r}_s, \hat{r}_s)\d W^r_s - \frac{1}{2}\int_0^T \betrag{g(M_s,\tilde{r}_s, \hat{r}_s)}^2\d s \right),
	\ee
	is integrable with respect to $\mathbb{P}$ and satisfies 
	\be \label{Eq_D_1_2}
	\mathbb{E}_{\mathbb{P}} [\mathfrak{D}_T] = 1.
	\ee
	Hence, $\mathfrak{D}_T$ induces a probability measure $\tilde{\mathbb{P}} = \tilde{\mathbb{P}}_T$ on $(\Omega, \mathcal{B})$ given by
	\be \label{Eq_P_tilde_girs_2}
	\tilde{\mathbb{P}}(\d \omega) = \mathfrak{D}_T(\omega)\mathbb{P}(\d \omega),\,\,\,\text{for}\,\,\,\,\omega \in \Omega.
	\ee
	Both, $\mathbb{P}$ and $\tilde{\mathbb{P}}$ are equivalent on the $\sigma$-algebra $\mathcal{B}_{T}$. Furthermore,
	the stochastic process $(\tilde{W}^r_t, \tilde{W}^{\theta}_t, \tilde{W}^3_t)_{t\in [0,T]}$ defined as
		 \bea \label{Eq_girsanov_tilde_W^2}
	\begin{pmatrix} 
		\tilde{W}^r_t \\
		\tilde{W}^{\theta}_t \\
		\tilde{W}^3_t 
	\end{pmatrix} = 
	\int_0^t G(M_s,\tilde{r}_s, \hat{r}_s)\d s +
	\begin{pmatrix} 
		W^r_t\\
		W^{\theta}_t \\
		W^3_t 
	\end{pmatrix}, \,\,\, t \in [0,T],
	\eea
	where 
	\bes
	G(M_t, \tilde{r}_t, \hat{r}_t)  = \begin{pmatrix} 
		g(M_t,\tilde{r}_t, \hat{r}_t) \\
		0 \\
		0
	\end{pmatrix}, \,\,\, t \in [0,T], 
	\ees
	is an $\mathbb{R}^3$-valued Brownian motion. Furthermore, it holds that
	\be \label{Eq_trans_prob_girs_2}
	\tilde{\mathbb{P}}\Big(\tilde{R}_t^{R_0} \in \Gamma\Big) = \mathbb{P}\Big(R_t^{R_0} \in \Gamma\Big),\,\,\, \Gamma \in \mathcal{B}( \mathcal{M}), \; R_0 \in \mathcal{M}, \,\,\,t\in [0,T].
	\ee
\el
\begin{proof}
	The proof is analogous to the proof of Lemma \ref{Lemma_girsanov}, where we showed the applicability of Girsanov's theorem \ref{Lemma_girsanov_general} to the transformed  system \eqref{Eq_3DHopf_polar_base_transf}. In particular, here it is sufficient to prove that there exists a $\delta>0$ such that 
	\be \label{Eq_novikov_rep_22}
	\sup_{t\in[0,T]} \mathbb{E}_{\mathbb{P}}\left[\exp\left(\delta |g(M_t, \tilde{r}_t, \hat{r}_t)|^2\right)\right] < \infty.
	\ee
	The validity of inequality \eqref{Eq_novikov_rep_22} can be proven as in Step 1 of the proof of Lemma \ref{Lemma_girsanov}, where the exponential moment bound of Lemma \ref{Lemma_moment_bounds} is replaced by the corresponding bound shown in Lemma \ref{Lemma_moment_bounds2} below.	
\end{proof}

\needspace{1\baselineskip}
\br \label{Remark_ex_inv_m_2}
\hspace*{2em}  \vspace*{-0.4em}
\bi
\item[(i)]
Due to Lemma \ref{Lemma_invm_ex_org_aug} in Appendix \ref{Appendix_E} below, there exists a unique ergodic invariant measure $\bar{\mu}$ in $Pr(\mathcal{M}\times \mathbb{R})$ with smooth Lebesgue-density for the augmented original system \eqref{Eq_3DHopf_polar_base_aug}. The fact that Eq.~\eqref{Eq_M_aug} is not coupled to equations \eqref{Eq_r_aug}-\eqref{Eq_z_aug} allows us to conclude that $\bar{\mu}$ is of the form
\be \label{Eq_prod_meas_st_pm}
\bar{\mu} = \mu \otimes \rho\,\,\, \text{on}\,\,\, \mathcal{M}\times \mathbb{R},
\ee
where $\mu$ is the unique ergodic invariant measure for the original system \eqref{Eq_3DHopf_polar_base_stoch} as discussed in Section \ref{Sec_slow_pm} (see also Remark \ref{Remark_ex_inv_m_1}). Moreover, $\rho$ denotes the unique ergodic invariant measure for the process $M_t$  solving Eq.~\eqref{Eq_M_aug}, where
\be \label{Eq_def_rho}
\rho(\d m) = \frac{1}{Z} e^{-\frac{m^4}{2\sigma^2}} \d m\,\,\,\text{with}\,\,\, Z = \int_{\mathbb{R}} e^{-\frac{m^4}{2\sigma^2}}\d m,
\ee
see e.g. \cite[Proposition 4.2, p. 110]{Pavliotis:2014}.
\item[(ii)]
For the reduced system \eqref{Eq_2DHopf_polar_base_aug} there exists a unique ergodic invariant measure $\bar{\nu}$ in $Pr(\mathcal{\hat{M}} \times \mathbb{R})$ with smooth Lebesgue-density (see Lemma \ref{Lemma_invm_exuniq_red_stoch} below).
\ei  
\er

\br \label{Remark_OU}
In Section \ref{Sec_numerics_2} below we present numerical results concerning the approximation properties of the marginal measure $\bar{\nu}_r$ relative to $\mu_r$. There, the stochastic parameterizing manifold \eqref{Eq_stoch_mnf} used to obtain the reduced system \eqref{Eq_2DHopf_polar_base_aug}, is driven by the Ornstein-Uhlenbeck (OU)-process solving 
\be \label{Eq_def_OU}
\d M_t = -\frac{1}{\epsilon} M_t\d t + \frac{\sigma}{\sqrt{\epsilon}}\d W^4_t, \,\,\, t \geq0.
\ee
 instead of the process $M_t$ defined in Eq.~\eqref{Eq_def_M}. 
Compared to an OU-process,  the stochastic process $M_t$ solving Eq.~\eqref{Eq_def_M} enjoys better dissipation properties. This allows us to prove inequality \eqref{Eq_novikov_rep_22} in Lemma \ref{Lemma_girsanov_stoch_mnf}, which ensures the applicability of Girsanov's theorem to the transformed system \eqref{Eq_3DHopf_polar_base_transf_aug}.
\er

\subsection{Statement and proof of Theorem \ref{Thm_hopf_stoch}} 
\label{Sec_main_thm_stoch_pm}
\bt \label{Thm_hopf_stoch}
Let $\mu(\d r,\d \theta, \d z)$ be the unique ergodic invariant measure for the original system \eqref{Eq_3DHopf_polar_base_stoch}.
The marginals of $\mu$ in the $(r,z)$-plane and along the radial component $r$ in the phase space $\mathcal{M}$ are denoted by $\mu_{r,z}$ and $\mu_r$, respectively.

Let $\bar{\nu}(\d r,\d \theta,\d m)$ be the unique ergodic invariant measure for the augmented reduced system \eqref{Eq_2DHopf_polar_base_aug} and $\bar{\nu}_r$ its marginal along $r$.
Furthermore, let
\be \label{Eq_r_*_stoch}
r_* = \left(\frac{\lambda}{2\gamma} + \frac{1}{2}\sqrt{\left(\frac{\lambda}{\gamma}\right)^2 + \frac{2\sigma^2}{\gamma}}\right)^\frac{1}{2},
\ee
which is the unique positive root of 
\bes
\lambda r - \gamma r^3 + \frac{\sigma^2}{2r} =0.
\ees
We consider the stochastic parameterizing manifold $h_{\tau}$ defined in \eqref{Eq_stoch_mnf} for $\tau$ in $(0,\infty)$. Then there exist constants $c>0$  and $C(r_*, \mu_r, \bar{\nu}_r)>0$ such that
\be \label{Eq_main_result_1_stoch}
D_{\mathcal{F}}(\mu_r, \bar{\nu}_r) \leq C(r_*, \mu_r, \bar{\nu}_r) + c \left(\int_ {\mathbb{R}_+ \times \mathbb{R}^2}|z - h_{\tau}(m,r)|^4\,\mu_{r, z}(\d r, \d z) \otimes \rho(\d m)\right)^{\frac{1}{4}},
\ee
where $\mathcal{F}$ is as in \eqref{Eq_test_fcts_1} and the constants $c$ and $C(r_*, \mu_r, \bar{\nu}_r)$ are explicitly given by
\bes
c = \left(\frac{\gamma}{q}\right)\left(\int_0^{\infty}r^4\mu_r(\d r)\right)^{\frac{1}{4}},
\ees
with $q$ as in \eqref{Eq_def_q_2} and
\be \label{Eq_def_C_*_2}
C(r_*, \mu_r, \bar{\nu}_r) = 2r_* + \int_{r_*}^{\infty} r \bar{\nu}_r(\d r) +\int_{r_*}^{\infty} r \mu_r(\d r).
\ee
\et

\begin{proof}
	The proof of Theorem \ref{Thm_hopf_stoch} is analogous to that of Theorem \ref{Thm_hopf_1}.
	We start by introducing Markov semigroups for the original system \eqref{Eq_3DHopf_polar_base_stoch}, the corresponding reduced system \eqref{Eq_2DHopf_polar_base_aug} and the transformed  system \eqref{Eq_3DHopf_polar_base_transf_aug}.
	
	For $R_0 $ in $\mathcal{M}$ let $(R_t^{R_0})_{t\geq0}$ be the unique strong solution to the original system \eqref{Eq_3DHopf_polar_base_stoch}. The associated Markov semigroup $P_t$ is defined by
	\be \label{Eq_P_aug_proof}
	[P_t\psi](R_0) = \mathbb{E}_{\mathbb{P}}\big[\psi\big(R_t^{R_0}\big)\big], \,\,\, \psi \in C_b(\mathcal{M}), \,\,\, R_0 \in \mathcal{M},\,\,\, t \geq 0.
	\ee 
	Due to Remark \ref{Remark_ex_inv_m_1}, there exists a unique ergodic invariant measure $\mu$ for $P_t$ on $\mathcal{M}$. Therefore, by virtue of \cite[Theorem 5.8, p. 73]{daprato_inf_dim} it follows that $P_t$ can be uniquely extended to $L^2_{\mu}(\mathcal{M})$ for all $t\geq0$.
	
	Similarly, for the reduced system \eqref{Eq_2DHopf_polar_base_aug}, the Markov semigroup $Q_t$ is given by
	\be \label{Eq_Q_aug_proof}
	[Q_t\psi](\hat{R}_0) = \mathbb{E}_{\mathbb{P}}\big[\psi\big(\hat{R}_t^{\hat{R}_0}\big)\big], \,\,\, \psi \in C_b(\hat{\mathcal{M}}\times \mathbb{R}), \,\,\, \hat{R}_0 \in \hat{\mathcal{M}}\times\mathbb{R},\,\,\, t \geq 0,
	\ee	
where  $ \hat{R}_t^{\hat{R}_0}= (\hat{r}_t, \hat{\theta}_t, M_t)$ denotes the unique strong solution to Eq.~\eqref{Eq_2DHopf_polar_base_aug}.
This Markov semigroup also possesses  a unique ergodic invariant measure $\bar{\nu}$ in $Pr(\hat{\mathcal{M}}\times \mathbb{R})$; see Remark \ref{Remark_ex_inv_m_2}. Hence, \cite[Theorem 5.8, p. 73]{daprato_inf_dim} ensures that $Q_t$ is uniquely extendible to $L^2_{\bar{\nu}}(\hat{\mathcal{M}}\times \mathbb{R})$.
	
	Next we turn to the transformed  system \eqref{Eq_3DHopf_polar_base_transf_aug} with a unique strong solution $(\tilde{R}_t^{R_0})_{t\geq 0}$.
	For every $T_*>0$, Lemma \ref{Lemma_girsanov_stoch_mnf} ensures the existence of a probability measure $
	\tilde{\mathbb{P}} = \tilde{\mathbb{P}}_{T_*}$ on $(\Omega, \mathcal{B})$, which is equivalent to $\mathbb{P}$ on $\mathcal{B}_{T_*}$ such that 
	\be \label{Eq_trans_aug_prob_proof}
	\tilde{\mathbb{P}}\Big(\tilde{R}_t^{R_0} \in \Gamma\Big) = \mathbb{P}\Big(R_t^{R_0} \in \Gamma\Big),\,\,\, \Gamma \in \mathcal{B}( \mathcal{M}), \; R_0 \in \mathcal{M}, \,\,\,t\in [0,T_*].
	\ee
	As for \eqref{Eq_trans_prob_proof}, we conclude from identity \eqref{Eq_trans_aug_prob_proof} that transition probabilities of system \eqref{Eq_3DHopf_polar_base_stoch} are preserved in its transformed counterpart \eqref{Eq_3DHopf_polar_base_transf_aug} up to a final time $T_*$. However, this property only holds, if the transformed system \eqref{Eq_3DHopf_polar_base_transf_aug} is viewed under the new probability measure $\tilde{\mathbb{P}}$.

We define
	\be \label{Eq_P_T_*_aug}
	[\tilde{P}^{T_*}_t\psi](R_0) = \mathbb{E}_{\tilde{\mathbb{P}}}\big[\psi\big(\tilde{R}_t^{R_0}\big)\big], \,\,\, \psi \in C_b(\mathcal{M}), \,\,\, R_0 \in \mathcal{M},\,\,\, t \in [0,T_*].
	\ee
By virtue of identity \eqref{Eq_trans_aug_prob_proof} and the fact that $P_t$ is defined for all $\psi$ in $L^2_{\mu}(\mathcal{M})$, we can set
	\be \label{Eq_proof_1_aug}
	[\tilde{P}^{T_*}_t\psi](R_0) = [P_t\psi](R_0),\,\,\, \psi \in L^2_{\mu}(\mathcal{M}), \,\,\, R_0 \in \mathcal{M},\,\,\, t \in [0,T_*]. 
	\ee 
	Let $e_1, e_2, e_3$ be the canonical basis vectors of $\mathbb{R}^3$.
	We consider subspaces of $\mathbb{R}^3$ of the form
	\bes
	V_{r,z} = \span\{e_1, e_3\} \,\,\,\text{and}\,\,\, V_r = \span\{e_1\}.
	\ees
	Associated with these subspaces we introduce projections $\Pi_V$, for $V=V_{r,z}$ and $V=V_r$, which are defined as in \eqref{Def_proj}.
	Furthermore, let 
	\bes
	\mu_{r,z} = \Pi_{V_{r,z}}^* \mu, \,\, \mu_r = \Pi_{V_r}^* \mu, \,\, \bar{\nu}_r = \Pi_{V_r}^* \bar{\nu}.
	\ees

We replace $f$ in \eqref{Eq_f_proof} by the following function defined as
	\bea \label{Eq_f_proof_aug}
	f: \, & \mathbb{R}_+\times \mathbb{R} \times \mathbb{R} \longrightarrow \mathbb{R},\\
	&(r,z,m	) \mapsto  \left(\frac{\gamma}{q} \right)^2 r^2 (z - h_{\tau}(m,r))^2,
	\eea
	where $h_{\tau}$ is given in \eqref{Eq_stoch_mnf} for $\tau$ in $(0,\infty)$. 
The proof then proceeds analogously to that of Theorem \ref{Thm_hopf_1}, commencing from \eqref{Eq_erg_1} and utilizing the definitions established in \eqref{Eq_P_aug_proof}-\eqref{Eq_f_proof_aug}. The usage of Lemma \ref{Lemma_aux_ineq_2} in Eq.~\eqref{Eq_proof_thm_1_aux_1} is replaced by Lemma \ref{Lemma_aux_ineq_4}. Additionally,
	Lemma \ref{Lemma_aux_ineq_3} replaces Lemma \ref{Lemma_aux_ineq_1} in \eqref{Eq_proof_thm_1_dep_phi}.
	
Hence, for an arbitrary $\varphi$ in the function space $\mathcal{F}$ and $\delta>0$, we arrive at an inequality analogous to \eqref{Eq_proof_final_1} of the form
	\bea \label{Eq_proof_final_2}
	\left|\int_{\R_+} \varphi(r) \, \mu_r(\d r) - \int_{\R_+} \varphi(r) \, \bar{\nu}_r(\d r)\right| 
	&\leq 4\delta +2r_* +\int_{r_*}^{\infty} r \bar{\nu}_r(\d r) +\int_{r_*}^{\infty} r \mu_r(\d r)\\
	&+\Bigg(\delta + \left(\frac{\gamma}{q}\right)^2  \int_{\mathbb{R}_+ \times \mathbb{R}} 
	r^2\left(z-h_{\tau}(m,r)\right)^2 \mu_{r,z}(\d r, \d z)  \\
	& +m \frac{2\gamma}{q} \frac{1}{T_*} \int_0^{T_*} \mathbb{E}_{\tilde{\mathbb{P}}}[\tilde{r}_t(\hat{r}_t - \tilde{r}_t)]\d t \Bigg)^{\frac{1}{2}}.
	\eea
	Note that inequality \eqref{Eq_proof_final_2} holds for all $m$ in $\mathbb{R}$. By integrating both sides of \eqref{Eq_proof_final_2} against $\rho(\d m)$ defined in \eqref{Eq_def_rho} over the interval $[-a,a]$ ($a>0$) yields
	\bea \label{Eq_proof_final_3}
	\Bigg|\int_{\R_+} \varphi(r) \, \mu_r(\d r) &- \int_{\R_+} \varphi(r) \, \bar{\nu}_r(\d r)\Bigg|  \rho([-a,a])\\
	&\leq \left(4\delta + C(r_*,\mu_r,\bar{\nu}_r)\right)\rho([-a,a])\\
	&+\Bigg(\delta \rho([-a,a]) + \left(\frac{\gamma}{q}\right)^2 \int_{-a}^{a} \left( \int_{\mathbb{R}_+ \times \mathbb{R}} 
	r^2\left(z-h_{\tau}(m,r)\right)^2 \mu_{r,z}(\d r, \d z) \right) \rho(\d m) 
	 \\
	& +\int_{-a}^a m \rho(\d m) \frac{2\gamma}{q} \frac{1}{T_*} \int_0^{T_*} \mathbb{E}_{\tilde{\mathbb{P}}}[\tilde{r}_t(\hat{r}_t - \tilde{r}_t)]\d t \Bigg)^{\frac{1}{2}},
	\eea
where the latter inequality is due to the Cauchy-Schwarz inequality and the constant $C(r_*, \mu_r, \bar{\nu}_r)$ is as in \eqref{Eq_def_C_*_2}.
	The symmetry of $\rho$ implies that the last term in inequality \eqref{Eq_proof_final_3} vanishes. 
	Letting $a$ tend to infinity on both sides of \eqref{Eq_proof_final_3} and using the monotone convergence theorem for the last term in \eqref{Eq_proof_final_3} results in 
	\bea \label{Eq_proof_final_4}
						\Bigg|\int_{\R_+} \varphi(r) \, \mu_r(\d r) &- \int_{\R_+} \varphi(r) \, \bar{\nu}_r(\d r)\Bigg|\\
	&\leq 4\delta + C(r_*,\mu_r,\bar{\nu}_r)\\
	&+\left(\delta  + \left(\frac{\gamma}{q}\right)^2 \int_{\mathbb{R}} \left( \int_{\mathbb{R}_+ \times \mathbb{R}} 
	r^2\left(z-h_{\tau}(m,r)\right)^2 \mu_{r,z}(\d r, \d z)\right) \rho(\d m) \right)^{\frac{1}{2}}.
	\eea
	Taking the supremum over all $\varphi$ in $\mathcal{F}$ on both sides of \eqref{Eq_proof_final_4} and letting $\delta$ tend to zero yields
	\be \label{Eq_final_thm2}
	D_{\mathcal{F}}(\mu_r,\bar{\nu}_r) \leq C(r_*, \mu_r,\bar{\nu}_r) + c\left(\int_{\mathbb{R}_+\times \mathbb{R}^2} \left|z-h_{\tau}(m,r)\right|^4 \mu_{r, z}(\d r, \d z) \otimes \rho(\d m)\right)^{\frac{1}{4}},
	\ee
	after application of the Cauchy-Schwarz inequality. The constant $c$ in \eqref{Eq_final_thm2} is given as
	\bes
	c = \left(\frac{\gamma}{q}\right) \left(\int_0^{\infty} r^4\mu_r(\d r)\right)^{\frac{1}{4}}.
	\ees
	This completes the proof of Theorem \ref{Thm_hopf_stoch}.
\end{proof}

\subsection{Statement and proof of Theorem \ref{Thm_hopf_22}} \label{Sec_main_thm_stoch_pm2}

The following Theorem \ref{Thm_hopf_22} is the counterpart of Theorem \ref{Thm_hopf_2}  stated in the case of the slow parameterizing manifold of Section \ref{Sec_slow_pm}.
\bt \label{Thm_hopf_22}
	Let $\mu$ be the ergodic invariant measure for the original system \eqref{Eq_3DHopf_polar_base_stoch} and denote by $\mu_{r,z}$ and $\mu_r$ its marginals in the $(r,z)$-plane and along the radial part $r$, respectively. Let $\bar{\nu}$ be the ergodic invariant measure for the reduced system \eqref{Eq_2DHopf_polar_base_aug} and $\bar{\nu}_r$ its marginal along $r$.
	Furthermore, let
	\be \label{Eq_r_det2}
	r_{\mathrm{det}} = \sqrt{\frac{\lambda}{\gamma}},
	\ee
	which is the unique positive root of 
	\bes
	\lambda r - \gamma r^3 =0.
	\ees
	Recall the definition of function space $\mathcal{F}$ in \eqref{Eq_test_fcts_1}. 
	Then there exist constants $c>0$ and $C(\lambda, \gamma, \sigma, \tau,\bar{\nu}_r)>0$ such that
	\be \label{Eq_main_result_22}
	D_{\mathcal{F}}(\mu_r, \bar{\nu}_r) \leq C(\lambda, \gamma, \sigma,\tau, \bar{\nu}_r) + c \left(\int_ {\mathbb{R}_+ \times \mathbb{R}^2}|z - h_{\tau}(m,r)|^4\,\mu_{r, z}(\d r, \d z)\otimes \rho(\d m)\right)^{\frac{1}{4}},
	\ee
	where 
	\be \label{Eq_closeness_to_crit_22}
	C(\lambda, \gamma, \sigma,\tau, \bar{\nu}_r) = \sqrt{2\left(\frac{r_{\mathrm{det}}^2}{c_{\tau}} + \frac{\sigma^2}{\lambda}\right)}+\int_0^{\infty} r \bar{\nu}_r(\d r),
	\ee
	and
	\bes
	c =\frac{q+2\lambda}{qr_{\mathrm{det}}^2}\left(\int_0^{\infty}r^4\mu_r(\d r)\right)^{\frac{1}{4}},
	\ees
	with $q$ as in \eqref{Eq_def_q_2}.
\et
\begin{proof}
The proof  proceeds analogously to that  of Theorem \ref{Thm_hopf_2} in Section  \ref{Sec_proof_thm_3}. We only replace the application of Lemma \ref{lemma_aux_ineq_new_lemma} in \eqref{Eq_proof_thm_1_aux_2}  therein, by Lemma \ref{lemma_aux_ineq_new_lemma2} here.
\end{proof}

\subsection{Stochastic parameterization for "inverted" timescale separation}
\label{Sec_numerics_2}
We extend our analysis to regimes exhibiting an inverted timescale separation, where the variable to be parameterized evolves on a slower timescale than the resolved variables. To investigate this scenario, we examine the parameterization performance of our Stochastic Parameterizing Manifold (SPM), $(h_{\tau}(M_t, \cdot))_{t\geq0}$, where $h_\tau$ is defined in \eqref{Eq_stoch_mnf} and the path-dependent coefficient $M_t$ defined in  \eqref{Eq_def_OU}; see Remark \ref{Remark_OU}.

The parameter regime in Table \ref{Table_Case}, with $\epsilon >1$, exemplifies this inverted timescale separation. In this regime, the variable $z(t)$ evolves more slowly than $x(t)$ and $y(t)$ (see Fig.~\ref{Fig_sol}), a scenario that contradicts the typical slow-fast timescale separation assumed in Cases I and II (Section \ref{Sec_numerics_1}) and in conventional slow manifold theory.
\begin{table}[h] 
	\caption{Parameter regime: Case III}
	\label{Table_Case}
	\centering
	\begin{tabular}{llllll}
		\toprule\noalign{\smallskip}
		& $\lambda$ & $f$ &  $\gamma$  & $\epsilon$ & $\sigma$\\ 
		\noalign{\smallskip}\hline\noalign{\smallskip}
		Case III  & $10^{-3}$  & $10$ & $1$ & $10$ & $0.3$ \\ 
		\noalign{\smallskip} \bottomrule 
	\end{tabular}
\end{table}

\begin{figure}[htbp]
	\centering
	\includegraphics[width=.95\textwidth,height=.4\textwidth]{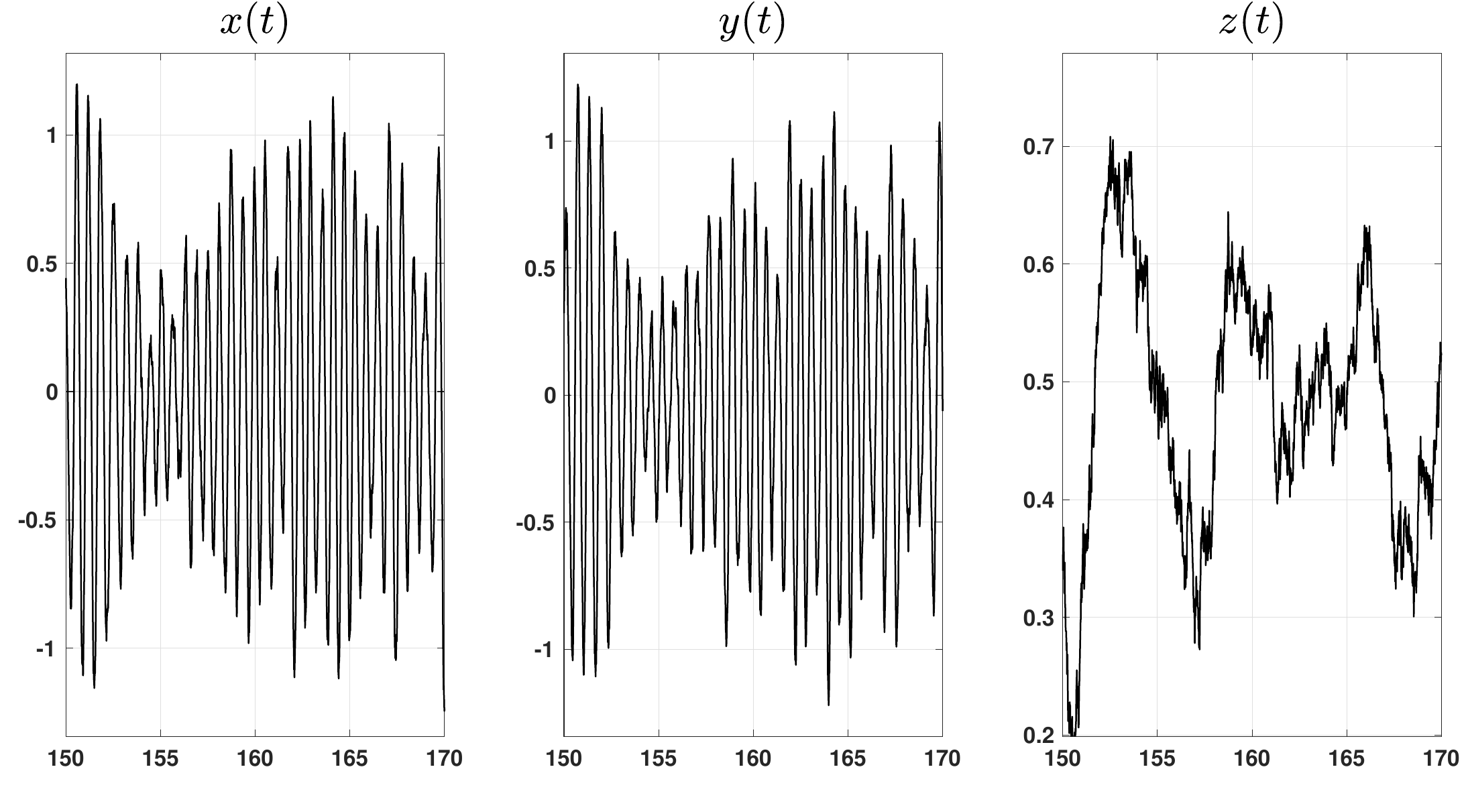}
	\caption{\footnotesize {\bf Solution to Eq.~\eqref{Eq_3DHopf_polar_base_stoch} for Case III}. Here, the variable $z(t)$ to parameterize evolves on a slower timescale than $x(t)$ and $y(t)$, a completely reverse situation compared to Cases I and II considered above.}
	\label{Fig_sol} 
\end{figure}

Following the approach in \cite{CLM19_closure,chekroun2023optimal}, we minimize the time-averaged parameterization defect $\overline{|z(t)^2- h_\tau(M_t,x(t),y(t))|^2}$ over a training path and characteristic timescale, where  $(x(t),y(t),z(t))$ represents data from the full system. This optimization procedure yields an SPM, $(h_{\tau}(M_t, \cdot))_{t\geq0}$, that effectively captures the underlying fluctuations of the $z$-variable in terms of the $x$- and $y$-variables. Remarkably, this predictive capability extends to out-of-sample paths. Figure \ref{Fig_PDF-section3} visually illustrates this by showcasing for an out-of-sample path, three snapshots of the SPM (red and bluish curves) spanning a diverse range of states within the $(r,z)$-plane of the original system as time evolves.

As illustrated in Figure \ref{Fig_PDF-section3}, the system exhibits a cluster of states within the $(r,z)$-plane characterized by significant variance along the $r$-direction. This radial elongation leads to pronounced transversality of the dynamics to the deterministic slow manifold (yellow curve). As highlighted in \cite{chekroun2021stochastic}, this high degree of transversality poses a significant challenge for the deterministic slow manifold to accurately parameterize the $z$-variable as a function of $x$ and $y$.

In contrast, the stochastic nature of the SPM, encapsulated by the coefficient $M_t$, and its inherent flexibility in shape, controlled by the parameter $\tau$, allow it to effectively capture these transversal fluctuations. This enhanced parameterization performance is evident in the normalized parameterization defects: $Q=3.9\times 10^{-1}$ for the SPM and $Q=1.1867$ ($>1$) for the deterministic slow manifold. Here, $Q=\overline{|z(t)^2- h_\tau(M_t,x(t),y(t))|^2}/\overline{|z^2(t)|}$ is a metric analogous to that defined in \eqref{Eq_param_def_cont}.
It assesses the performance of the stochastic parameterization, $h_\tau$, relative to the deterministic slow manifold, $h$, from a practical perspective, i.e., by directly evaluating the performance on actual data.

Further evidence of the SPM's superior performance is provided in Figure \ref{Fig_PDF-SPM}. The right panel of this figure compares the radial probability density functions (PDFs) of the original system (Eq.~\eqref{Eq_3DHopf_polar_base_stoch}) with those of the reduced systems. Notably, the PDF predicted by the reduced system based on the SPM (Eq.~\eqref{Eq_2DHopf_polar_base_aug}) exhibits reasonable agreement with the original system's PDF. In contrast, the PDF predicted by the reduced system based on the deterministic slow manifold (Eq.~\eqref{Eq_2DHopf_polar_base}) significantly deviates from the true distribution. These results demonstrate the enhanced ability of the SPM reduced system to capture the large radial excursions observed in the original system.

\begin{figure}[htbp]
	\centering
	\includegraphics[width=.6\textwidth,height=.5\textwidth]{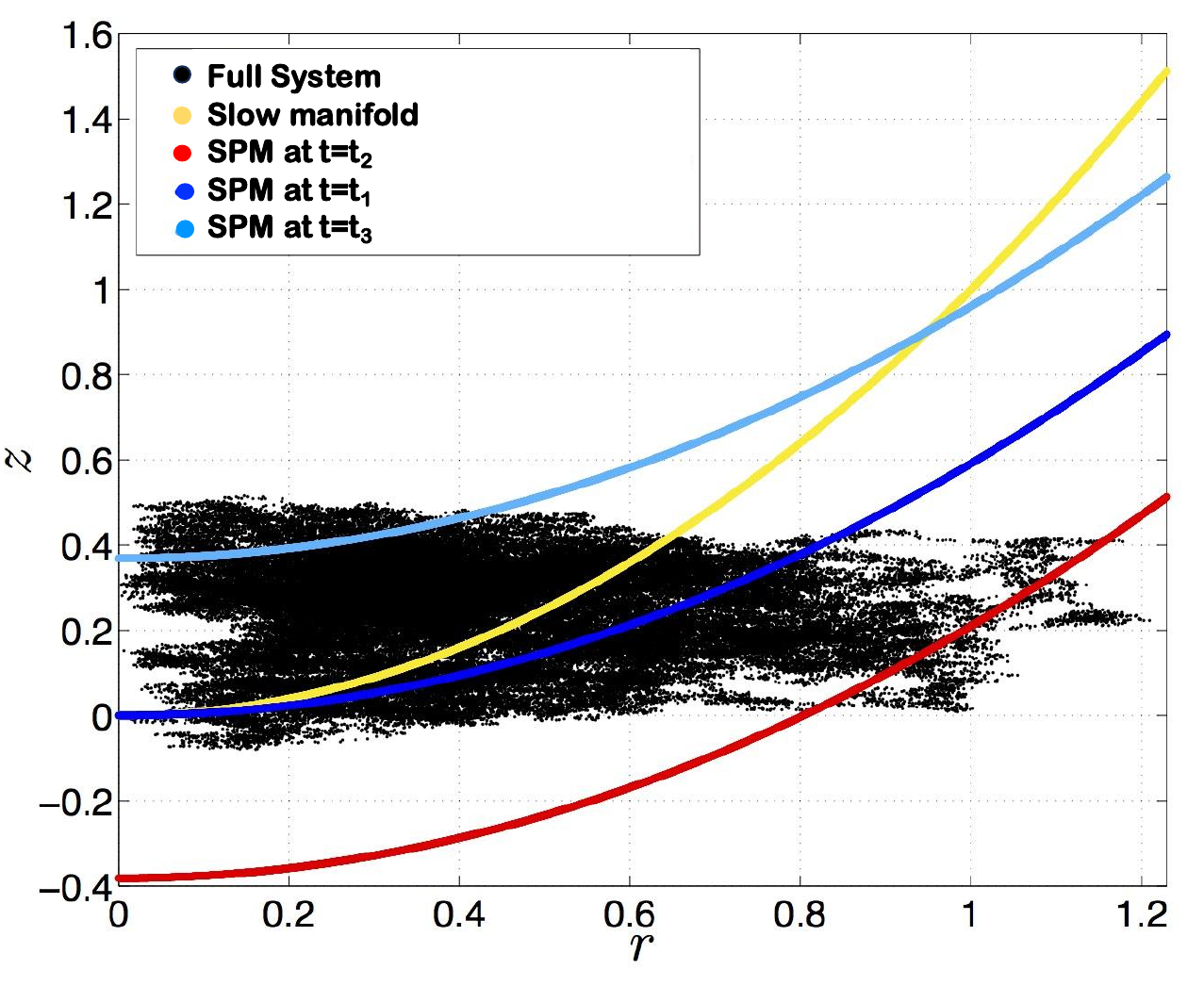}
	\caption{\footnotesize {\bf Capturing system variability}. This figure depicts three snapshots of the Stochastic Parameterizing Manifold (SPM), $(h_{\tau}(M_t, \cdot))_{t\geq0}$ (with $M_t$ defined in \eqref{Eq_def_OU}), represented by cyan, blue, and red curves. For comparison, the deterministic slow manifold is shown in yellow. Driven by the stochastic process $M_t$, the SPM dynamically evolves, effectively capturing a significant portion of the state space occupied by the full system, which exhibits substantial variance along the $r$-direction.}
	\label{Fig_PDF-section3} 
\end{figure}

\begin{figure}[htbp]
	\centering
	\includegraphics[width=1\textwidth,height=.5\textwidth]{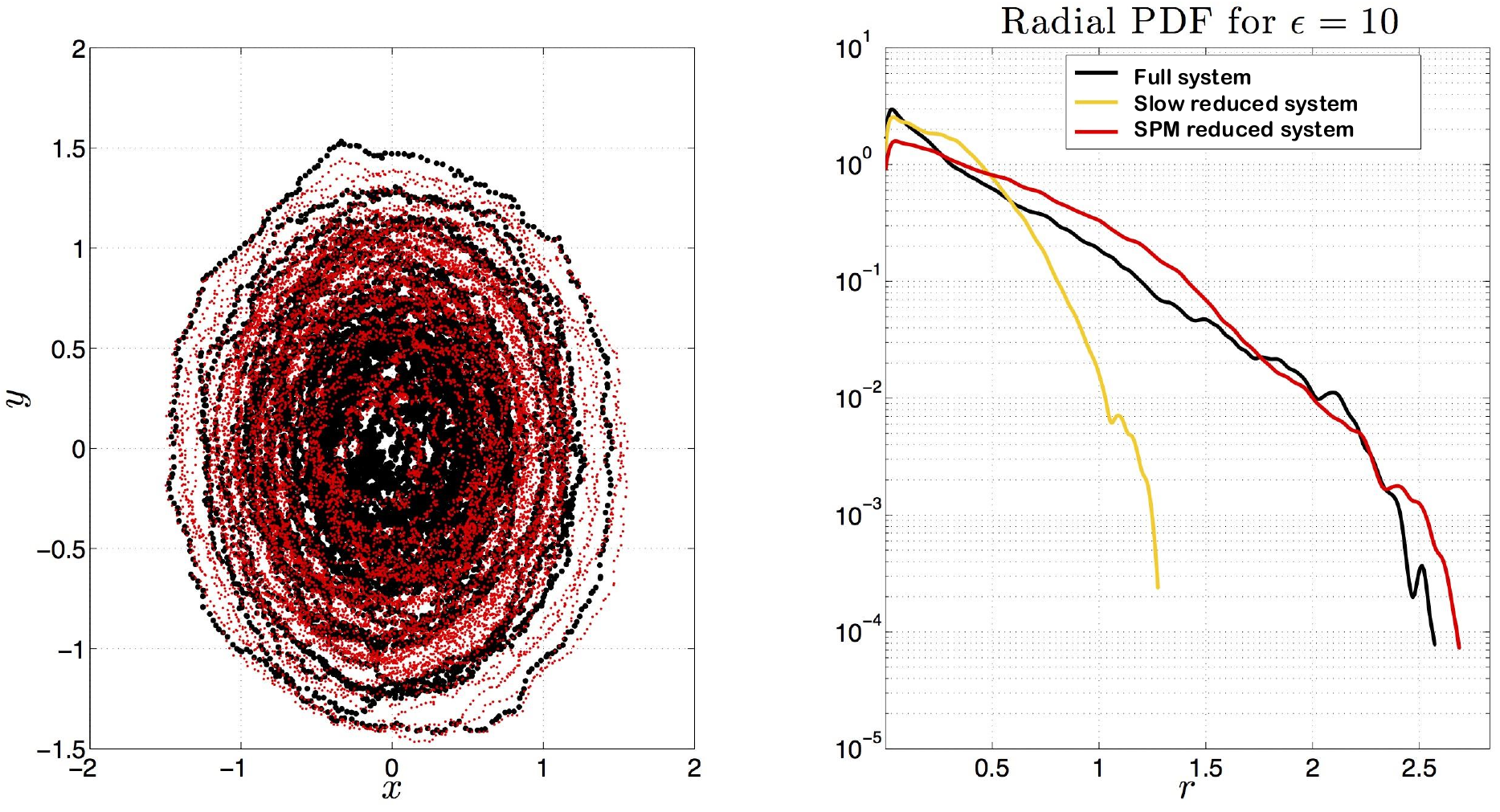}
	\caption{\footnotesize {\bf SPM reduced system skills.} The SPM reduced system shows a better ability than the slow reduced system in capturing the large radial excursions.}
	\label{Fig_PDF-SPM}
\end{figure}

\section*{Acknowledgments}
This work was completed while Christian J. Pangerl (CJP) was affiliated with Imperial College London.
CJP, JSWL and MR have received funding from the European Union’s Horizon 2020 research and innovation
programme under the Marie Sk\l odowska-Curie grant agreement No 643073. This work has been also supported by the Office of Naval Research (ONR) Multidisciplinary University Research Initiative (MURI) grant N00014-20-1-2023, by the National Science Foundation grant DMS-2407484, and by the European Research Council (ERC) under the European Union’s Horizon 2020 research and innovation program (Grant Agreement No. 810370).

\section*{Data Availability}
Data sets generated during the current study are available from the corresponding author on reasonable request.


\appendix
{\small

\section{Exponential moment bounds for systems \eqref{Eq_2DHopf_polar_base} and \eqref{Eq_3DHopf_polar_base_transf} in Section \ref{Sec_slow_pm}}
\label{Appendix_lemmas}

\bl \label{Lemma_moment_bounds}
	For $\hat{R}_0 = (\hat{r}_0, \hat{\theta}_0)$ in $\hat{\mathcal{M}}$ let $(\hat{R}_t^{\hat{R}_0})_{t\geq 0}$ with 
	\bes
	\hat{R}_t^{\hat{R}_0} = (\hat{r}_t, \hat{\theta}_t),\,\,\, t \geq 0,
	\ees
	be the unique strong solution to the reduced system \eqref{Eq_2DHopf_polar_base}. Furthermore, let $(\tilde{R}_t^{\tilde{R}_0})_{t\geq 0}$ with 
	\bes
	\tilde{R}_t^{\tilde{R}_0} = (\tilde{r}_t, \tilde{\theta}_t, \tilde{z}_t),\,\,\, t \geq 0,
	\ees
	be the unique strong solution to the transformed  system \eqref{Eq_3DHopf_polar_base_transf} for an initial value 
	\bes
	\tilde{R}_0 = (\tilde{r}_0, \tilde{\theta}_0, \tilde{z}_0)\,\,\, \text{in}\,\,\, \mathcal{M}.
	\ees
	For constants $\alpha, \beta, \eta>0$, we introduce the function $V$ by
	\bea \label{Eq_V_girs}
	V:\,& \mathbb{R}_+\times \mathbb{R} \times \mathbb{R}_+ \longrightarrow \mathbb{R},\\
	&(\tilde{r}, \tilde{z}, \hat{r}) \mapsto  \alpha \tilde{r}^2 + \beta \tilde{z}^2 + \alpha \hat{r}^2.
	\eea
	Furthermore, let
	\bes
	Q(\tilde{r}, \tilde{z}, \hat{r}) = \exp\left(V(\tilde{r}, \tilde{z}, \hat{r})\right).
	\ees
	Then, there exist parameter choices $\alpha, \beta, \eta>0$ such that for all $T>0$ there is a constant 
	\bes
	C = C(T, \hat{r}_0,\tilde{r}_0,\tilde{z}_0) > 0,
	\ees
	 such that  
	\be \label{Eq_to_show_B1_1_new}
	\mathbb{E}_{\mathbb{P}}\left[Q(\tilde{r}_t, \tilde{z}_t, \hat{r}_t) \right] \leq C,\,\,\,\text{for all}\,\,\,t \in [0,T].
	\ee
\el
\begin{proof}

In the following, our goal is to show that there are choices of parameters $\alpha$  $\beta$ and $\eta$ in the definition of $V$ in \eqref{Eq_V_girs} such that for every $T>0$ there exists a constant $C_T>0$, which is also dependent on the initial values $\tilde{r}_0$, $\tilde{z}_0$ and $\hat{r}_0$ such that
\be \label{Eq_to_show_B1}
\mathbb{E}_{\mathbb{P}}\left[Q(\tilde{r}_t, \tilde{z}_t, \hat{r}_t) \right] \leq C_T,\,\,\,\text{for all}\,\,\,t \in [0,T].
\ee
 Let $T>0$. In order to show inequality \eqref{Eq_to_show_B1} we apply It\^{o}'s formula (e.g. \cite[Theorem 5.1, p. 66-67]{watanabe_ikeda}) to the stochastic process $(Q(\tilde{r}_t, \tilde{z}_t, \hat{r}_t))_{t\in[0,T]}$ and obtain $\mathbb{P}$-a.s. for all $t$ in $[0,T]$
\bea \label{Eq_girs_nov_1}
\d Q(\tilde{r}_t, \tilde{z}_t, \hat{r}_t) &= \left\langle \nabla Q(\tilde{r}_t, \tilde{z}_t, \hat{r}_t),  \left(
\begin{array}{c}
	\d \tilde{r}_t\\
	\d \tilde{z}_t\\
	\d \hat{r}_t\\
\end{array}
\right)\right\rangle + \frac{1}{2} \left\langle D^2Q(\tilde{r}_t, \tilde{z}_t, \hat{r}_t) \left(
\begin{array}{c}
	\d \tilde{r}_t\\
	\d \tilde{z}_t\\
	\d \hat{r}_t\\
\end{array}
\right),\left(
\begin{array}{c}
	\d \tilde{r}_t\\
	\d \tilde{z}_t\\
	\d \hat{r}_t\\
\end{array}
\right) \right \rangle \\
&= A + B.
\eea
For term $A$ in \eqref{Eq_girs_nov_1} we have that
\bea \label{Eq_girs_A}
A &= \left \langle 2Q(\tilde{r}_t, \tilde{z}_t,\hat{r}_t) \left(
\begin{array}{c}
	\alpha \tilde{r}_t\\
	\beta \tilde{z}_t\\
	\eta \hat{r}_t\\
\end{array}
\right),
\left(
\begin{array}{c}
	- q\tilde{r}_t - \gamma \tilde{r}_t \tilde{z}_t + \frac{\sigma^2}{2 \tilde{r}_t} +(q+\lambda)\hat{r}_t \\
	-\frac{1}{\epsilon}\left(\tilde{z}_t - \tilde{r}^2_t \right)\\
	\lambda \hat{r}_t -\gamma\hat{r}_t^3 + \frac{\sigma^2}{2\hat{r}_t} \\
\end{array}
\right)\d t \right \rangle \\
&+  \left \langle 2Q(\tilde{r}_t,\tilde{z}_t,\hat{r}_t)\left(
\begin{array}{c}
	\alpha \tilde{r}_t\\
	\beta \tilde{z}_t\\
	\eta \hat{r}_t \\
\end{array}
\right),
\left(
\begin{array}{c}
	\sigma \d W^r_t\\
	\frac{\sigma}{\sqrt{\epsilon}}\d W^3_t\\
	\sigma\d W^r_t\\
\end{array}
\right) \right \rangle.
\eea
Evaluating the scalar products in \eqref{Eq_girs_A} yields
\bea \label{Eq_girs_A_2}
A &= 2Q(\tilde{r}_t,\tilde{z}_t,\hat{r}_t)\Big(-\alpha q \tilde{r}^2_t -\frac{ \beta}{\epsilon} \tilde{ z}^2_t + \left(\frac{\beta}{\epsilon} - \alpha \gamma \right) \tilde{ z}_t\tilde{r}^2_t + \alpha(q+\lambda)\tilde{r}_t\hat{r}_t  \\
&+ \lambda\eta \hat{r}_t^2 - \gamma \eta \hat{r}_t^4 +\frac{\alpha \sigma^2}{2} + \frac{\eta \sigma^2}{2}
\Big)\d t \\
&+ \d N^1_t  + \d N^2_t + \d N^3_t,
\eea
where the stochastic processes $N_1=(N^1_t)_{t\in[0,T]}$, $N_2=(N^2_t)_{t\in[0,T]}$ and $N_3=(N^3_t)_{t\in[0,T]}$ are given by
\bea \label{Eq_girs_loc_mart}
N^1_t &= \int_0^t2\alpha \sigma Q(\tilde{r}_s,\tilde{z}_s,\hat{r}_s) \tilde{r}_s\d W^r_s,\,\,\,t\in [0,T],\\
N^2_t &= \int_0^t \frac{2\beta \sigma}{\sqrt{\epsilon}}Q(\tilde{r}_s,\tilde{z}_s,\hat{r}_s)\tilde{z}_s \d W^3_s,\,\,\,t\in [0,T],\\
N^3_t &= \int_0^t 2 \eta \sigma Q(\tilde{r}_s,\tilde{z}_s,\hat{r}_s) \hat{r}_s \d W^r_s,\,\,\,t\in [0,T].
\eea
Since the stochastic processes $(2\alpha \sigma Q(\tilde{r}_t,\tilde{z}_t,\hat{r}_t) \tilde{r}_t)_{t\in [0,T]}$, $(\frac{2\beta \sigma}{\sqrt{\epsilon}}Q(\tilde{r}_t,\tilde{z}_t,\hat{r}_t)\tilde{z}_t)_{t\in[0,T]}$ and $(2\eta \sigma Q(\tilde{r}_t,\tilde{z}_t,\hat{r}_t) \hat{r}_t)_{t\in[0,T]}$ are continuous $\mathbb{P}$-a.s., it follows that
\bea \label{Eq_P_1_girs}
\mathbb{P}\left(\int_0^T \left|2 \alpha \sigma Q(\tilde{r}_t,\tilde{z}_t,\hat{r}_t) \tilde{r}_t\right|^2\d t < \infty \right) &= \mathbb{P}\left(\int_0^T \left|\frac{2\beta \sigma}{\sqrt{\epsilon}}Q(\tilde{r}_t,\tilde{z}_t,\hat{r}_t)\tilde{z}_t\right|^2\d t < \infty\right) \\
&= \mathbb{P}\left(\int_0^T \left|2\eta \sigma Q(\tilde{r}_t,\tilde{z}_t,\hat{r}_t) \hat{r}_t\right|^2 \d t<\infty \right) \\
&= 1.
\eea
Let the sequences of stopping times\footnote{With the convention $\inf \emptyset = \infty$.} $(\tau^1_n)_{n\in\mathbb{N}}$, $(\tau^2_n)_{n\in\mathbb{N}}$ and $(\tau^3_n)_{n\in\mathbb{N}}$ be defined as 
\bea \label{Eq_stopping_times_hilf}
\tau^1_n &= \inf\left\{t\in[0,T]: \int_0^t \left|2 \alpha \sigma Q(\tilde{r}_s,\tilde{z}_s,\hat{r}_s) \tilde{r}_s\right|^2\d s > n  \right\}\wedge T, \\
\tau^2_n &= \inf\left\{t\in[0,T]: \int_0^t \left|\frac{2\beta \sigma}{\sqrt{\epsilon}}Q(\tilde{r}_s,\tilde{z}_s,\hat{r}_s)\tilde{z}_s\right|^2 \d s > n \right\} \wedge T, \\
\tau^3_n &= \inf\left\{t \in [0,T]: \int_0^t \left|2\eta \sigma Q(\tilde{r}_s,\tilde{z}_s,\hat{r}_s) \hat{r}_s\right|^2
\ \d s > n \right\} \wedge T.
\eea
From \eqref{Eq_P_1_girs} it follows that $(\tau^i_n)_{n\in\mathbb{N}}$, $i=1,2,3$, are monotonically increasing in $n$ and 
\be \label{Eq_B1_tau_conv}
\lim_{n\rightarrow \infty} \tau^i_n = T, \,\,\, \mathbb{P}-a.s.,\,\,\,\text{for}\,\,\,i=1,2,3.
\ee
Therefore, the stopped processes $(N^i_{t\wedge \tau^i_n})_{t\in[0,T]}$, $i=1,2,3$, are martingales for all $n$ in $\mathbb{N}$. In addition, let $(\tau_n)_{n\in\mathbb{N}}$ be defined as
\be \label{Eq_B1_tau_n}
\tau_n = \tau^1_n \wedge \tau^2_n \wedge \tau^3_n,\,\,\,\text{for}\,\,\,n\in\mathbb{N}.
\ee
Note that the stopping times $(\tau_n)_{n\in\N}$ are again monotonically increasing and it holds that
	\be \label{Eq_conv_stop_expmom}
	\lim_{n\rightarrow \infty} \tau_n = T,\,\,\,\W-a.s.
	\ee
The stochastic processes $(N^{i}_{t\wedge \tau_n})_{t\in[0,T]}$, $i=1,2,3$, are also martingales for all $n$ in $\mathbb{N}$. In other words the processes $(N^i_t)_{t\in[0,T]}$, $i=1,2,3$, are local martingales (e.g.~\cite[Proposition 2.~24, p.~147]{bm_sc}P).

For the term $B$ in \eqref{Eq_girs_nov_1} we obtain for all $t$ in $[0,T]$ that
\bea \label{Eq_girs_B}
B = \left( 
2\alpha^2 \sigma^2 \tilde{r}_t^2 + 4 \alpha \eta \sigma^2 \tilde{r}_t \hat{r}_t + 2 \eta^2 \sigma^2 \hat{r}_t^2 +  \frac{2\sigma^2}{\epsilon}\beta^2 \tilde{z}_t^2 + \sigma^2 \left(\alpha + \eta + \frac{\beta}{\epsilon}\right) \right) Q(\tilde{r}_t,\tilde{z}_t,\hat{r}_t) \d t.
\eea
Combining $A$ in \eqref{Eq_girs_A_2} and $B$ in \eqref{Eq_girs_B} yields
\bea \label{Eq_girs_A_B_1}
A+B = Q(\tilde{r}_t,\tilde{z}_t,\hat{r}_t)\bigg(&\left(2\alpha^2\sigma^2 - 2q \alpha  \right) \tilde{r}_t^2 + \left(\frac{2\sigma^2}{\epsilon}\beta^2 - \frac{2\beta}{\epsilon}\right)\tilde{z}_t^2 
+\left(\frac{2\beta}{\epsilon}-2\alpha\gamma\right) \tilde{z}_t \tilde{r}_t^2\\
&+\left(4\alpha \eta \sigma^2 + 2\alpha(q+\lambda)\right)\tilde{r}_t\hat{r}_t +(2\lambda \eta \hat{r}_t^2 - 2\gamma \eta \hat{r}^4_t) +2\eta^2\sigma^2 \hat{r}_t^2+C_1\bigg)\d t \\
&+ \d N^1_t  + \d N^2_t + \d N^3_t,
\eea 
where 
\bes
C_1 = \sigma^2(\alpha + \eta) + \sigma^2\left(\alpha + \eta + \frac{\beta}{\epsilon}\right).
\ees
By choosing 
\be \label{Eq_param_choices_1}
\beta = \frac{1 }{\sigma^2}\,\,\,\text{and}\,\,\,\alpha = \frac{1}{\sigma^2\epsilon \gamma},
\ee
we achieve that the second and third term in the RHS of \eqref{Eq_girs_A_B_1} vanish.
Inserting \eqref{Eq_girs_A_B_1} into \eqref{Eq_girs_nov_1} with parameter choices as in \eqref{Eq_param_choices_1} and introducing the stopping time $\tau_n$ in \eqref{Eq_B1_tau_n} for an arbitrary $n$ in $\mathbb{N}$ yields $\mathbb{P}$-a.s. for all $t$ in $[0,T]$  
\bea \label{Eq_girs_main_2}
Q(\tilde{r}_{t\wedge \tau_n}, \tilde{z}_{t\wedge \tau_n}, \hat{r}_{t\wedge \tau_n}) &= Q(\tilde{r}_{0},\tilde{z}_{0}, \hat{r}_{0}) \\
&+ \int_0^{t\wedge \tau_n}	Q(\tilde{r}_{s}, \tilde{z}_{s}, \hat{r}_{s})\Big( \left(2\alpha^2\sigma^2 - 2q \alpha \right) \tilde{r}_s^2 +\left(4\alpha \eta \sigma^2 + 2\alpha(q+\lambda)\right)\tilde{r}_s\hat{r}_s \\
&+(2\lambda \eta \hat{r}_s^2 - 2\gamma \eta \hat{r}^4_s) +2\eta^2\sigma^2 \hat{r}_s^2+C_1 \Big) \d s \\
&+N^1_{t\wedge \tau_n} + N^2_{t\wedge \tau_n} + N^3_{t\wedge \tau_n}.	 
\eea
We set
\be \label{Eq_choice_eta}
\eta = \alpha.
\ee
Furthermore, a straight-forward computation shows that for every $\rho>2\lambda \eta$ there exists a constant $C_{\rho}>0$ such that
\be \label{Eq_girs_inequ_1}
2\lambda \eta r^2 - 2\gamma \eta r^4 \leq -\rho r^2 + C_{\rho},\,\,\,\text{for all}\,\,\,r>0,
\ee 
where $C_{\rho}$ is explicitly given by
\bes
C_{\rho} = \left(\frac{\rho}{2}-\lambda \eta\right)\left(\frac{\rho + 2\lambda \eta}{4\gamma \eta}\right).
\ees
Additionally, an application of Young's inequality yields $\W$-a.s. for all $t$ in $[0,T]$
\be \label{Eq_girs_ineq_2}
\left(4\alpha \eta \sigma^2 + 2\alpha(q+\lambda)\right)\tilde{r}_t\hat{r}_t \leq 2\alpha^2 \tilde{r}_t^2 + \frac{1}{2} (2\eta\sigma^2 + q +\lambda)^2 \hat{r}_t^2.
\ee
Using inequalities \eqref{Eq_girs_inequ_1} and \eqref{Eq_girs_ineq_2} in \eqref{Eq_girs_main_2} results $\W$-a.s. for all $t$ in $[0,T]$ in 
\bea \label{Eq_girs_main_3}
Q(\tilde{r}_{t\wedge \tau_n}, \tilde{z}_{t\wedge \tau_n}, \hat{r}_{t\wedge \tau_n}) &\leq Q(\tilde{r}_{0},\tilde{z}_{0}, \hat{r}_{0}) \\
&+ \int_0^{t\wedge \tau_n}	Q(\tilde{r}_{s}, \tilde{z}_{s}, \hat{r}_{s}) \tilde{r}_s^2  \left(2\alpha^2\sigma^2 - 2q \alpha + 2\alpha^2 \right) \d s \\
&+\int_0^{t\wedge \tau_n} Q(\tilde{r}_{s}, \tilde{z}_{s}, \hat{r}_{s}) \bigg( \Big(2\eta^2\sigma^2 - \rho + \frac{1}{2}(2\eta \sigma^2 + q + \lambda)^2\Big) \hat{r}_s^2 + C_{\rho} + C_1 \bigg) \d s \\
&+N^1_{t\wedge \tau_n} + N^2_{t\wedge \tau_n} + N^3_{t\wedge \tau_n}.	 
\eea
Recall our particular choices of $\alpha$ and $\eta$ in \eqref{Eq_param_choices_1} and \eqref{Eq_choice_eta}, respectively. By setting 
	\be \label{Eq_def_q}
	q = \frac{1}{\epsilon \gamma}\left(1 + \frac{1}{\sigma^2}\right) \,\,\,\text{and}\,\,\, \rho = \max\left\{2\eta^2\sigma^2 + \frac{1}{2}(2\eta \sigma^2 + q + \lambda)^2, 2\lambda \eta +1\right\},
	\ee
 due to cancellations in \eqref{Eq_girs_main_3}, we obtain $\W$-a.s. for all $t$ in $[0,T]$,
\bea \label{Eq_girs_main_4}
Q(\tilde{r}_{t\wedge \tau_n}, \tilde{z}_{t\wedge \tau_n}, \hat{r}_{t\wedge \tau_n}) \leq Q(\tilde{r}_{0},\tilde{z}_{0}, \hat{r}_{0}) 
&+ C_2\int_0^{t\wedge \tau_n}	Q(\tilde{r}_{s}, \tilde{z}_{s}, \hat{r}_{s}) \d s \\
&+N^1_{t\wedge \tau_n} + N^2_{t\wedge \tau_n} + N^3_{t\wedge \tau_n},	 
\eea
where 
\bes
C_2 = C_{\rho} + C_1.
\ees
Note that the process $(Q(\tilde{r}_t,\tilde{z}_t, \hat{r}_t))_{t\in [0,T]}$ is positive $\W$-a.s. for all $t$ in $[0,T]$.
	Therefore, we also have $\W$-a.s. for all $t$ in $[0,T]$ that
	\bea \label{Eq_girs_main_6}
	Q(\tilde{r}_{t\wedge \tau_n}, \tilde{z}_{t\wedge \tau_n}, \hat{r}_{t\wedge \tau_n}) \leq Q(\tilde{r}_{0},\tilde{z}_{0}, \hat{r}_{0}) 
	&+ C_2\int_0^{t}	Q(\tilde{r}_{s\wedge \tau_n}, \tilde{z}_{s\wedge \tau_n}, \hat{r}_{s\wedge \tau_n}) \d s \\
	&+N^1_{t\wedge \tau_n} + N^2_{t\wedge \tau_n} + N^3_{t\wedge \tau_n}.	 
	\eea
	For simplicity we define the family of process $(N_{t\wedge\tau_n})_{t\in [0,T]}$ for all $n$ in $\N$ by
	\be \label{Def_N_girs}
	N_{t\wedge \tau_n} = N^1_{t\wedge \tau_n} + N^2_{t\wedge \tau_n} + N^3_{t\wedge \tau_n},\,\,\,t\in [0,T].
	\ee
	Applying Gronwall's inequality to \eqref{Eq_girs_main_6} yields $\W$-a.s. for all $t$ in $[0,T]$ that
	\bea \label{Eq_girs_main_7}
	Q(\tilde{r}_{t\wedge \tau_n}, \tilde{z}_{t\wedge \tau_n}, \hat{r}_{t\wedge \tau_n}) &\leq \left(Q(\tilde{r}_{0},\tilde{z}_{0}, \hat{r}_{0}) +N_{t\wedge \tau_n}\right) \\
	&+ C_2\int_0^{t} e^{C_2(t-s)}\left(Q(\tilde{r}_{0},\tilde{z}_{0}, \hat{r}_{0})+N_{s\wedge \tau_n} \right)  	\d s. \\	 
	\eea
	Recall that the stochastic processes $(N^1_{t\wedge\tau_n})_{t\in [0,T]}$, $(N^2_{t\wedge\tau_n})_{t \in [0,T]}$ and $(N^3_{t\wedge\tau_n})_{t \in [0,T]}$ are martingales for all $n$ in $\mathbb{N}$. Hence,
	\be \label{Eq_N_null}
	\E_{\W}[N_{t\wedge \tau_n}] = 0,\,\,\,\text{for all}\,\,\,t\in[0,T]\,\,\,\text{and}\,\,\,n\in\N.
	\ee
	In addition, we have for all $n$ in $\mathbb{N}$ and $t$ in $[0,T]$ that
	\bea \label{Eq_girs_fub_ton}
	\int_0^t e^{-C_2s} \E_{\W}[\betrag{N_{s\wedge \tau_n}}] \d s &\leq \sum_{i=1}^3\int_0^t e^{-C_2s} \E_{\W}[\betrag{N^i_{s\wedge \tau_n}}]\d s \\
	&\leq
	\sum_{i=1}^3\int_0^t e^{-C_2s} \left(\E_{\W}[\betrag{N^i_{s\wedge \tau_n}}^2]\right)^{\frac{1}{2}}\d s &\text{(due to Cauchy-Schwarz ineq.)}\\
	& = \sum_{i=1}^3\int_0^t e^{-C_2s} \left(\E_{\W}\left[\left\langle N^i \right\rangle_{s\wedge \tau_n}\right]\right)^{\frac{1}{2}}\d s &\text{(due to It\^{o}'s isometry)}\\
	&\leq  \frac{3 \sqrt{n}}{C_2}.  &\text{(due to \eqref{Eq_stopping_times_hilf} and \eqref{Eq_B1_tau_n})}
	\eea
Therefore, Fubini's theorem is applicable and we obtain
	\bea \label{Eq_appl_tonelli}
	\E_{\W}\left[\int_0^t e^{-C_2 s} N_{s\wedge \tau_n}\d s\right] = \int_0^t e^{-C_2 s} \E_{\W}[N_{s\wedge \tau_n}] \d s = 0,\,\,\,\text{for all}\,\,\,t\in [0,T]\,\,\,\text{and}\,\,\,n \in \N.
	\eea
	By taking expectations on both sides of \eqref{Eq_girs_main_7} we obtain, using \eqref{Eq_N_null} and \eqref{Eq_appl_tonelli}, that
	\be \label{Eq_final_exp_Q}
	\mathbb{E}_{\mathbb{P}}\left[Q(\tilde{r}_{t\wedge \tau_n}, \tilde{z}_{t\wedge \tau_n}, \hat{r}_{t\wedge \tau_n})\right] \leq Q(\tilde{r}_{0},\tilde{z}_{0}, \hat{r}_{0}) e^{C_2T},\,\,\,\text{for all}\,\,\, t \in [0,T]\,\,\,\text{and}\,\,\,n\in \N.
	\ee
	Because of \eqref{Eq_conv_stop_expmom}, we obtain for all $t$ in $[0,T]$ by virtue of Fatou's lemma
	\bes
	\mathbb{E}_{\mathbb{P}}\left[Q(\tilde{r}_{t}, \tilde{z}_{t}, \hat{r}_{t})\right] \leq \liminf_{n\rightarrow \infty}\mathbb{E}_{\mathbb{P}}\left[Q(\tilde{r}_{t\wedge \tau_n}, \tilde{z}_{t\wedge \tau_n}, \hat{r}_{t\wedge \tau_n})\right] \leq Q(\tilde{r}_{0},\tilde{z}_{0}, \hat{r}_{0}) e^{C_2T}.
	\ees
	which completes the proof.
\end{proof}

\section{Exponential moment bounds for systems \eqref{Eq_2DHopf_polar_base_aug} and \eqref{Eq_3DHopf_polar_base_transf_aug} in Section \ref{Sec_stoch_pm}}
\label{Appendix_C}
\bl \label{Lemma_moment_bounds2}
	For $\hat{R}_0 = (\hat{r}_0, \hat{\theta}_0, M_0)$ in $\hat{\mathcal{M}}\times \mathbb{R}$ let $(\hat{R}_t^{\hat{R}_0})_{t\geq 0}$ with 
	\bes
	\hat{R}_t^{\hat{R}_0} = (\hat{r}_t, \hat{\theta}_t, M_t),\,\,\, t \geq 0,
	\ees
	be the unique strong solution to the reduced system \eqref{Eq_2DHopf_polar_base_aug}. Furthermore, let $(\tilde{R}_t^{\tilde{R}_0})_{t\geq 0}$ with 
	\bes
	\tilde{R}_t^{\tilde{R}_0} = (\tilde{r}_t, \tilde{\theta}_t, \tilde{z}_t, \tilde{M}_t),\,\,\, t \geq 0,
	\ees
	be the unique strong solution to the transformed  system \eqref{Eq_3DHopf_polar_base_transf_aug} for an initial value 
	\bes
	\tilde{R}_0 = (\tilde{r}_0, \tilde{\theta}_0, \tilde{z}_0, \tilde{M}_0)\,\,\, in\,\,\, \mathcal{M}\times \mathbb{R}.
	\ees
	For constants $\alpha_i>0$, $i=1,\dots,5$, we introduce the function $V$ by
	\bea \label{Eq_V_girs1}
	V:\,& \mathbb{R}_+\times \mathbb{R} \times \mathbb{R}_+ \times \mathbb{R} \longrightarrow \mathbb{R},\\
	&(\tilde{r}, \tilde{z}, \hat{r}, M) \mapsto  \alpha_1 \tilde{r}^2 + \alpha_2 \tilde{z}^2 + \alpha_3 \hat{r}^2 + \alpha_4 M^2 \hat{r}^2 + \alpha_5 M^2.
	\eea
	Furthermore, let
	\bes
	Q(\tilde{r}, \tilde{z}, \hat{r}, M) = \exp\left(V(\tilde{r}, \tilde{z}, \hat{r}, M)\right).
	\ees
	Then, there exist parameter choices $\alpha_i>0$, $i=1,\dots,5$, such that for all $T>0$ there is a constant $C = C(T, \hat{R}_0, \tilde{R}_0)$ such that  
	\be \label{Eq_to_show_B1_1_new2}
	\mathbb{E}_{\mathbb{P}}\left[Q(\tilde{r}_t, \tilde{z}_t, \hat{r}_t, M_t) \right] \leq C,\,\,\,\text{for all}\,\,\,t \in [0,T].
	\ee
\el
\begin{proof}
	
	In the following, our goal is to show that there are combinations of parameters $\alpha_i$, $i=1,\dots,5$, in the definition of $V$ such that for every $T>0$ there exists a constant $C_T>0$, which is also dependent on the initial values $\tilde{R}_0$ and $\hat{R}_0$ such that
	\be \label{Eq_to_show_B2}
	\mathbb{E}_{\mathbb{P}}\left[Q(\tilde{r}_t, \tilde{z}_t, \hat{r}_t, M_t) \right] \leq C_T,\,\,\,\text{for all}\,\,\,t \in [0,T].
	\ee
	In order to show inequality \eqref{Eq_to_show_B2} we apply It\^{o}'s formula to the stochastic process $(Q(\tilde{r}_t, \tilde{z}_t, \hat{r}_t, M_t))_{t\in[0,T]}$ and obtain $\mathbb{P}$-a.s. for all $t$ in $[0,T]$
	\bea \label{Eq_girs_nov_2}
	\d Q(\tilde{r}_t, \tilde{z}_t, \hat{r}_t, M_t) &= \left\langle \nabla Q(\tilde{r}_t, \tilde{z}_t, \hat{r}_t, M_t),  \left(
	\begin{array}{c}
		\d \tilde{r}_t\\
		\d \tilde{z}_t\\
		\d \hat{r}_t\\
		\d M_t \\
	\end{array}
	\right)\right\rangle + \frac{1}{2} \left\langle D^2Q(\tilde{r}_t, \tilde{z}_t, \hat{r}_t, M_t) \left(
	\begin{array}{c}
		\d \tilde{r}_t\\
		\d \tilde{z}_t\\
		\d \hat{r}_t\\
		\d M_t \\
	\end{array}
	\right),\left(
	\begin{array}{c}
		\d \tilde{r}_t\\
		\d \tilde{z}_t\\
		\d \hat{r}_t\\
		\d M_t \\
	\end{array}
	\right) \right \rangle \\
	&= A + B.
	\eea
	For term $A$ in \eqref{Eq_girs_nov_2} we have that
	\bea \label{Eq_girs_A2}
	A &= \left \langle 2Q(\tilde{r}_t, \tilde{z}_t,\hat{r}_t, M_t) \left(
	\begin{array}{c}
		\alpha_1 \tilde{r}_t\\
		\alpha_2 \tilde{z}_t\\
		\alpha_3 \hat{r}_t + \alpha_4M_t^2 \hat{r}_t\\
		\alpha_4 M_t \hat{r}^2_t +	\alpha_5 M_t\\
	\end{array}
	\right),
	\left(
	\begin{array}{c}
		- q\tilde{r}_t - \gamma \tilde{r}_t \tilde{z}_t + \frac{\sigma^2}{2 \tilde{r}_t} +(q+\lambda)\hat{r}_t -\gamma M_t \hat{r}_t\\
		-\frac{1}{\epsilon}\left(\tilde{z}_t - \tilde{r}^2_t \right)\\
		\lambda \hat{r}_t -\gamma c_{\tau}\hat{r}_t^3 + \frac{\sigma^2}{2\hat{r}_t} -\gamma M_t \hat{r}_t \\
		-\frac{1}{\epsilon}M_t^3 \\
	\end{array}
	\right)\d t \right \rangle \\
	&+  \left \langle 2Q(\tilde{r}_t,\tilde{z}_t,\hat{r}_t, M_t           )\left(
		\begin{array}{c}
		\alpha_1 \tilde{r}_t\\
		\alpha_2 \tilde{z}_t\\
		\alpha_3 \hat{r}_t + \alpha_4M_t^2 \hat{r}_t\\
		\alpha_4 M_t \hat{r}^2_t +	\alpha_5 M_t\\
	\end{array}
	\right),
	\left(
	\begin{array}{c}
		\sigma \d W^r_t\\
		\frac{\sigma}{\sqrt{\epsilon}}\d W^3_t\\
		\sigma\d W^r_t\\
		\frac{\sigma}{\sqrt{\epsilon}} \d W^4_t\\
	\end{array}
	\right) \right \rangle.
	\eea
	Evaluating the scalar products in \eqref{Eq_girs_A2} yields
	\bea \label{Eq_girs_A_3}
	A &= 2Q(\tilde{r}_t,\tilde{z}_t,\hat{r}_t,M_t)\Big(-\alpha_1 q \tilde{r}^2_t -\frac{ \alpha_2}{\epsilon} \tilde{ z}^2_t + \left(\frac{\alpha_2}{\epsilon} - \alpha_1 \gamma \right) \tilde{ z}_t\tilde{r}^2_t + \alpha_1(q+\lambda)\tilde{r}_t\hat{r}_t  \\
	&+ \lambda\alpha_3 \hat{r}_t^2 - \gamma c_{\tau}  \alpha_3 \hat{r}_t^4 \\
	&+M_t^2(\lambda \alpha_4 \hat{r}^2_t - \gamma c_{\tau} \alpha_4 \hat{r}^4_t)\\
	&+ \frac{\alpha_4\sigma^2}{2} M_t^2 - \frac{\alpha_5}{\epsilon}M_t^4 \\
	&-\alpha_1\gamma \tilde{ r}_t \hat{r}_t M_t - \alpha_3 \gamma M_t \hat{r}_t^2 - \alpha_4 \gamma \hat{r}_t^2 M_t^3 - \frac{\alpha_4}{\epsilon}\hat{r}^2_t M_t^4 \\ &+\frac{\alpha_1 \sigma^2}{2} + \frac{\alpha_3 \sigma^2}{2}
	\Big)\d t \\
	&+ \d N^1_t  + \d N^2_t + \d N^3_t  + \d N^4_t,
	\eea
	where the stochastic processes $N_1=(N^1_t)_{t\in[0,T]}$, $N_2=(N^2_t)_{t\in[0,T]}$,  $N_3=(N^3_t)_{t\in[0,T]}$ and $N_4=(N^4_t)_{t\in[0,T]}$ are given by
	\bea \label{Eq_girs_loc_mart2}
	N^1_t &= \int_0^t2\alpha_1 \sigma \tilde{r}_s Q(\tilde{r}_s,\tilde{z}_s,\hat{r}_s, M_s) \d W^r_s,\,\,\,t\in [0,T],\\
	N^2_t &= \int_0^t \frac{2\alpha_2 \sigma}{\sqrt{\epsilon}}\tilde{z}_sQ(\tilde{r}_s,\tilde{z}_s,\hat{r}_s, M_s) \d W^3_s,\,\,\,t\in [0,T],\\
	N^3_t &= \int_0^t 2 \sigma (\alpha_3 \hat{r}_s + \alpha_4 M_s^2 \hat{r}_s) Q(\tilde{r}_s,\tilde{z}_s,\hat{r}_s, M_s) \d W^r_s, \,\,\,t\in [0,T],\\
	N^4_t &= \int_0^t  \frac{2\sigma}{\sqrt{\epsilon}} (\alpha_4 M_s \hat{r}^2_s + \alpha_5 M_s) Q(\tilde{r}_s,\tilde{z}_s,\hat{r}_s, M_s) \d W^4_s,
	\,\,\,t\in [0,T].
	\eea
Similarly to \eqref{Eq_B1_tau_n} in the proof of Lemma \ref{Lemma_moment_bounds}, it follows that there exists an increasing sequence of stopping times $(\tau_n)_{n\in\mathbb{N}}$ such that the stopped processes $(N^i_{t\wedge \tau_n})_{t\in[0,T]}$, $i=1,2,3,4$, are martingales for all $n$ in $\mathbb{N}$.
For the term $B$ in \eqref{Eq_girs_nov_2} we obtain for all $t$ in $[0,T]$ that
	\bea \label{Eq_girs_B2}
	B &= \sigma^2\Big( \alpha_1 + \alpha_3 + \frac{\alpha_2}{\epsilon} + \frac{\alpha_5}{\epsilon} + 2 \alpha_1^2 \tilde{ r}_t^2 + \frac{2 \alpha_2^2}{\epsilon} \tilde{ z}_t^2 + \left(\alpha_4 + \frac{2 \alpha_5^2}{\epsilon}\right) M_t^2 \\
	&+ \left(2\alpha_3^2 + \frac{\alpha_4}{\epsilon}\right) \hat{r}_t^2 + 4 \left(\frac{\alpha_4 \alpha_5}{\epsilon} + \alpha_3 \alpha_4\right) M_t^2 \hat{r}^2_t \\
	&+ \frac{2\alpha_4^2}{\epsilon} M_t^2 \hat{r}_t^4 + 2\alpha_4^2 M_t^4 \hat{r}_t^2 \\
	&+ 4 \alpha_1 \alpha_3 \tilde{ r}_t \hat{r}_t + 4 \alpha_1\alpha_4 M_t^2 \hat{r}_t \tilde{ r}_t\Big) Q(\tilde{r}_t,\tilde{z}_t,\hat{r}_t, M_t) \d t.
	\eea
	We set
	\be \label{Eq_def_q2}
	q = \frac{1}{\epsilon \gamma}\left(1 + \frac{5}{2\sigma^2} \right), \,\,\,\alpha_1 = \frac{1}{\sigma^2\epsilon  \gamma},\,\,\, \alpha_2 = \frac{1}{\sigma^2},\,\,\, \alpha_3 = \alpha_5 = 1\,\,\,\text{and}\,\,\, \alpha_4 \in (0,\alpha_4^*),
	\ee
	where 
	\bes
	\alpha_4^* = \min\left\{\frac{\gamma c_{\tau} \epsilon}{\sigma^2}, \frac{2}{\epsilon(2\sigma^2 + 2\sigma^4 + \gamma)}\right\}.
	\ees
	Analogously to Eqns.~\eqref{Eq_girs_A_B_1}-\eqref{Eq_girs_main_3} we deduce, using the parameter choices in \eqref{Eq_def_q2}, the existence of a positive constant $C$ (independent of $n$ in $\mathbb{N}$) such that $\W$-a.s. for all $t$ in $[0,T]$ it holds that
	\bea \label{Eq_girs_main_42}
	Q(\tilde{r}_{t\wedge \tau_n}, \tilde{z}_{t\wedge \tau_n}, \hat{r}_{t\wedge \tau_n}, M_{t\wedge \tau_n}) \leq Q(\tilde{r}_{0},\tilde{z}_{0}, \hat{r}_{0}, M_0) 
	&+ C\int_0^{t\wedge \tau_n}	Q(\tilde{r}_{s}, \tilde{z}_{s}, \hat{r}_{s}, M_s) \d s \\
	&+N^1_{t\wedge \tau_n} + N^2_{t\wedge \tau_n} + N^3_{t\wedge \tau_n} + N^4_{t\wedge \tau_n}.	 
	\eea
Applying an analogous reasoning to that used in lines \eqref{Eq_girs_main_6}--\eqref{Eq_appl_tonelli} for the proof of Lemma \ref{Lemma_moment_bounds}, yields
	\bes
	\mathbb{E}_{\mathbb{P}}\left[Q(\tilde{r}_{t\wedge \tau_n}, \tilde{z}_{t\wedge \tau_n}, \hat{r}_{t\wedge \tau_n}, M_{t\wedge \tau_n})\right] \leq Q(\tilde{r}_{0},\tilde{z}_{0}, \hat{r}_{0}, M_0) e^{C T},\,\,\,\text{for all}\,\,\,n\in\N\,\,\,\text{and}\,\,\, t \in [0,T].
	\ees
	By virtue of Fatou's lemma we obtain for all $t$ in $[0,T]$ that
	\bes
	\mathbb{E}_{\mathbb{P}}\left[Q(\tilde{r}_{t}, \tilde{z}_{t}, \hat{r}_{t}, M_t)\right] \leq \liminf_{n\rightarrow \infty}\mathbb{E}_{\mathbb{P}}\left[Q(\tilde{r}_{t\wedge \tau_n}, \tilde{z}_{t\wedge \tau_n}, \hat{r}_{t\wedge \tau_n}, M_{t\wedge \tau_n})\right] \leq Q(\tilde{r}_{0},\tilde{z}_{0}, \hat{r}_{0}, M_0) e^{CT},
	\ees
	which yields \eqref{Eq_to_show_B2} and thus completes the proof.
\end{proof}

\section{Existence and uniqueness of invariant measures for systems \eqref{Eq_3DHopf_polar_base} and \eqref{Eq_2DHopf_polar_base}}
\label{App_ex_uniq_org}

\bl \label{Lemma_invm_ex_org}
	Recall the set $\mathcal{M}$ defined in \eqref{Def_M_hat_M}. For the original system \eqref{Eq_3DHopf_polar_base} we introduce the following Lyapunov function\footnote{The sub-level sets
		$K_a = \{(r,\theta,z) \in \mathcal{M}: V(r,\theta, z) \leq a\}$
		are compact in $\mathcal{M}$ for all $a >0$ and $V(r,\theta, z) \rightarrow \infty$ as $|(r,\theta,z)| \rightarrow \infty$.} $V$ on $\mathcal{M}$. Let
	\bea \label{Eq_V_lyap}
	V:\,& \mathcal{M} \longrightarrow \mathbb{R}_+,\\
	&(r,\theta, z) \mapsto  \alpha_1 r^2+ \alpha_2 \theta + (\alpha_3 z-p)^2 + 1,
	\eea
	with 
	\be \label{Eq_choices_of_param_invm_org}
	\alpha_1 = \frac{1}{\gamma \epsilon},\,\alpha_2, \alpha_3 = 1 \,\, \text{and} \,\,\, p = \frac{1 +2\lambda }{2\gamma}.
	\ee
	 Associated with $V$ we consider the space of functions $\mathcal{F}_V$ given by
					\bes
\mathcal{F}_V = \left\{ \varphi: \mathcal{M} \rightarrow \mathbb{R}  \Vert \varphi \,\,\,\text{is measurable}\,\,\,\text{and}\,\,\, \|\varphi\|_V < \infty \right\},
\ees
	where 
	\bes
	\|\varphi\|_V = \sup_{(r, \theta, z) \in \mathcal{M}} \frac{|\varphi(r,\theta,z)|}{V(r,\theta,z)}.
	\ees
	Let $P_t$ be the Markov semigroup induced by the original system \eqref{Eq_3DHopf_polar_base}.
	Then, there exists a unique ergodic invariant measure $\mu$ for $P_t$ with smooth density with respect to the Lebesgue measure on $\mathcal{M}$.
	In addition, it holds that 
	\be \label{Eq_second_moments_r}
	\int_{\mathbb{R}_+} r^2 \mu_r(\d r) < \infty,
	\ee
	where $\mu_r$ denotes the marginal of $\mu$ along the radial part in the phase space $\mathcal{M}$.
	Furthermore, there exist constants $C, \kappa > 0$ and an initial time $t_0>0$ such that for all initial values $(r,\theta, z)$ in $\mathcal{M}$ it holds that
	\be \label{Eq_exp_conv_1}
	\left|P_t \varphi (r,\theta,z) - \int_{\mathcal{M}} \varphi \d \mu \right| \leq C e^{-\kappa t}  V(r,\theta,z), \,\,\, \text{for all}\,\,\, t > t_0.
	\ee
	
\el
\begin{proof}
	{\bf Step 1: Existence of invariant measures.} Let $P_t$ denote the Markov semigroup for the original system \eqref{Eq_3DHopf_polar_base} as defined in \eqref{Eq_P_proof} above.
	Our goal is to apply \cite[Proposition 7.10, p. 99]{daprato_inf_dim}. For this purpose, 
	  note that for the choices of parameters $\alpha_1, \alpha_2, \alpha_3, p > 0$ in \eqref{Eq_choices_of_param_invm_org} the sub-level sets
	 \be \label{Eq_lemma_ex_invm_org_1}
	 K_a = \{(r,\theta,z) \in \mathcal{M}: V(r,\theta, z) \leq a\},\,\,\, a>0,
	 \ee
	 are compact in $\mathcal{M}$ for all $a >0$.
	 For $R_0 = (r_0,\theta_0,z_0)$ in $\mathcal{M}$ let $(R_t^{R_0})_{t\geq0}$ with
	 \bes
	 R_t^{R_0}= (r_t, \theta_t,z_t), \,\,\, t \geq 0,
	 \ees
	  be the unique strong solution to system \eqref{Eq_3DHopf_polar_base}. By applying It\^{o}'s formula to the process $(V(r_t, \theta_t, z_t))_{t\geq0}$ we obtain, using Gronwall's inequality, 
	 \be \label{Eq_exniq_ineq_2}
	 \mathbb{E}_{\mathbb{P}}[V(r_t,\theta_t,z_t)] \leq e^{-\rho t}V(r_0,\theta_0,z_0) + C, \,\,\, \text{for all } t \geq 0,
	 \ee
	 where 
	 \bes
	 \rho = \min\left\{1, \epsilon^{-1}\right\}\,\,\,\,\text{and}\,\,\,\, C = 2\pi + \frac{1}{\rho}\left(\left(2\alpha_1 + \frac{\alpha_3^2}{\epsilon}\right)\sigma^2 + \frac{p^2}{\epsilon}\right).
	 \ees
	 We note that inequality \eqref{Eq_exniq_ineq_2} can also be expressed in terms of the Markov semigroup $P_t$ as
	 \be \label{Eq_P_t_ineq_V}
	 P_tV(r_0,\theta_0,z_0) \leq e^{-\rho t} V(r_0,\theta_0,z_0) + C,\,\,\, \text{for all}\,\,\, (r_0,\theta_0,z_0) \in \mathcal{M}\,\,\,\text{and}\,\,\,t \geq0.
	 \ee
	 The existence of at least one invariant measure for $P_t$ now follows from \cite[Proposition 7.10, p. 99]{daprato_inf_dim}. 
	 \\
	{\bf Step 2: Strong Feller property.} We refer to \cite[Definition 5.2(ii), Chapter 5, p. 70]{daprato_inf_dim} for this notion.  To prove this property, we write the original system \eqref{Eq_3DHopf_polar_base} in the form
		\be{\label{Eq_3DHopf_polar_base_strat}}
		\d \begin{pmatrix} 
		r\\
		\theta\\
		z
		\end{pmatrix} = V_0(r,\theta, z)\d t + V_1\d W^r_t + V_2(r)\d W^{\theta}_t + V_3\d W^3_t,
		\ee
		where 
		\bes
		V_0(r,\theta,z) = \begin{pmatrix} 
			\lambda r -\gamma r z + \frac{\sigma^2}{2 r}\\
			f\\
			-\frac{1}{\epsilon}\left(z -r^2\right)
		\end{pmatrix}, 
		V_1 = \sigma \begin{pmatrix} 
			1\\
			0\\
			0
		\end{pmatrix},
	 	V_2(r) = \sigma \begin{pmatrix} 
	 		0\\
	 		\frac{1}{r}\\
	 		0
	 	\end{pmatrix},
		V_3 = \frac{\sigma}{\sqrt{\epsilon}} \begin{pmatrix} 
			0\\
			0	\\
			1
		\end{pmatrix}.
		\ees
		The generator of the Markov semigroup $P_t$ is then given by the second-order differential operator $L$ with
		\be \label{Def_K}
		L\varphi = V_0 \cdot \nabla \varphi + \frac{1}{2} \sum_{j = 1}^{3} V_{jj}^2 \partial_j^2 \varphi, \,\,\, \varphi \in C^2_{b}(\mathcal{M}),
		\ee
		where $V_{jj}$ denotes the $j$-th entry of $V_j$, $j=1,2,3$, and
		\bes
		\partial_1 = \partial_r, \,\, \partial_2 = \partial_{\theta},\,\, \partial_3 = \partial_z\,\,\,\text{and}\,\,\, \nabla = (\partial_1, \partial_2,\partial_3)^T.
		\ees
In the following we refer to $L$ as the (backward) Kolmogorov operator.

		We can verify H\"{o}rmander's bracket condition for $L$ in \eqref{Def_K} (see e.g. \cite{hoermander1967}), which implies that $L$ is hypoelliptic. This yields the strong Feller property for the Markov semigroup $P_t$; see e.g. \cite[Theorem 6.1, p. 188]{Pavliotis:2014}, \cite[Theorem 6.3, p. 26]{Hairer10convergenceof} or \cite[Theorem 5.6, Chapter VIII, p. 214]{Bass}.\\
	{\bf Step 3: Irreducibility.} We refer to \cite[Definition 5.2(iii), Chapter 5, p. 70]{daprato_inf_dim} for this notion. 
	Let $T>0$. To prove this property, we consider the following control system for $t$ in $[0,T]$ associated with the original system \eqref{Eq_3DHopf} in Cartesian coordinates. Let
	\begin{subnumcases}{\label{Eq_3DHopf_ctrl}}
	\frac{\d \phi_1}{\d t} =\lambda \phi_1- f \phi_2 -\gamma \phi_1 \phi_3 + \sigma u^1_t\\
	\frac{\d \phi_2}{\d t} =f \phi_1+\lambda \phi_2 -\gamma \phi_2 \phi_3 + \sigma u^2_t\\
	\frac{\d \phi_3}{\d t} =-\frac{1}{\epsilon}( \phi_3-\phi_1^2-\phi_2^2) + \frac{\sigma}{\sqrt{\epsilon}} u^3_t \label{Eq_3DHopf_z_ctrl},
	\end{subnumcases}
	where $u^i$, $i=1,2,3$, are elements of $C([0,T], \mathbb{R})$. For arbitrary $x_0$ and $x_1$ in $\mathbb{R}^3$ let 
	\bea \label{Eq_Phi_ctrl}
	\Phi: \,& [0,T] \longrightarrow \mathbb{R}^3,\\
	& t \mapsto  	 \begin{pmatrix} 
		 \phi_1(t)\\
		\phi_2(t)\\
		\phi_3(t)
	\end{pmatrix},
	\eea
such that $\Phi(0) = x_0$ and $\Phi(T) = x_1$. Moreover, we assume that $\Phi$ is an element of $C^1([0,T],\R^3)$.
	For $\Phi$ in \eqref{Eq_Phi_ctrl} we define $u^i$, $i=1,2,3$, for all $t$ in $[0,T]$ as
	\bea \label{Eq_u_ctrl}
	u^1_t & =\frac{1}{\sigma}\left(\frac{\d \phi_1}{\d t}- \lambda \phi_1 + f \phi_2 +\gamma \phi_1 \phi_3\right) \\
	u^2_t & =\frac{1}{\sigma}\left(\frac{\d \phi_2}{\d t} - f\phi_1 - \lambda \phi_2 + \gamma \phi_2 \phi_3\right) \\
	u^3_t & =\frac{\sqrt{\epsilon}}{\sigma}\left(\frac{\d \phi_3}{\d t} + \frac{1}{\epsilon}\phi_3 - \frac{1}{\epsilon} \phi_1^2 - \frac{1}{\epsilon} \phi_2^2 \right).
	\eea
	It then follows that $\Phi$ given in \eqref{Eq_Phi_ctrl} satisfies the control problem \eqref{Eq_3DHopf_ctrl} for the controls $u^i$, $i=1,2,3$, defined in \eqref{Eq_u_ctrl}.
	The irreducibility of the Markov semigroup $P_t$ then follows from \cite[Theorem 5.2]{stroock1972}.
	
	Doob's theorem \cite[Theorem 4.2.1(ii), Chapter 4, p. 43]{da_prato_zabczyk_1996} together with \cite[Theorem 3.2.6, Chapter 3, p. 28]{da_prato_zabczyk_1996} now ensures that there exists a unique ergodic invariant measure $\mu$ in $Pr(\mathcal{M})$ for the Markov semigroup $P_t$. Furthermore, due to \cite[(5.6) Theorem, Chapter VIII, p. 214]{Bass} and \cite[4.2.1(iii), Chapter 4, p. 43]{da_prato_zabczyk_1996}, $\mu$ has a smooth density with respect to the Lebesgue measure on $\mathcal{M}$.\\
	{\bf Step 4: Proof of \eqref{Eq_exp_conv_1}.} In order to show the exponential convergence to equilibrium  \eqref{Eq_exp_conv_1}, we rely on standard methods from stochastic analysis involving Lyapunov functions (e.g. ~\cite[Theorem 5]{Chekroun_al_RP2}). In that respect, we aim at proving, given the Kolmogorov operator defined in \eqref{Def_K}, and $V$ the function defined in \eqref{Eq_V_lyap}, that there exists a constant $a > 0$ such that 
	\be \label{Eq_KV_leq_aV}
	L V \leq a V. 
	\ee

The proof of this inequality consists of simple calculus. Indeed, for $(r,\theta,z)$ in $\mathcal{M}$, we obtain, by applying $L$ to $V$
	\bes
	LV(r,\theta, z) \leq \tilde{a} (V(r,\theta,z) - 1) + c,
	\ees
	with
	\bes
	\tilde{a} = \max\{ 1, \epsilon^{-1} \} \,\,\,\text{and}\,\,\, c = f +  \frac{2\sigma^2}{\gamma \epsilon } + \frac{1}{\epsilon}\left(\sigma^2 + \left(\frac{1+2\lambda}{2\gamma}\right)^2\right).
	\ees
	Hence,
	\bes
	LV(r,\theta,z) \leq a V(r,\theta,z), \,\,\,\text{with}\,\,\, a = \max\{\tilde{a}, c\},
	\ees
	which is \eqref{Eq_KV_leq_aV}.
	Inequalities \eqref{Eq_P_t_ineq_V} and \eqref{Eq_KV_leq_aV} together with the fact that $P_t$ is strongly Feller and irreducible yields inequality \eqref{Eq_exp_conv_1}; see \cite[Theorem 5]{Chekroun_al_RP2}. Finally, the moment bound \eqref{Eq_second_moments_r} follows from the estimate  \eqref{Eq_exp_conv_1}.
\end{proof}

\br \label{Remark_f_mom_bound}
Let $\mu$ be the unique ergodic invariant measure for the Markov semigroup $P_t$ of the original system \eqref{Eq_3DHopf_polar_base}. We remark that it also holds
\bes
\int_{\mathbb{R}_+ \times \mathbb{\mathbb{R}}} r^4 (z - r^2)^4 \mu_{r,z}(\d r, \d z)< \infty,
\ees
where $\mu_{r,z}$ denotes the marginal of $\mu$ along the $(r,z)$-plane (see \eqref{Eq_marginal_proof}).
\er

\bl \label{Lemma_invm_exuniq_red}
	Let $Q_t$ be the Markov semigroup induced by the reduced system \eqref{Eq_2DHopf_polar_base}. Then there exists a unique ergodic invariant measure $\nu$ in $Pr(\hat{\mathcal{M}})$ with smooth Lebesgue-density on $\hat{\mathcal{M}}$.
	In addition, a bound analogous to \eqref{Eq_exp_conv_1} holds for $Q_t$. 
	Furthermore,
	\bes
	\int_{\mathbb{R}_+} r^2 \nu_r(\d r) < \infty,
	\ees
	where $\nu_r$ denotes the marginal of $\nu$ along the radial direction in the phase space $\hat{\mathcal{M}}$.
\el
\begin{proof}
The proof  proceeds analogously to that  of Lemma \ref{Lemma_invm_ex_org}. 
\end{proof}

\section{Existence and uniqueness of invariant measures for systems \eqref{Eq_3DHopf_polar_base_aug} and \eqref{Eq_2DHopf_polar_base_aug}}
\label{Appendix_E}
\bl \label{Lemma_invm_ex_org_aug}
	There exists a unique ergodic invariant measure $\bar{\mu}$ in $Pr(\mathcal{M}\times \mathbb{R})$ for the augmented original system \eqref{Eq_3DHopf_polar_base_aug}. Moreover, $\bar{\mu}$ has a smooth density with respect to the Lebesgue measure on $\mathcal{M}\times \mathbb{R}$. In addition, a bound analogous to \eqref{Eq_exp_conv_1} holds.
\el
\begin{proof}
	{\bf Step 1: Existence of invariant measures.}
	The proof is analogous to the first step used in the proof of Lemma \eqref{Lemma_invm_ex_org}. 
	We define the function
	\bea \label{Eq_V_lyap_stoch}
	V:\,& \mathcal{M} \times \mathbb{R} \longrightarrow \mathbb{R}_+,\\
	&(r,\theta, z, M) \mapsto  \alpha_1 r^2+ \alpha_2 \theta + (\alpha_3 z-p)^2 + \alpha_4 M^2 +1,
	\eea
	where 
	\be \label{Eq_choices_of_param_invm_org_aug}
	\alpha_1 = \frac{1}{\gamma \epsilon},\,\alpha_2, \alpha_3, \alpha_4 = 1\,\, \text{and} \,\, p = \frac{1 +2\lambda }{2 \gamma}.
	\ee
	Note that all sub-level sets of $V$ in \eqref{Eq_V_lyap_stoch} are compact in $\mathcal{M}\times \mathbb{R}$ for the choice of parameters in \eqref{Eq_choices_of_param_invm_org_aug}. Using It\^{o}'s formula and Gronwall's inequality we obtain that
	\bes
	\mathbb{E}_{\mathbb{P}}[V(r_t,\theta_t,z_t,M_t)] \leq e^{-\rho t} V(r_0,\theta_0,z_0,M_0) + C,
	\ees
	where
	\bes
	\rho = \min\{1, \epsilon^{-1}\} \,\,\, \text{and}\,\,\, C = 2\pi + \frac{1}{\rho}\left(\left(2 \alpha_1 + \frac{\alpha_3^2}{\epsilon} + \frac{\alpha_4}{\epsilon} \right)\sigma^2 + \frac{2\alpha_4}{\epsilon}+ \frac{p^2}{\epsilon} \right).
	\ees
	The existence of at least one invariant measure for the Markov semigroup $P_t$ associated with the augmented original system \eqref{Eq_3DHopf_polar_base_aug} follows then from \cite[Proposition 7.10, p. 99]{daprato_inf_dim}. \\
	{\bf Step 2: Strong Feller property and irreducibility.} The strong Feller property and irreducibility for $P_t$ follow by the same arguments as in the second and third step in the proof of Lemma \ref{Lemma_invm_ex_org}. 
	
	By Doob's theorem \cite[Theorem 4.2.1(ii), Chapter 4, p. 43]{da_prato_zabczyk_1996} together with \cite[Theorem 3.2.6, Chapter 3, p. 28]{da_prato_zabczyk_1996} and \cite[(5.6) Theorem, Chapter VIII, p. 214]{Bass} there exists a unique ergodic invariant measure with smooth Lebesgue-density for the Markov semigroup $P_t$ associated with system \eqref{Eq_3DHopf_polar_base_aug}.
	
	Using the same argument as in the fourth step of Lemma \ref{Lemma_invm_ex_org}, we can show a bound analogous to \eqref{Eq_exp_conv_1} for the augmented original system \eqref{Eq_3DHopf_polar_base_aug}.
\end{proof}

\bl \label{Lemma_invm_exuniq_red_stoch}
	The Markov semigroup $Q_t$ associated with the reduced system \eqref{Eq_2DHopf_polar_base_aug} has a unique ergodic invariant measure $\bar{\nu}$ with smooth Lebesgue-density on $\hat{\mathcal{M}}\times \mathbb{R}$. Furthermore, a bound analogous to \eqref{Eq_exp_conv_1} holds for $Q_t$.
\el
\begin{proof}
	{\bf Step 1: Existence of invariant measures.} In the first step we prove the existence of at least one invariant measure for the Markov semigroup $Q_t$ associated with the reduced system \eqref{Eq_2DHopf_polar_base_aug}. We define
	\bea \label{Eq_V_lyap_stoch_2}
	V:\,& \hat{\mathcal{M}} \times \mathbb{R} \longrightarrow \mathbb{R}_+,\\
	&(r,\theta, M) \mapsto  \alpha_1 r^2+ \alpha_2 \theta + \alpha_3 M^2 +1 ,
	\eea
	where 
	\bes
	\alpha_1, \alpha_2 = 1 \,\,\, \text{and}\,\,\, \alpha_3 = \frac{\gamma \epsilon}{c_{\tau}}. 
	\ees
	Proceeding as in the proof of Lemmas \eqref{Lemma_invm_ex_org} and \eqref{Lemma_invm_ex_org_aug} yields the existence of constants $\rho>0$ and $C>0$ such that
	\bes
	\mathbb{E}_{\mathbb{P}}[V(r_t,\theta_t,M_t)] \leq e^{-\rho t} V(r_0,\theta_0,M_0) + C,\,\,\, \text{for all }t \geq 0,
	\ees
	which implies the existence of at least one invariant measure for $Q_t$, due to \cite[Proposition 7.10, p. 99]{daprato_inf_dim}.\\
	{\bf Step 2: Strong Feller property and irreducibility.} The strong Feller property and irreducibility for $Q_t$ can be shown as in the second and third step in the proof of Lemma \ref{Lemma_invm_ex_org}. 
	
	Lastly, Doob's theorem \cite[Theorem 4.2.1(ii), Chapter 4, p. 43]{da_prato_zabczyk_1996} together with \cite[Theorem 3.2.6, Chapter 3, p. 28]{da_prato_zabczyk_1996} and \cite[(5.6) Theorem, Chapter VIII, p. 214]{Bass} ensure the existence of a unique ergodic invariant measure for $Q_t$ with smooth Lebesgue-density on $\hat{\mathcal{M}}\times \mathbb{R}$.
	
	A bound analogous to \eqref{Eq_exp_conv_1} can be shown to hold for $Q_t$ using the same argument as in the fourth step of Lemma \ref{Lemma_invm_ex_org}. 
\end{proof}

\section{Auxiliary inequalities for Section \ref{Sec_slow_pm}}
\label{Appendix_aux_ineq}
\bl \label{Lemma_aux_ineq_1}
	Let $T>0$. Consider the unique strong solution $(\hat{R}_t^{\hat{R}_0})_{t \in [0,T]}$ of the reduced system \eqref{Eq_2DHopf_polar_base} to an initial value $\hat{R}_0 = (\hat{r}_0, \hat{\theta}_0)$ in $\hat{\mathcal{M}}$. Given the radial component $(\hat{r}_t)_{t\in [0,T]}$ of $(\hat{R}_t^{\hat{R}_0})_{t \in [0,T]}$ we have the transformed  system \eqref{Eq_3DHopf_polar_base_transf}. Let the stochastic process $(\tilde{R}_t^{R_0})_{t\in [0,T]}$ with 
	\bes
	\tilde{R}_t^{R_0} = (\tilde{r}_t, \tilde{\theta}_t, \tilde{z}_t), \,\,\,t \in [0,T],
	\ees
	be the unique strong solution of the transformed  system \eqref{Eq_3DHopf_polar_base_transf} to an initial value $R_0$ in $\mathcal{M}$, where $R_0= (r_0,\theta_0,z_0)$. In addition, let $\tilde{\mathbb{P}} = \tilde{\mathbb{P}}_T$ be the probability measure given by \eqref{Eq_P_tilde_girs} in Lemma \ref{Lemma_girsanov}. 
  Then, we have for all $t$ in $[0,T]$ that 
  \be  \label{Eq_lemma_aux_ineq_1}
  \mathbb{E}_{\tilde{\mathbb{P}}} \left[|\hat{r}_t - \tilde{r}_t|^2\right] \leq e^{-qt}|\hat{r}_0 - r_0|^2 + \frac{\gamma^2}{q} \int_0^t e^{-q(t-s)} \mathbb{E}_{\tilde{\mathbb{P}}}\left[ \tilde{r}^2_s(\tilde{z}_s-\tilde{r}^2_s)^2\right]\d s.
  \ee
\el
\begin{proof}
	Let $t$ be in $[0,T]$. 	For brevity of notation we introduce the function
	\bea \label{Eq_F_gw_lemma}
	F: \, &\mathbb{R}_+ \longrightarrow \mathbb{R},\\
	&r \mapsto  - \gamma r^3 + \frac{\sigma^2}{2 r}.
	\eea 
	The radial component $(\hat{r}_t)_{t\in[0,T]}$ of the reduced system \eqref{Eq_2DHopf_polar_base} satisfies
	\be \label{Eq_rhat_lemma_gronwall}
	\hat{r}_t = \hat{r}_0 + \int_0^t \left( \lambda \hat{r}_s + F(\hat{r}_s)\right)\d s  + \int_0^t \sigma \d W^r_s,\,\,\, t \in [0,T].
	\ee	
	Recall the definition of $g$ in \eqref{Eq_g}. The radial part $(\tilde{r}_t)_{t\in[0,T]}$ of the transformed system \eqref{Eq_3DHopf_polar_base_transf} satisfies
	\bea \label{Eq_r3_lemma_gronwall}
	\tilde{r}_t = r_0 + \int_0^t \left( \lambda \tilde{r}_s + F(\tilde{r}_s) \right) \d s &- \int_0^t \gamma \tilde{r}_s(\tilde{z}_s-\tilde{r}^2_s) \d s \\
	&+ \int_0^t \sigma g(\tilde{r}_s, \hat{r}_s)\d s	+ \int_0^t \sigma \d W^r_s, \,\,\, t \in [0,T].
	\eea
	Note that, by assumption, the stochastic process $(\hat{r}_t)_{t\in[0,T]}$ in the RHS of Eq.~\eqref{Eq_r3_lemma_gronwall} satisfies Eq.~\eqref{Eq_rhat_lemma_gronwall}.
	Since Eqns.~\eqref{Eq_rhat_lemma_gronwall} and \eqref{Eq_r3_lemma_gronwall} are driven by the same Brownian motion $W^r_t$, the difference between $\hat{r}_t$ and $\tilde{r}_t$, $t$ in $[0,T]$, satisfies
	\beas
	\hat{r}_t - \tilde{r}_t = (\hat{r}_0 - r_0) &+ \int_0^t \left(-q(\hat{r}_s - \tilde{r}_s) + F(\hat{r}_s) - F(\tilde{r}_s) \right) \d s \\
	&+ \int_0^t \gamma \tilde{r}_s(\tilde{z}_s-\tilde{r}^2_s) \d s.
	\eeas
	We also have $\mathbb{P}$-a.s. for all $t$ in $[0,T]$ that
	\bea \label{Eq_lemma_grW^2}
	|\hat{r}_t - \tilde{r}_t|^2 = |\hat{r}_0 - r_0|^2 &+ \int_0^t \left( -2q|\hat{r}_s - \tilde{r}_s|^2 + 2\left(F(\hat{r}_s) - F(\tilde{r}_s)\right)\left(\hat{r}_s - \tilde{r}_s\right) \right)\d s \\
	&+  \int_0^t 2\gamma\tilde{r}_s(\tilde{z}_s-\tilde{r}^2_s)( \hat{r}_s - \tilde{r}_s)\d s.
	\eea
	According to Lemma \ref{Lemma_girsanov}, the probability measures $\mathbb{P}$ and $\tilde{\mathbb{P}}$ are equivalent, i.e. they share the same null sets on $(\Omega, \mathcal{B}_T)$. Therefore, identity \eqref{Eq_lemma_grW^2} also holds $\tilde{\mathbb{P}}$-a.s. for all $t$ in $[0,T]$.
	Differentiating both sides of \eqref{Eq_lemma_grW^2} results in
	\bea \label{Eq_lemma_grW^3}
	\frac{\d}{\d t}|\hat{r}_t - \tilde{r}_t|^2 
	= -2q
	|\hat{r}_t - \tilde{r}_t|^2 + 2(F(\hat{r}_t) - F(\tilde{r}_t))(\hat{r}_t - \tilde{r}_t ) + 2\gamma \tilde{r}_t(\tilde{z}_t-\tilde{r}^2_t) (\hat{r}_t - \tilde{r}_t) .
	\eea
	Since the function $F$ in \eqref{Eq_F_gw_lemma} is decreasing on $\mathbb{R}_+$, it follows for all $t$ in $[0,T]$ that 
	\bes
	(F(\hat{r}_t) - F(\tilde{r}_t))(\hat{r}_t - \tilde{r}_t ) \leq 0,\,\,\, \mathbb{P} \text{  and  } \tilde{\mathbb{P}} \text{  almost surely.}
	\ees
	Thus, $\W$ and $\tilde{\W}$-a.s. for all $t$ in $[0,T]$ we have that
	\be \label{Eq_lemma_grw_4}
	\frac{\d}{\d t}|\hat{r}_t - \tilde{r}_t|^2 
	\leq -2q
	|\hat{r}_t - \tilde{r}_t|^2  + 2\gamma \tilde{r}_t(\tilde{z}_t-\tilde{r}^2_t) (\hat{r}_t - \tilde{r}_t) .
	\ee
	Applying Young's inequality to the last term in the RHS of \eqref{Eq_lemma_grw_4} yields
	\be \label{Eq_lemma_grw_5}
	\frac{\d}{\d t}|\hat{r}_t - \tilde{r}_t|^2  \leq  -q|\hat{r}_t - \tilde{r}_t|^2  + \frac{\gamma^2}{q}  \tilde{r}^2_t\left(\tilde{z}_t-\tilde{r}^2_t\right)^2,\,\,\,\text{for all}\,\,\,t \in [0,T],\,\,\,\W\,\,\,\text{and}\,\,\,\tilde{\W}-a.s.
	\ee
	Due to Gronwall's inequality applied to \eqref{Eq_lemma_grw_5}, we arrive at
	\be \label{Eq_lemma_grw_7}
	|\hat{r}_t - \tilde{r}_t|^2 \leq  e^{-qt} |\hat{r}_0 - r_0|^2 + \frac{\gamma^2}{q} \int_0^t e^{-q(t-s)} \tilde{r}^2_s\left(\tilde{z}_s-\tilde{r}^2_s\right)^2 \d s,\,\,\,t \in [0,T],\,\,\,\W\,\,\,\text{and}\,\,\,\tilde{\W}-a.s.
	\ee
	Taking expectations with respect to $\tilde{\W}$ on both sides of inequality \eqref{Eq_lemma_grw_7} and applying Fubini's theorem finally yields
	\be \label{Eq_lemma_grw_6}
	\E_{\tilde{\W}}\left[|\hat{r}_t - \tilde{r}_t|^2\right] \leq  e^{-qt} |\hat{r}_0 - r_0|^2 + \frac{\gamma^2}{q} \int_0^t e^{-q(t-s)} \E_{\tilde{\W}}\left[\tilde{r}^2_s\left(\tilde{z}_s-\tilde{r}^2_s\right)^2\right] \d s,\,\,\,t \in [0,T],
	\ee
	which is inequality \eqref{Eq_lemma_aux_ineq_1} and completes the proof.
\end{proof}

\bl \label{Lemma_aux_ineq_2}
	Let $T>0$. 
	We consider the unique strong solution $(\hat{R}_t^{\hat{R}_0})_{t\in [0,T]}$ of the reduced system \eqref{Eq_2DHopf_polar_base} for an initial value $\hat{R}_0 = (\hat{r}_0, \hat{\theta}_0)$ in $\hat{\mathcal{M}}$. Given the radial part $(\hat{r}_t)_{t\in[0,T]}$ of $(\hat{R}_t^{\hat{R}_0})_{t\in [0,T]}$ we have the transformed  system \eqref{Eq_3DHopf_polar_base_transf}.
	The stochastic process $(\tilde{R}_t^{R_0})_{t\in [0,T]}$ with 
	\bes
	\tilde{R}_t^{R_0} = (\tilde{r}_t, \tilde{\theta}_t, \tilde{z}_t),\,\,\, t \in [0,T],
	\ees	
	is the unique strong solution to system \eqref{Eq_3DHopf_polar_base_transf} for an initial value $R_0 = (r_0, \theta_0, z_0)$ in $\mathcal{M}$. Similarly, the unique strong solution for the original system \eqref{Eq_3DHopf_polar_base} to the initial value $R_0$ shall be denoted by $(R_t^{R_0})_{t\in [0,T]}$ with
	\bes
	R_t^{R_0} = (r_t, \theta_t, z_t),\,\,\, t \in [0,T].
	\ees
	Let the probability measure $\tilde{\mathbb{P}} = \tilde{\mathbb{P}}_{T}$ with density $\mathfrak{D}_{T}$ with respect to $\mathbb{P}$ be defined as in \eqref{Eq_D_girs} of Lemma \ref{Lemma_girsanov}.
	Recall the definitions of $r_*$ in \eqref{Eq_r_*} and the space of test functions $\mathcal{F}$ in \eqref{Eq_test_fcts_1}.
	Then, for all $\varphi$ in $\mathcal{F}$ and $t$ in $[0,T]$ the following inequality holds
	\bea \label{Eq_aux_ineq_new_2}
	\left| \mathbb{E}_{\mathbb{P}}[\varphi(\hat{r}_t)] - \mathbb{E}_{\tilde{\mathbb{P}}} [\varphi(\tilde{r}_t)]\right| \leq 2 r_* &+ \mathbb{E}_{\tilde{\mathbb{P}}} [\left|\varphi(\hat{r}_t) - \varphi(\tilde{r}_t) \right|] \\
	&+\mathbb{E}_{\mathbb{P}} [|\varphi(\hat{r}_t)| \mathds{1}_{\{\hat{r}_t > r_*\}}] \\
	&+ \mathbb{E}_{\mathbb{P}} [|\varphi(r_t)| \mathds{1}_{\{r_t > r_*\}}].
	\eea
\el
\begin{proof}
	Let $t$ be in $[0,T]$. We define the following events in $(\Omega, \mathcal{B}_T)$,
	\beas
	A_1(t) &= \{ \tilde{r}_t, \hat{r}_t \leq r_*\},  \\
	A_2(t) &= \{ \tilde{r}_t \leq r_*, \hat{r}_t >r_*\},  \\
	A_3(t) &= \{ \tilde{r}_t > r_*, \hat{r}_t \leq r_*\},  \\
	A_4(t) &= \{ \tilde{r}_t, \hat{r}_t > r_*\}.  
	\eeas
	It follows that the sets $A_i(t)$, $i=1,\dots,4$, are disjoint and $\Omega = \bigcup_{i=1}^4 A_{i}(t)$.
	Therefore,
	\bea \label{Lemma_aux_ineq_2_10}
	\left| \mathbb{E}_{\mathbb{P}}[\varphi(\hat{r}_t)] - \mathbb{E}_{\tilde{\mathbb{P}}} [\varphi(\tilde{r}_t)]\right| &= \left| \mathbb{E}_{\mathbb{P}}[\varphi(\hat{r}_t) - \varphi(\tilde{r}_t)\mathfrak{D}_{T}]\right|  \\
	& = \left| \sum_{i = 1}^4 \mathbb{E}_{\mathbb{P}} \left[(\varphi(\hat{r}_t) - \varphi(\tilde{r}_t)\mathfrak{D}_{T}) \mathds{1}_{A_i(t)} \right]\right|. 
	\eea
	For $i = 1$ we obtain 
	\bea \label{Lemma_aux_ineq_2_1}
	\mathbb{E}_{\mathbb{P}} \left[(\varphi(\hat{r}_t) - \varphi(\tilde{r}_t)\mathfrak{D}_{T}) \mathds{1}_{A_1(t)}\right] &= \mathbb{E}_{\mathbb{P}} \left[(\varphi(\hat{r}_t) - \varphi(\hat{r}_t) \mathfrak{D}_{T} + \varphi(\hat{r}_t) \mathfrak{D}_{T} - \varphi(\tilde{r}_t)\mathfrak{D}_{T}) \mathds{1}_{A_1(t)} \right]  \\
	& = \mathbb{E}_{\mathbb{P}} [(1-\mathfrak{D}_{T}) \varphi(\hat{r}_t) \mathds{1}_{A_1(t)}] + \mathbb{E}_{\mathbb{P}}[(\varphi(\hat{r}_t) - \varphi(\tilde{r}_t)) \mathfrak{D}_{T} \mathds{1}_{A_1(t)}] \\
	&= \mathbb{E}_{\mathbb{P}} [(1-\mathfrak{D}_{T}) \varphi(\hat{r}_t) \mathds{1}_{A_1(t)}] + \mathbb{E}_{\tilde{\mathbb{P}}}[(\varphi(\hat{r}_t) - \varphi(\tilde{r}_t)) \mathds{1}_{A_1(t)}]. 
	\eea	
	Let $i = 2$. It is immediate that
	\bea \label{Lemma_aux_ineq_2_2}
	\mathbb{E}_{\mathbb{P}} \left[(\varphi(\hat{r}_t) - \varphi(\tilde{r}_t)\mathfrak{D}_{T}) \mathds{1}_{A_2(t)} \right] = - \mathbb{E}_{\mathbb{P}}[\varphi(\tilde{r}_t)\mathfrak{D}_T \mathds{1}_{A_2(t)}] +
	\mathbb{E}_{\mathbb{P}}[\varphi(\hat{r}_t) \mathds{1}_{A_2(t)}].
	\eea
	Let $i = 3$. Then,
	\bea \label{Lemma_aux_ineq_2_3}
	\mathbb{E}_{\mathbb{P}} \left[(\varphi(\hat{r}_t) - \varphi(\tilde{r}_t)\mathfrak{D}_{T}) \mathds{1}_{A_3(t)} \right] &= \mathbb{E}_{\mathbb{P}}[(1 - \mathfrak{D}_{T})\varphi(\hat{r}_t) \mathds{1}_{A_3(t)}] +\mathbb{E}_{\mathbb{P}}[(\varphi(\hat{r}_t) - \varphi(\tilde{r}_t))\mathfrak{D}_{T} \mathds{1}_{A_3(t)}] \\
	&=\mathbb{E}_{\mathbb{P}}[(1 - \mathfrak{D}_{T})\varphi(\hat{r}_t) \mathds{1}_{A_3(t)}] +\mathbb{E}_{\tilde{\mathbb{P}}}[(\varphi(\hat{r}_t) - \varphi(\tilde{r}_t)) \mathds{1}_{A_3(t)}].
	\eea
	For $i = 4$ it holds that
	\be \label{Lemma_aux_ineq_2_4}
	\mathbb{E}_{\mathbb{P}} \left[(\varphi(\hat{r}_t) - \varphi(\tilde{r}_t)\mathfrak{D}_{T}) \mathds{1}_{A_4(t)} \right]
	= \mathbb{E}_{\mathbb{P}}[\varphi(\hat{r}_t) \mathds{1}_{A_4(t)}] - \mathbb{E}_{\mathbb{P}}[\varphi(\tilde{r}_t)\mathfrak{D}_{T}\mathds{1}_{A_4(t)}]. \\
	\ee
	Inserting  Eqns.~\eqref{Lemma_aux_ineq_2_1}, \eqref{Lemma_aux_ineq_2_2}, \eqref{Lemma_aux_ineq_2_3} and \eqref{Lemma_aux_ineq_2_4} into Eq.~\eqref{Lemma_aux_ineq_2_10} and applying the triangular inequality yields
	\bea \label{Eq_aux_ineq_new_1}
	\left| \mathbb{E}_{\mathbb{P}}[\varphi(\hat{r}_t)] - \mathbb{E}_{\tilde{\mathbb{P}}} [\varphi(\tilde{r}_t)]\right| \leq 2r_* &+ \mathbb{E}_{\tilde{\mathbb{P}}}[|\varphi(\hat{r}_t) - \varphi(\tilde{r}_t)|] \\
	& +\mathbb{E}_{\mathbb{P}}[|\varphi(\hat{r}_t)|\mathds{1}_{\{\hat{r}_t>r_*\}}] \\
	&+ \mathbb{E}_{\tilde{\mathbb{P}}}[|\varphi(\tilde{r}_t)|\mathds{1}_{\{\tilde{r}_t>r_*\}}].
	\eea
	Inequality \eqref{Eq_aux_ineq_new_2}  follows then from identity \eqref{Eq_trans_prob_girs} in Lemma \ref{Lemma_girsanov} applied to the last term in Eq.~\eqref{Eq_aux_ineq_new_1}. 
\end{proof}
\bl \label{lemma_aux_ineq_new_lemma}
	Under the assumptions of Lemma \ref{Lemma_aux_ineq_2} it holds for all $t$ in $[0,T]$ that
	\bea \label{Eq_new_aux_ineq_1_to_show}
	\left| \mathbb{E}_{\mathbb{P}}[\varphi(\hat{r}_t)] - \mathbb{E}_{\tilde{\mathbb{P}}} [\varphi(\tilde{r}_t)]\right| \leq e^{-\frac{\lambda}{2} t} \hat{r}_0 &+\mathbb{E}_{\tilde{\mathbb{P}}}[|\tilde{r}_t - \hat{r}_t|] + \frac{(q+\lambda)}{\sqrt{\lambda}}\left(\int_0^t e^{-\lambda(t-s)} \mathbb{E}_{\tilde{\mathbb{P}}}[|\hat{r}_s - \tilde{r}_s|^2] \d s\right)^{\frac{1}{2}} \\
	&+\sqrt{2\left(r_{\mathrm{det}}^2 + \frac{\sigma^2}{\lambda}\right)}   +\mathbb{E}_{\mathbb{P}}[\hat{r}_t] ,
	\eea
	where 
	\bes
	r_{\mathrm{det}} = \sqrt{\frac{\lambda}{\gamma}},
	\ees
	denotes the radius of the limit cycle associated with systems \eqref{Eq_3DHopf_polar_base} and \eqref{Eq_2DHopf_polar_base} for $\sigma=0$.
\el
\begin{proof}
Let $t \in [0,T]$. Due to the definitions of $\tilde{\mathbb{P}}$ in \eqref{Eq_P_tilde_girs} and $\mathcal{F}$ in \eqref{Eq_test_fcts_1} above, we have for all $\varphi$ in $\mathcal{F}$ that
		\bea \label{Eq_new_aux_ineq_1}
		\left| \mathbb{E}_{\mathbb{P}}[\varphi(\hat{r}_t)] - \mathbb{E}_{\tilde{\mathbb{P}}} [\varphi(\tilde{r}_t)]\right| &= \left|\mathbb{E}_{\mathbb{P}}[\varphi(\hat{r}_t) - \varphi(\hat{r}_t)\mathfrak{D}_T + \varphi(\hat{r}_t)\mathfrak{D}_T - \varphi(\tilde{r}_t)\mathfrak{D}_T]\right| \\
		&\leq \mathbb{E}_{\mathbb{P}}[\hat{r}_t] + \mathbb{E}_{\tilde{\mathbb{P}}}[\hat{r}_t] + \mathbb{E}_{\tilde{\mathbb{P}}}[\left|\hat{r}_t-\tilde{r}_t\right|].
		\eea
		In order to derive inequality \eqref{Eq_new_aux_ineq_1_to_show} we need to estimate the expected value of $\hat{r}_t$ for the probability $\tilde{\mathbb{P}}$.
		Recall that, by virtue of Lemma \ref{Lemma_girsanov}, we have that the process $(\tilde{W}^r_t)_{t\in[0,T]}$ given by
		\be \label{Eq_hilf_2}
		\tilde{W}_t^r = \int_0^t g(\tilde{r}_s, \hat{r}_s)\d s + W_t^r, \,\,\,t \in [0,T],
		\ee
		with $g$ defined in \eqref{Eq_g}, is a real-valued Brownian motion on $(\Omega, \mathcal{B}_T, \{\mathcal{B}_t\}_{t\in[0,T]}, \tilde{\mathbb{P}})$. Furthermore, the radial component $(\hat{r}_t)_{t\in [0,T]}$ of the reduced system \eqref{Eq_2DHopf_polar_base} satisfies the integral equation
		\be \label{Eq_hilf_6}
		\hat{r}_t = \hat{r}_0 + \int_0^t \lambda \hat{r}_s - \gamma \hat{r}_s^3 + \frac{\sigma^2}{2\hat{r}_s} \d s + \int_0^t \sigma \d W^r_s,\,\,\,\W-a.s.\,\,\,\text{for all}\,\,\,t\in[0,T].
		\ee
		In addition we have, due to \eqref{Eq_hilf_2}, for all $t\in [0,T]$ that
		\be \label{Eq_hilf_7}
		\hat{r}_t = \hat{r}_0 + \int_0^t \lambda \hat{r}_s - \gamma \hat{r}_s^3 + \frac{\sigma^2}{2\hat{r}_s} \d s - \int_0^t \sigma g(\tilde{r}_s,\hat{r}_s)\d s + \int_0^t \sigma \d \tilde{W}^r_s,\,\,\,\W\,\,\,\text{and}\,\,\,\tilde{\W}-a.s.
		\ee
		Note that the process $(\hat{r}_t)_{t\in[0,T]}$ in \eqref{Eq_hilf_7} is a semi-martingale (e.g.~\cite[Chapter II, p.~66]{watanabe_ikeda}) on $(\Omega, \mathcal{B}_T, \{\mathcal{B}_t\}_{t\in[0,T]}, \tilde{\mathbb{P}})$. Applying It\^{o}'s formula to $(\hat{r}_t^2)_{t\in[0,T]}$ yields the following integral equation for all $t$ in $[0,T]$
		\be \label{Eq_hilf_1}
		\hat{r}_t^2 = \hat{r}_0^2 + \int_0^t \left(2\lambda \hat{r}_s^2 - 2\gamma \hat{r}^4_s + 2\sigma^2 - 2\sigma \hat{r}_s g(\tilde{r}_s,\hat{r}_s)\right)\d s + 2\sigma \int_0^t \hat{r}_s \d \tilde{W}^r_s,\,\,\,\W\,\,\,\text{and}\,\,\,\tilde{\W}-a.s.
		\ee
		We remark that the process 
		\bes
		\left(\int_0^t \hat{r}_s \d \tilde{W}^r_s\right)_{t\in[0,T]}\,\,\,\text{is a martingale on}\,\,\, (\Omega, \mathcal{B}_T, \{\mathcal{B}_t\}_{t\in[0,T]}, \tilde{\mathbb{P}}).
		\ees
		According to \cite[Proposition 2.2, p.~55-56]{watanabe_ikeda} the proof of this statement consists in verifying that 
		\bes
		\E_{\tilde{\W}}\left[\int_0^T \betrag{\hat{r}_t}^2\d t\right] < \infty.
		\ees
		This can be shown from Eq.~\eqref{Eq_hilf_1} and a localization argument similar to that employed in the proof of Lemma \ref{Lemma_moment_bounds}.	
		Next, by taking expectations {\mkr for the probability} $\tilde{\mathbb{P}}$ in \eqref{Eq_hilf_1}, applying Fubini's theorem and differentiating on both sides, we obtain
		\bea \label{Eq_hilf_4}
		\frac{\d}{\d t} \mathbb{E}_{\tilde{\mathbb{P}}}[\hat{r}_t^2] = 2\lambda \mathbb{E}_{\tilde{\mathbb{P}}}[\hat{r}_t^2] - 2\gamma \mathbb{E}_{\tilde{\mathbb{P}}}[\hat{r}_t^4] + 2\sigma^2 - 2\sigma \mathbb{E}_{\tilde{\mathbb{P}}}[\hat{r}_t g(\tilde{r}_t,\hat{r}_t)],\,\,\,t\in[0,T].
		\eea
		A straight-forward computation shows that
		\be \label{Eq_hilf_3}
		2\lambda r^2 - 2\gamma r^4 \leq -2\lambda r^2 + \frac{2\lambda^2}{\gamma},\,\,\,
		\text{for all}\,\,\,r\geq0.
		\ee
		Inserting inequality \eqref{Eq_hilf_3} into \eqref{Eq_hilf_4} and applying Young's inequality to the last term in the RHS of Eq.~\eqref{Eq_hilf_4} results in 
		\be \label{Eq_hilf_5}
		\frac{\d}{\d t} \mathbb{E}_{\tilde{\mathbb{P}}}[\hat{r}_t^2] \leq -\lambda \mathbb{E}_{\tilde{\mathbb{P}}}[\hat{r}_t^2] + 2\left(\frac{\lambda^2}{\gamma} + \sigma^2\right) + \frac{(q+\lambda)^2}{\lambda} \mathbb{E}_{\tilde{\mathbb{P}}}[|\hat{r}_t - \tilde{ r}_t|^2],\,\,\,t\in[0,T].
		\ee
		Using Gronwall's inequality in \eqref{Eq_hilf_5} yields for all $t$ in $[0,T]$ that
		\bes
		\mathbb{E}_{\tilde{\mathbb{P}}}[\hat{r}_t^2] \leq e^{-\lambda t} \hat{r}_0^2 + 2\left(r_{\mathrm{det}}^2 + \frac{\sigma^2}{\lambda}\right) + \frac{(q+\lambda)^2}{\lambda}\int_0^t e^{-\lambda(t-s)} \mathbb{E}_{\tilde{\mathbb{P}}}[|\hat{r}_s - \tilde{r}_s|^2] \d s.
		\ees
		Furthermore, by means of the Cauchy-Schwarz inequality it holds for all $t$ in $[0,T]$ that
		\bea \label{Eq_lemma_aux_final_1}
		\mathbb{E}_{\tilde{\mathbb{P}}}[\hat{r}_t]&\leq \sqrt{\mathbb{E}_{\tilde{\mathbb{P}}}[\hat{r}_t^2]}
		\leq e^{-\frac{\lambda}{2} t} \hat{r}_0 + \sqrt{2\left(r_{\mathrm{det}}^2 + \frac{\sigma^2}{\lambda}\right)} + \frac{(q+\lambda)}{\sqrt{\lambda}}\left(\int_0^t e^{-\lambda(t-s)} \mathbb{E}_{\tilde{\mathbb{P}}}[|\hat{r}_s - \tilde{r}_s|^2] \d s\right)^{\frac{1}{2}}.	 
		\eea
Inserting inequality \eqref{Eq_lemma_aux_final_1} into \eqref{Eq_new_aux_ineq_1}
			yields \eqref{Eq_new_aux_ineq_1_to_show} and completes the proof.
\end{proof}

\section{Auxiliary inequalities for Section \ref{Sec_stoch_pm}}
\label{Appendix_aux_ineq_2}
The next lemma is analogous to Lemma \ref{Lemma_aux_ineq_1} for the stochastic parameterizing manifold case discussed in Section \ref{Sec_stoch_pm}.

\bl \label{Lemma_aux_ineq_3}
	Let $T>0$. 
	Consider the unique strong solution $(\hat{R}_t^{\hat{R}_0})_{t \in [0,T]}$ of the reduced system \eqref{Eq_2DHopf_polar_base_aug} to an initial value $\hat{R}_0 = (\hat{r}_0, \hat{\theta}_0,M_0)$ in $\hat{\mathcal{M}}\times \mathbb{R}$. Given the radial component $(\hat{r}_t)_{t\in [0,T]}$ and the process $(M_t)_{t\in[0,T]}$ of $(\hat{R}_t^{\hat{R}_0})_{t \in [0,T]}$ we obtain the transformed  system \eqref{Eq_3DHopf_polar_base_transf_aug}.
	The stochastic process $(\tilde{R}_t^{R_0})_{t\in [0,T]}$ with 
	\bes 
	\tilde{R}_t^{R_0} = (\tilde{r}_t, \tilde{\theta}_t, \tilde{z}_t), \,\,\, t \in [0,T],
	\ees
	is the unique strong solution of system \eqref{Eq_3DHopf_polar_base_transf_aug} to the initial value $R_0 = (r_0,\theta_0,z_0)$ in $\mathcal{M}$. Let $\tilde{\mathbb{P}} = \tilde{\mathbb{P}}_T$ be the probability measure given by identity \eqref{Eq_P_tilde_girs_2} of Lemma \ref{Lemma_girsanov_stoch_mnf}.
	Then, we have for all $m$ in $\mathbb{R}$ and $t$ in $[0,T]$ that 
	\bea  \label{Eq_lemma_aux_ineq_3}
	\mathbb{E}_{\tilde{\mathbb{P}}} \left[|\hat{r}_t - \tilde{r}_t|^2\right] \leq e^{-qt}|\hat{r}_0 - r_0|^2 &+ \frac{\gamma^2}{q} \int_0^t e^{-q(t-s)} \mathbb{E}_{\tilde{\mathbb{P}}}\left[ \tilde{r}^2_s(\tilde{z}_s-h_{\tau}(m,\tilde{r}_s))^2\right]\d s\\
	&+ 2m \gamma \int_0^t e^{-q(t-s)}\mathbb{E}_{\tilde{\mathbb{P}}}[\tilde{r}_s(\hat{r}_s - \tilde{r}_s)]\d s,
	\eea
	for $h_{\tau}$ defined in \eqref{Eq_stoch_mnf} with $\tau$ taking values in $(0,\infty)$.
\el
\begin{proof}
	The proof of this lemma for Eqns.~\eqref{Eq_2DHopf_polar_base_aug} and \eqref{Eq_3DHopf_polar_base_transf_aug} is analogous to that of Lemma \ref{Lemma_aux_ineq_1} for Eqns.~\eqref{Eq_2DHopf_polar_base} and \eqref{Eq_3DHopf_polar_base_transf}, in which we replace  the usage of Lemma \ref{Lemma_girsanov} by its counterpart for SPMs (Lemma \ref{Lemma_girsanov_stoch_mnf}).
\end{proof}

The following lemma is the counterpart, for SPM reduced systems, to Lemma \ref{Lemma_aux_ineq_2}.
\bl \label{Lemma_aux_ineq_4}
Let $T>0$. 
Let $(\hat{R}_t^{\hat{R}_0})_{t\in [0,T]}$ be the unique strong solution of the reduced system \eqref{Eq_2DHopf_polar_base_aug} for an initial value $\hat{R}_0 = (\hat{r}_0, \hat{\theta}_0,M_0)$ in $\hat{\mathcal{M}}\times \mathbb{R}$. Given the radial component $(\hat{r}_t)_{t\in[0,T]}$ and the process $(M_t)_{t\in[0,T]}$ of $(\hat{R}_t^{\hat{R}_0})_{t\in [0,T]}$ we have the transformed  system \eqref{Eq_3DHopf_polar_base_transf_aug} with unique strong solution
$(\tilde{R}_t^{R_0})_{t\in [0,T]}$ of the form 
\bes
\tilde{R}_t^{R_0} = (\tilde{r}_t, \tilde{\theta}_t, \tilde{z}_t),\,\,\, t \in [0,T].
\ees	
The initial value $R_0$ is an element of $\mathcal{M}$. Furthermore, the unique strong solution for the original system \eqref{Eq_3DHopf_polar_base_stoch} to the initial value $R_0$ is denoted by $(R_t^{R_0})_{t\in [0,T]}$ with
\bes
R_t^{R_0} = (r_t, \theta_t, z_t),\,\,\, t \in [0,T].
\ees
Let the probability measure $\tilde{\mathbb{P}} = \tilde{\mathbb{P}}_{T}$ with density $\mathfrak{D}_{T}$ with respect to $\mathbb{P}$ be defined as in \eqref{Eq_P_tilde_girs_2} of Lemma \ref{Lemma_girsanov_stoch_mnf}.
Recall the definitions of $r_*$ in \eqref{Eq_r_*_stoch} and the space of test functions $\mathcal{F}$ in \eqref{Eq_test_fcts_1}.
Then, for all $\varphi$ in $\mathcal{F}$ and $t$ in $[0,T]$ the following inequality holds
\bea \label{Eq_aux_ineq_new_3}
\left| \mathbb{E}_{\mathbb{P}}[\varphi(\hat{r}_t)] - \mathbb{E}_{\tilde{\mathbb{P}}} [\varphi(\tilde{r}_t)]\right| \leq 2 r_* &+ \mathbb{E}_{\tilde{\mathbb{P}}} [\left|\varphi(\hat{r}_t) - \varphi(\tilde{r}_t) \right|] \\
&+\mathbb{E}_{\mathbb{P}} [|\varphi(\hat{r}_t)| \mathds{1}_{\{\hat{r}_t > r_*\}}] \\
&+ \mathbb{E}_{\mathbb{P}} [|\varphi(r_t)| \mathds{1}_{\{r_t > r_*\}}].
\eea
\el
\begin{proof}
The proof  proceeds analogously to that  of Lemma \ref{Lemma_aux_ineq_2}.
\end{proof}

\bl \label{lemma_aux_ineq_new_lemma2}
	Under the assumptions of Lemma \ref{Lemma_aux_ineq_4} it holds for all $t$ in $[0,T]$ that
	\bea \label{Eq_new_aux_ineq_1_to_show2}
	\left| \mathbb{E}_{\mathbb{P}}[\varphi(\hat{r}_t)] - \mathbb{E}_{\tilde{\mathbb{P}}} [\varphi(\tilde{r}_t)]\right| \leq e^{-\frac{\lambda}{2} t} \hat{r}_0 &+\mathbb{E}_{\tilde{\mathbb{P}}}[|\tilde{r}_t - \hat{r}_t|] + \frac{(q+\lambda)}{\sqrt{\lambda}}\left(\int_0^t e^{-\lambda(t-s)} \mathbb{E}_{\tilde{\mathbb{P}}}[|\hat{r}_s - \tilde{r}_s|^2] \d s\right)^{\frac{1}{2}} \\
	&+\sqrt{2\left(\frac{r_{\mathrm{det}}^2}{c_{\tau}} + \frac{\sigma^2}{\lambda}\right)}   +\mathbb{E}_{\mathbb{P}}[\hat{r}_t] ,
 	\eea
	where $q$ is defined in \eqref{Eq_def_q_2} and 
	\bes
	r_{\mathrm{det}} = \sqrt{\frac{\lambda}{\gamma}},
	\ees
	denotes the radius of the limit cycle associated with systems \eqref{Eq_3DHopf_polar_base} and \eqref{Eq_2DHopf_polar_base} for $\sigma=0$.
\el
\begin{proof}
The proof  proceeds analogously to that  of Lemma \ref{lemma_aux_ineq_new_lemma}.
\end{proof}

}

\bibliographystyle{amsalpha}
\bibliography{Slow_PM_rev}

\end{document}